\documentclass[twoside,reqno]{amsart}

\usepackage[dvips,usenames,dvipsnames,svgnames,table]{xcolor}
\usepackage[dvips,bookmarks,colorlinks=true,linkcolor=blue]{hyperref}

\usepackage{euscript}
\usepackage{amssymb}
\usepackage{latexsym}
\usepackage{amsmath,amsfonts}
\usepackage{amsthm,amstext}
\usepackage{amsbsy,amsxtra,amsopn}

\usepackage{gensymb}
\newcommand{\ssim}{\overset{\ensuremath{\raisebox{-1pt}{\tiny$\mathsf{s}$}}}{\sim}}

\newcommand{\compcent}[1]{\vcenter{\hbox{$#1\circ$}}}
\newcommand{\comp}{\mathbin{\mathchoice
{\compcent\scriptscriptstyle}{\compcent\scriptscriptstyle}
{\compcent\scriptscriptstyle}{\compcent\scriptscriptstyle}}}

\newtheorem{theorem}{Theorem}[section]
\newtheorem{proposition}[theorem]{Proposition}
\newtheorem{lemma}[theorem]{Lemma}
\newtheorem{corollary}[theorem]{Corollary}
\newtheorem*{ThA}{Theorem A}
\newtheorem*{ThB}{Theorem B}

\theoremstyle{definition}
\newtheorem{definition}[theorem]{Definition}
\newtheorem{example}[theorem]{Example}

\newtheorem{remark}[theorem]{Remark}

\numberwithin{equation}{section}

\newcommand{\dR}{{\mathbb R}}
\newcommand{\dC}{{\mathbb C}}

\newcommand{\kip}{[\,\cdot\, , \cdot\,]}

\newcommand{\ahip}{\langle \,\cdot\, , \cdot\,\rangle}

\newcommand{\la}{\langle} \newcommand{\ra}{\rangle}

\newcommand{\cH}{{\mathcal H}}

\newcommand{\cK}{{\mathcal K}}
\newcommand{\cL}{{\mathcal L}}

\newcommand{\cS}{{\mathcal S}}

\def\wt#1{{{\widetilde #1} }}
\newcommand{\wh}[1]{\widehat{#1}}

 \DeclareMathOperator{\sgn}{sgn}
 
 \DeclareMathOperator{\ran}{ran}
 \DeclareMathOperator{\dom}{dom}
 \DeclareMathOperator{\re}{Re}
 \DeclareMathOperator{\im}{Im}
 
 \DeclareMathOperator{\loc}{loc}
\DeclareMathOperator{\lspan}{span}

\allowdisplaybreaks

\makeatletter
\def\l@section{\@tocline{2}{1pt plus1pt}{1.5pc}{0pt}{}}
\def\l@subsection{\@tocline{1}{1pt plus1pt}{3.2pc}{0pt}{}}
\makeatother

\makeatletter
\renewcommand\p@enumii{}
\makeatother

\begin{document}

\title[Indefinite Sturm-Liouville operators]%
{Indefinite Sturm-Liouville operators\\ in polar form}

\author{Branko \'{C}urgus}
\address{Department of Mathematics, Western Washington University,
Bellingham, Washington 98226, USA}
 \email{curgus@wwu.edu}

\author{Volodymyr Derkach}
\address{
	Institut f\"ur  Mathematik,
	Technische Universit\"{a}t Ilmenau,
	Germany}
\address{
	Vasyl Stus Donetsk National University,
	Vinnytsya, Ukraine}

\email{volodymyr.derkach@tu-ilmenau.de}

 \author{Carsten Trunk}
\address{Institut f\"ur  Mathematik,
Technische Universit\"{a}t Ilmenau,
Postfach 100565, D-98684 Ilmenau,
Germany}
 \email{carsten.trunk@tu-ilmenau.de}

\begin{abstract}
We consider the indefinite Sturm-Liouville differential expression
\[
\mathfrak{a}(f) := - \frac{1}{w}\left( \frac{1}{r} f' \right)',
\]
where $\mathfrak{a}$ is defined on a finite or infinite open interval $I$ with $0\in I$ and
the coefficients $r$ and $w$ are locally summable and such that $r(x)$ and $(\sgn x) w(x)$ are positive a.e.\ on $I$. With the differential expression $\mathfrak{a}$ we associate a nonnegative self-adjoint operator $A$ in the Krein space $L^2_w(I)$ which is viewed as a coupling of symmetric operators in Hilbert spaces related to the intersections of $I$ with the positive and the negative semi-axis. For the operator $A$ we derive conditions in terms of the coefficients $w$ and $r$ for the existence of a Riesz basis consisting of generalized eigenfunctions of $A$ and for the similarity of $A$ to a self-adjoint operator in a Hilbert space $L^2_{|w|}(I)$. These results are obtained as consequences of abstract results about the  regularity of critical points of nonnegative self-adjoint operators in Krein spaces which are couplings of two symmetric operators acting in Hilbert spaces.
\end{abstract}

\subjclass[2000]{Primary 47B50, 34B24; Secondary 46C20, 47B25, 34L10}



\keywords{self-adjoint extension, symmetric operator, Krein space,
Riesz basis, coupling of operators, boundary triple, Weyl function, regular critical point}

\maketitle

 \tableofcontents

\section{Introduction}

Let $I = (b_-,b_+)$ be a finite or infinite interval such that $-\infty \leq b_- < 0 < b_+ \leq + \infty$. We consider the indefinite Sturm-Liouville differential expression $\mathfrak{a}$ on $I$ that is given in polar form
\begin{equation} \label{eq2ndq69A}
\bigl(\mathfrak{a}f \bigr)(x):= -  \frac{1}{w(x)} \frac{d}{dx} \left( \frac{1}{r(x)} \frac{d}{dx} f(x) \right),
\end{equation}
where the coefficients $r$ and $w$ are real functions on $I$ satisfying the conditions
\begin{equation}\label{eq:rw}
 r, w \in L^1_{\loc}(I) \quad \text{and} \quad  r(x) , \ (\sgn x) w(x) > 0 \quad \text{for a.a.} \quad x \in I.
\end{equation}

With the differential expression $\mathfrak{a}$ we associate a closed linear operator $A$ in the weighted Hilbert space $L^2_{|w|}(I)$. The operator $A$ is not self-adjoint in $L^2_{|w|}(I)$ but it is self-adjoint and nonnegative in the Krein space $L^2_{w}(I)$ which coincides with $L^2_{|w|}(I)$ as a normed vector space and has indefinite inner product
\begin{equation*} 
   [f,g]_{w} := \int_I f(x)\overline{g(x)} w(x) dx,
\end{equation*}
see~\cite{DL77} for a similar setting.

We are interested in the following two  properties of the differential operator $A$:
\begin{enumerate}
  \item [(Ri)]
  \textit{Riesz basis property}, that is, the existence of a Riesz basis of the Hilbert space  $L^2_{|w|}(I)$ which consists of eigenfunctions and generalized eigenfunctions of $A$;
  \item [(Si)]
  \textit{Similarity} of $A$ to a self-adjoint operator in the Hilbert space $L^2_{|w|}(I)$, that is, the existence of a bounded and boundedly invertible operator $T$ such that $T^{-1}AT$ is self-adjoint in the Hilbert space  $L^2_{|w|}(I)$.
\end{enumerate}

Our results will be formulated in terms of the functions
\begin{equation}\label{eq:WR}
  W_\pm(x):=\int_0^x w_\pm(\xi)d \xi,\quad
R_\pm(x):=\int_0^x r_\pm(\xi)d\xi,\quad x\in I_\pm,
\end{equation}
where $I_- = (b_-,0)$, $I_+ = (0,b_+)$, $w_{-}$ is the restriction of $-w$ onto $I_-$ and $w_{+}$ is the restriction of $w$ onto $I_+$.

The first result of this kind was given by Beals in~\cite{Bea85}, where the Riesz basis property was proved for the constant function $r= 1$ and for a weight $w$ which behaves as a power at $0$, see Example~\ref{ex:Beals}. This result was generalized to ordinary differential equations of higher order  by \'{C}urgus and Langer \cite{Cur}, \cite{CL89}, and to partial differential equations  by Pyatkov~\cite{Py3}, \'{C}urgus and Najman~\cite{CN94}. The first proof of the existence of a weight $w$, with $r=1$, for which $A$ does not have the Riesz basis property
was given by Volkmer in \cite{Vol} by using Baire category arguments. Explicit examples of such a weight were found by Fleige, \cite{Fl98}, and Abasheeva, Pyatkov \cite{AP}. A full characterization of the Riesz basis property for the operator $A$ was given by Parfenov \cite{Par03} in the case when $w$ is odd and $r=1$. Using Pyatkov's approach via interpolation spaces \cite{Py1,Py3}, Parfenov proved that the Riesz basis property for the operator $A$  holds if and only if the function $W_+$ is positively increasing at $0_+$. Recall,  see \cite[Definition~3.26]{BIKS}, that a nondecreasing positive function $\varphi$ is called {\it positively increasing at} $0_+$ if there exists $\lambda \in (0,1)$ such that $\limsup_{x\downarrow 0}\bigl(\varphi(\lambda x)/\varphi(x)\bigr) < 1$; $\psi$ is {\em positively increasing at} $0_-$ if $x\mapsto -\psi(-x)$ is positively increasing at $0_+$.

In~\cite{Kost13} Kostenko used a different method to characterize the properties (Si) and (Ri) for differential operator $A$ with odd $w$ and even $r$. In particular, it was shown in~\cite{Kost13}  that the Riesz basis property for the operator $A$  holds if and only if the function $W_+\circ R_+^{-1}$ is positively increasing at $0_+$.

One of the main results of this paper is the following theorem in which we give a sufficient condition for the Riesz basis property, in the spirit of Parfenov's and Kostenko's results, but without the assumptions that $w$ is odd and $r$ is even. We also give a new kind of characterization of the Riesz basis property when $W_\pm\circ R_\pm^{-1}$ are slowly varying functions. Recall that a measurable positive function $\varphi$ is said to be {\em slowly varying at $0_+$} if for all $\lambda > 0$ we have $\lim_{x\downarrow 0} \bigl(\varphi(\lambda x)/\varphi(x)\bigr) = 1$; $\psi$ is {\em slowly varying at} $0_-$ if $x\mapsto -\psi(-x)$ is slowly varying at $0_+$, for more about slowly varying functions see Appendix~\ref{appen}.

\begin{ThA} 
Let the differential expression $\mathfrak{a}$ satisfy \eqref{eq:rw} and let $W_\pm$ and $R_\pm$ be the functions defined in \eqref{eq:WR}. Assume that the spectrum of the operator $A$ associated with the differential expression $\mathfrak{a}$ in the Hilbert space $L_{|w|}^2(I)$ is discrete. Then the eigenvalues of $A$ accumulate on both sides of $\infty$ and the following two statements hold.
\begin{enumerate}
\renewcommand*\theenumi{\alph{enumi}}
\renewcommand*\labelenumi{\rm{(\theenumi)}}
\item \label{thmA:Riesz-B}
If either $W_+\circ R_+^{-1}$ is positively increasing at $0_+$ or $W_-\circ R_-^{-1}$ is positively increasing at $0_-$, then the operator $A$ has the Riesz basis property {\rm (Ri)}.

\item \label{thmA:Riesz-A}
If $W_+ \circ R_+^{-1}$ is slowly varying at $0_+$ and $W_- \circ R_-^{-1}$ is slowly varying at $0_-$, then the Riesz basis property {\rm (Ri)} is equivalent to the condition
\begin{equation}\label{eq:FarFromOdd}
\left( 1 + \frac{W_-\bigl(R_-^{-1}(-x)\bigr)}{W_+\bigl(R_+^{-1}(x)\bigr)} \right)^{\!\!-1} = O(1) \quad \text{as}\quad  x\downarrow 0.
\end{equation}
\end{enumerate}
\end{ThA}

The main tool that we use in this paper is Langer's spectral theory of definitizable operators in Krein spaces, see~\cite{La}. Our differential operator $A$ is a nonnegative self-adjoint operator with a nonempty resolvent set in the Krein space $L^2_{w}(I)$. This is a special kind of a definitizable operator that admits a spectral function  $E$ which behaves similarly to the spectral function of a self-adjoint operator in a Hilbert space with a possible exception at its critical points which are contained in the set $\{0,\infty\}$, for details see Section~\ref{Sec2}. A critical point is called regular if $E$ is bounded in a neighbourhood of that point. Otherwise, a critical point is called singular. By $c_s(A) \subseteq \{0,\infty\}$ we denote the set of singular critical points of $A$ and by $c_r(A) \subseteq \{0,\infty\}$ the set of regular critical points of $A$.

In the case of discrete spectrum of the differential operator $A$, the Riesz basis property of $A$ is equivalent to the regularity of the critical point $\infty$, see \cite[Proposition~4.1]{CL89}. This fact and the paper of Beals~\cite{Bea85} were motivation for \cite{CL89, CN95, Fl95, Vol, BC06, BC09, Fl08, Kost11, CFK13} to study definitizability and the regularity of the critical point infinity for differential operators; see also a detailed survey by Fleige~\cite{Fl15}.

The regularity of both critical points of $A$ is equivalent to $A$ being similar to a self-adjoint operator in a Hilbert space. This fact was used by \'{C}urgus and Najman in \cite{CN94} to prove that the operator associated with \eqref{eq2ndq69A} where $w(x) = \sgn(x)$, $r=1$ and $I = \mathbb{R}$ is similar to a self-adjoint operator in the Hilbert space $L^2(\mathbb{R})$. This result was reproved and generalized by Krein space and other methods by several authors, see \cite{CN96, CN98, FN98, Kar00, FS00, KaKost08, KaKoMal09, KT09, Kost13}.


We use the resolvent criterion of Veseli\'{c} \cite{Ves}, \cite{Ako80}, to study regularity of critical points of the operator $A$ in terms of the Weyl functions $m_+$ and $m_-$ of some symmetric operators generated by $\mathfrak{a}$ on intervals $I_+$ and $I_-$. It was shown in~\cite{KaKost08} that the so-called $D_\infty$-\textit{property} (resp. $D_0$-\textit{property})
\begin{equation} \label{Intro:Dinf}
\frac{\max\bigl\{\im m_+(iy),\im m_-(iy)
\bigr\}}{\bigl|m_+(iy)+m_-(-iy)\bigr|}=O(1) \quad
\text{as} \quad y \to +\infty\quad ({\rm resp}.\, y \downarrow 0)
\end{equation}
is necessary for $\infty\not\in c_s(A)$ (resp. $0\not\in c_s(A)$).
In the case when $w$ is odd and $r$ is even the functions $m_+$ and $m_-$ coincide. In this case,  conditions~\eqref{Intro:Dinf}  can be rewritten as
\begin{equation} \label{Intro:ImRem_pm}
\im m_+(iy)=O(\re m_+(iy))\quad
\text{as} \quad y \to +\infty\quad ({\rm resp}.\ y \downarrow 0)
\end{equation}
and are proved in~\cite{Kost13} to be equivalent both to the similarity property (Si) and to the validity of the Hardy, Everitt, Littlewood and Polya (HELP) inequality, see~\cite{Ev72}.
Moreover, in the general case it was proved in~\cite{Kost14} that the $D_\infty$-property together with the $D_0$-property is equivalent to the validity of the so-called Volkmer inequality \cite{Vol}, an indefinite analog of the HELP inequality.

As another main result of our paper, in Theorem~\ref{prop:Semibound2} we prove that the $D_\infty$-property is necessary and sufficient for $\infty\not\in c_s(A)$
provided that the Weyl functions $m_+$ and $m_-$ satisfy the assumption:
\begin{equation} \label{Intro:Rem_pm}
\text{For some $y_0>0$} \quad  \re m_+(iy) \re m_-(iy) > 0 \quad \text{for all} \quad y>y_0.
\end{equation}
The equivalence in Theorem~A(\ref{thmA:Riesz-A}) is obtained by combining the characterization of the regularity of the critical point $\infty$ for the operator $A$ from Theorem~\ref{prop:Semibound2} with the Atkinson-Bennewitz asymptotic formula for the Weyl functions $m_+(iy)$ and $m_-(iy)$ proved in~\cite{Atk85} and~\cite{Ben89} and presented in Theorem~\ref{thm:Asym_m2}.


The questions of similarity of a differential operator to a self-adjoint operator and the existence of a Riesz basis consisting of its eigenfunctions arise in problems of numerical computation of eigenvalues. For example, in \cite[Subsection 4.1.2]{HT} the authors study the differential expression \eqref{eq2ndq69A} with $w(x) = x^3$, $r=1$ and $I = [-1,1]$. To construct an efficient and accurate eigensolver for the associated differential operator it was important that the operator is similar to a self-adjoint operator and that its eigenfunctions form a Riesz basis of the Hilbert space $L^2_{|w|}[-1,1]$.

This paper is organized as follows. In Sections~\ref{Sec2} and~\ref{Sec3} we establish conditions for the regularity of the critical points $0$ and $\infty$ for a nonnegative self-adjoint operator $A$ with a nonempty resolvent set in an abstract Krein space $\cK$. We use a boundary triple approach to extension theory developed in \cite{Koch79, GG, DM91}. We construct $A$ as a coupling of two abstract symmetric operators $A_+=B_+$ and $A_-=-B_-$, where $B_+$ and $B_-$ are nonnegative symmetric operators with defect numbers $(1,1)$ acting in Hilbert spaces $\cH_+$ and $\cH_-$ which form a fundamental decomposition for $\cK$. When boundary triples $\bigl(\mathbb{C},\Gamma_0^+,\Gamma_1^+\bigr)$ and $\bigl(\mathbb{C},\Gamma_0^-,\Gamma_1^-\bigr)$ for the operators $B_+$ and $B_-$ are fixed the coupling $A$ of the operators $A_+$ and $A_-$ relative to these boundary triples is uniquely defined as a self-adjoint operator acting in the Krein space $\cK$ with the fundamental decomposition $\cK=\cH_+[+]\cH_-$, see Theorem~\ref{tHsc}.


The origins of the coupling method are twofold. On one side, it is an abstract version of an idea used by H.~Weyl \cite{Weyl12}, called Dirichlet-Neumann decoupling by B.~Simon in \cite{Simon2015}. On the other side, it is a generalization of I.M.~Glazman's decomposition method \cite{Glaz50}. The coupling method was recently extended to self-adjoint operators in Hilbert spaces in~\cite{DHMS00}. In the Krein space setting, it was used in \cite{Kar10, CD15, DT17} to study the problem of the similarity of differential operators with indefinite weights to self-adjoint operators in Hilbert spaces.


In Theorems~\ref{prop:Semibound2} and~\ref{thm:a_cs} we prove that the $D_\infty$-property is sufficient for $\infty\not\in c_s(A)$ under the assumption~\eqref{Intro:Rem_pm} and that the $D_0$-property is sufficient for $0\not\in c_s(A)$ provided that \eqref{Intro:Rem_pm} is true for all $0<y<y_0$. In Theorem~\ref{prop:Semibound2C} we prove that under the assumption~\eqref{Intro:Rem_pm} the one-sided  condition~\eqref{Intro:ImRem_pm} at $+\infty$ is sufficient for $\infty\not\in c_s(A)$. In Theorem~\ref{thm:a_cs02}  we prove analogous results for $0\not\in c_s(A)$. These results are the key stones in the proof of Theorem~A (and Theorem~B below) and are of independent interest for the coupling of two nonnegative operators in Krein spaces.


In Section~\ref{Sec:S-L}, the abstract results from Section~\ref{Sec3} are adapted to indefinite Sturm-Liouville operators. Let  $\cH_\pm$  be the weighted spaces $\cH_\pm:=L^2_{w_\pm}(I_\pm)$ and let $B_\pm$ be nonnegative symmetric operators  generated in  $\cH_\pm$ by the differential expressions
\begin{equation} \label{eq:Bpm}
\mathfrak{b}_{\pm}(f):= -  \frac{1}{w_\pm} \left( \frac{1}{r_\pm}f' \right)'
\quad \text{on} \quad I_\pm.
 \end{equation}
Using the above scheme we represent the operator $A$ as a coupling of two  symmetric operators $A_+:=B_+$  and $A_-:=-B_-$. Conditions for regularity of critical points $0$ and $\infty$ of the differential operator are formulated in terms of the functions \eqref{eq:WR}. We use the results of Bennewitz~\cite{Ben89} and Kostenko~\cite{Kost13} to reformulate one-sided sufficient conditions for regularity of critical points $\infty$ or $0$ from Theorem~\ref{prop:Semibound2C} in terms of the functions $W_\pm$ and $R_\pm$. Specifically, in Theorem~\ref{thm:Reg_infty_PIS}, we show that if either $W_+\circ R_+^{-1}$ is positively increasing at $0_+$ or $W_-\circ R_-^{-1}$ is positively increasing at $0_-$, then $\infty$ is a regular critical point for the operator $A$ associated with indefinite differential expression~\eqref{eq2ndq69A}. In Theorem~\ref{thm:Reg_infty_PI} we prove that in the case of slowly varying functions $W_\pm \circ R_\pm^{-1}$ the condition  $\infty\in c_r(A)$  is equivalent to the condition~\eqref{eq:FarFromOdd} in Theorem~A.

To show the strength of our results, in Example~\ref{ex:Reg1} we present an indefinite Sturm-Liouville operator $A$ with a non-odd weight for which Theorem~\ref{thm:Reg_infty_PI} guarantees that $\infty$ is a regular critical point, but other known criteria for regularity such as Volkmer's condition from~\cite{Vol}, Fleige's condition from~\cite{CFK13}, Parfenov's condition~\cite{Par05} cannot be applied.

In Theorem~\ref{thm:Reg0} we give a list of sufficient conditions under which we have $0\notin c_s(A)$ for the differential operator $A$.  In particular, it is shown that in the case when $w_+\in L_1(I_+)$ and $w_-\in L_1(I_-)$ the following equivalence holds
\begin{equation}\label{eq:Intro_0}
  0\not\in c_s(A) \ \ \text{and} \ \ \ker A=\ker A^2 \quad \Leftrightarrow \quad
  W_+(b_+)+W_-(b_-)\ne 0.
\end{equation}
The proof of this theorem is based on abstract results from Theorems~\ref{thm:a_cs} and~\ref{thm:a_cs02} and asymptotic formulas for the Weyl functions of the operators $B_+$ and $B_-$ from Lemmas~\ref{cor:R_0} and~\ref{lem:PolF}.

In Theorem~\ref{thm:PolF} we combine the regularity results for the points $0$ and $\infty$ to obtain new results about similarity of the operator $A$ to a self-adjoint operator in a Hilbert space. In the particular case when $w_+\in L_1(I_+)$ and $w_-\in L_1(I_-)$ these results take the folowing form

\begin{ThB} 
Let the differential expression $\mathfrak{a}$ satisfy \eqref{eq:rw} and let $W_\pm$, $R_\pm$ be the functions defined in \eqref{eq:WR}. Assume that $w_+ \in L_1(I_+)$, $w_- \in L_1(I_-)$ and one of the equivalent conditions in \eqref{eq:Intro_0} is satisfied. Then the following statements hold.
\begin{enumerate}
\renewcommand*\theenumi{\roman{enumi}}
\renewcommand*\labelenumi{\rm{(\theenumi)}}
\item \label{thmB:Riesz-B}
If either $W_+\circ R_+^{-1}$ is positively increasing at $0_+$ or $W_-\circ R_-^{-1}$ is positively increasing at $0_-$, then {\rm (Si)} holds for $A$.

\item \label{thmB:Riesz-A}
If $W_+ \circ R_+^{-1}$ is slowly varying at $0_+$ and $W_- \circ R_-^{-1}$ is slowly varying at $0_-$,
then similarity property {\rm (Si)} for $A$ is equivalent to condition \eqref{eq:FarFromOdd}.
\end{enumerate}
\end{ThB}

In Section~\ref{Sec:S-L} we systematically use results from Karamata theory about positively increasing and slowly varying functions which are presented and developed for our purposes in Appendix~\ref{appen}. In particular, it is shown that the condition for the function $W_\pm \circ R_\pm^{-1}$ to be slowly varying is equivalent to Atkinson-Bennewitz condition~\eqref{eq:Sub_RW}, see Corollary~\ref{cor:AtkCon}.

\subsection{Notation.} By $\mathbb{C}$ we denote the set of complex numbers and by $\mathbb{R}$ the set of real numbers. By $\mathbb{C}_+$ (resp. $\mathbb{C}_-$) we denote the set of all $z \in \mathbb{C}$ with positive (resp. negative) imaginary part. Similarly, $\mathbb{R}_+$ (resp. $\mathbb{R}_-$) stands for the set of all positive (resp. negative) reals. For $z \in \mathbb{C}$, $\overline{z}$, $\re z$ and $\im z$ denote the complex conjugate, real and imaginary part of $z$.

All operators in this paper are closed densely defined linear operators. For such an operator $T$, we use the common notation $\rho(T)$, $\dom(T)$, $\ran(T)$ and $\ker(T)$ for the resolvent set, the domain, the range and the null-space, respectively, of $T$.

We use the asymptotic notation little-$o$, big-$O$ and $\sim$ defined at $+\infty$ as follows:
$f(x) = o\bigl(g(x)\bigr)$ as $x \to +\infty$ if and only if $\lim_{x\to +\infty} f(x)/g(x) = 0$;
$f(x) = O\bigl(g(x)\bigr)$  as $x \to +\infty$ if and only if there exist $M, a \in \mathbb{R}_+$ such that $|f(x)| \leq M |g(x)|$ for all $x\geq a$; $f(x) \sim g(x)$ as $x \to +\infty$ if and only if $\lim_{x\to +\infty} f(x)/g(x) = 1$. Similar notation is used in the right and left neighborhood of $0$ with analogous definitions.

\section{Preliminaries} \label{Sec2}

\subsection{Definitizable operators in  Krein  spaces} 
A Krein space  $\bigl(\cK,\kip_{\cK}\bigr)$ is a complex vector space $\cK$ with a sesquilinear form $\kip_{\cK}$ such that there exist subspaces $\cH_+$ and $\cH_-$ of $\cK$ with $\bigl(\cH_+, \kip_{\cK}\bigr)$ and $\bigl(\cH_-, -\kip_{\cK}\bigr)$ being Hilbert spaces and $\cK = \cH_+[\dot{+}]   \cH_-$ is a direct and orthogonal sum; this direct orthogonal sum is called a {\em fundamental decomposition} of a Krein space $\cK$. Let $P_+$ and $P_-$ be projections associated with the direct sum $\cK = \cH_+ \dot{+} \cH_-$. The operator $J := P_+ - P_-$ is called a {\em fundamental symmetry} of a Krein space.  The space $\cK$ with the inner product  $\la x, y \ra_{\cK} = [Jx,y]_{\cK}$, $x,y \in \cK$, is a Hilbert space.  All topological notions in a Krein space refer to the topology of the Hilbert space $\bigl(\cK, \ahip_{\cK}\bigr)$.  For the general theory of Krein spaces and operators acting in them we refer to the monographs \cite{AI, B74}.
For a subspace $\cL\subset \cK$ denote by $\kappa_+(\cL)$ (resp.\  $\kappa_-(\cL)$) the least upper bound of the dimensions of positive (resp.\ negative) subspaces of $\cL$.

Let  $A$ be a linear operator in a Krein space $\bigl(\cK,\kip_{\cK}\bigr)$ with a dense domain $\dom A$. The {\em adjoint} of $A$ with respect to the inner product $\kip_{\cK}$ is denoted by $A^{[*]}$. The operator $A$ is called {\em symmetric} in $\bigl(\cK,\kip_{\cK}\bigr)$ if $A^{[*]}$ is an extension of $A$ and $A$ is called {\em self-adjoint} in $\bigl(\cK,\kip_\cK\bigr)$ if $A = A^{[*]}$.  The operator $A$ is called {\em nonnegative}  in $\bigl(\cK,\kip_{\cK}\bigr)$ if $[Af,f]_\cK\ge 0$ for all $f\in\dom A$.

According to~\cite{La} a self-adjoint operator $A$ is called {\em definitizable}, if its resolvent set $\rho(A)$ is nonempty and there exists a real polynomial $p$ such that $p(A)$ is nonnegative. Such polynomial $p$ is called definitizing polynomial of $A$. The non-real spectrum of a definitizable operator consists of a finite set of points symmetric with respect to $\mathbb{R}$.  A real number $\lambda\in\sigma(A)$  is said to be a {\em critical point} of $A$ if $p(\lambda)=0$ for every definitizing polynomial $p$ of $A$. Similarly, $\infty$  is a {\em critical point} of $A$, if at least one of its definitizing polynomials $p$ is of odd degree and the real spectrum of $A$ is neither bounded from below, nor bounded from above.  The set of critical points of $A$ is denoted by $c(A)$.

In particular, a self-adjoint  nonnegative operator $A$ with nonempty resolvent set $\rho(A)$ is definitizable with definitizing polynomial $p(\lambda)=\lambda$ and the only possible critical points of $A$ are $0$ and $\infty$.

A definitizable operator $A$ admits a spectral function $E$,
see~\cite[Theorem~II.3.1]{La}, defined on the semiring $\mathcal{R}$ generated by all intervals
whose endpoints are not critical points of $A$ with $E(\Delta)$ being self-adjoint projection in $\bigl(\cK, \kip_{\cK} \bigr)$ for every $\Delta\in \mathcal{R}$. Moreover,
\begin{equation}\label{2.5}
\bigl(E(\Delta)\cK, \kip_{\cK} \bigr) \text{ is a Hilbert space whenever }
\Delta \subset \{ t \in \mathbb{R} : p(t) > 0 \}.
\end{equation}
It follows from the properties of the spectral function $E$, see~\cite{La}, that
the restriction of $A$ to its spectral subspace $E(\Delta)\cK$ in \eqref{2.5}
is a self-adjoint operator in the Hilbert space $\bigl(E(\Delta)\cK, \kip_{\cK} \bigr)$.
A similar statement holds for intervals in $\{ t \in \mathbb R : p(t) < 0 \}$.
However, if one of the endpoints of the interval approaches a critical point, it may happen
that the norms of the corresponding spectral projections are unbounded. More precisely,
 a  point $\alpha\in c(A)$  is called a {\em
regular critical point} of $A$, if there exists a neighbourhood $G$ of $\alpha$ such that
\[
\text{the set of projections} \quad
\bigl\{
E(\Delta): \Delta \in \mathcal{R},\; \overline{\Delta}\subset G\setminus\{\alpha\} \bigr\} \quad \text{is bounded.}
\]
The set of all regular critical points of $A$ is denoted by $c_r(A)$. A critical point of $A$ which is not regular is called {\em singular critical point} of $A$. The set of all singular critical points of $A$ is denoted by $c_s(A)$. It is often difficult
to decide whether a critical point is singular or regular. A widely used characterization for $\infty\not\in c_s(A)$ is from K.~Veseli\'{c} \cite{Ves}, see also~\cite{Ako80}, \cite[Corollaries~1.5 and~1.6]{J84}. Due to the Uniform Boundedness Principle it can be reformulated as follows.

\begin{theorem} \label{Veselic}
Let $A$ be a definitizable operator in a Krein
space $({\cK}, \kip_{\cK})$ and  $\alpha\in \mathbb R$. Then:
\begin{enumerate}
\renewcommand*\theenumi{\alph{enumi}}
\renewcommand*\labelenumi{{\rm (\theenumi)}}
  \item
    $\infty\not\in c_s(A)$ if and only if there exists $\eta_0>0$ such
    that for every $f \in\cK$
  \begin{equation}\label{Veselicky}
  \int_{\eta_0}^{\eta}\re \bigl[(A-{iy})^{-1}f,f \bigr]_{\cK} dy=O(1) \quad \text{as} \quad
  \eta\to +\infty.
 \end{equation}
  \item
   $\alpha\not\in c_s(A)$ and
  $\ker(A-\alpha)=\ker\bigl((A-\alpha)^2\bigr)$ if and only if there exists $\eta_0>0$
   such that for every $f \in\cK$
  \begin{equation}\label{Veselicky2}
 \int_{\eta}^{\eta_0}\re\bigl[(A-\alpha-iy)^{-1}f,f \bigr]_{\cK} dy=O(1)
 \quad \text{as} \quad  \eta\downarrow 0.
 \end{equation}
\end{enumerate}
\end{theorem}
Let us consider
a nonnegative operator $A$ in a Krein space  $\bigl(\cK,\kip_{\cK}\bigr)$.
Then, as mentioned above, the only possible critical points are $0$ and $\infty$.
Assume that~\eqref{Veselicky} holds.
Then, by  the proof of~\cite[Lemma~1]{J82}
this implies that the set of projections $E((1,n))$ and $E((-n,-1))$,
$n\in \mathbb{N}$, is bounded,
which in turn implies $\infty\not\in c_s(A)$.
Moreover, it is easy to see that the space $ (I-E([-1,1]))\cK$  with the new inner product
\begin{equation*}
\langle f,g\rangle_{\operatorname{new}} := \lim_{n\to \infty} [(E(1,n))-E((-n,-1)f,g]
\end{equation*}
is a Hilbert space 
and that the restriction of $A$
to $(I-E([-1,1]))\cK$ is self-adjoint in the corresponding Hilbert space. A similar reasoning using \eqref{Veselicky2} holds for the point zero. This implies the following well-known statement (see, e.g.\ \cite{La}).

\begin{theorem} \label{WellKnown}
A nonnegative operator $A$ in a Krein space has similarity property {\rm (Si)} if and only if
 $\rho(A) \neq \emptyset$, $\ker(A) = \ker(A^2)$ and  $0, \infty \not\in c_s(A)$.
\end{theorem}

\subsection{Boundary triples and Weyl functions of symmetric operators} 
In this subsection $S$ is a closed densely defined symmetric operator in the Krein space $\bigl(\cK,\kip_{\cK} \bigr)$. Let $\wh\rho(S)$ denote the set of
points of regular type of $S$, see~\cite{AG}, and let $\mathfrak{N}_z$ denote
the defect subspace of the operator $S$
\begin{equation*} 
\mathfrak{N}_z : = \ran(S - \overline{z})^{[\perp]},
\quad z \in\wh\rho(S).
\end{equation*}
The numbers $n_\pm(S):=\dim(\mathfrak{N}_z)$ are
constant for all $z \in \wh\rho(S)\cap\mathbb{C}_\pm$ and are called defect numbers of ${S}$.

In what follows we assume that  the operator ${S}$  admits a
self-adjoint extension $\wt {S}$ in $\bigl(\cK,\kip_\cK\bigr)$ with a
nonempty resolvent set $\rho(\wt {S})$. Then for all $z\in\rho(\wt {S})$
we have
 \begin{equation} \label{eS*ds}
\dom(S^{[*]}) = \dom (\wt {S}) \dotplus \mathfrak{N}_z \quad \text{ direct sum in }
\quad \cH.
 \end{equation}
This implies, in particular, that the dimension $\dim(\mathfrak{N}_z)$ is
constant for all $z \in \rho(\wt {S})$  and hence $n_+(S)=n_-(S)$. Moreover, we assume everywhere in this paper that $n_\pm(S)=1$.
Notice that the equality $n_+(S)=n_-(S)$
does not imply the existence of a self-adjoint extension $\wt {S}$ of ${S}$, see~\cite{Shm74}.

\begin{definition}\label{def:BTriple}
Let $\Gamma_0$ and $\Gamma_1$ be linear mappings from $\dom(S^{[*]})
$ to ${\mathbb C}$ such that
\begin{enumerate}
\def\labelenumi{\rm (\roman{enumi})}
\item the mapping $\Gamma:\, f \to \begin{pmatrix}
                                     \Gamma_0f\\
                                     \Gamma_1f
                                   \end{pmatrix}$
      from $\dom(S^{[*]})$ to ${\mathbb C}^{2}$ is surjective;
\item the abstract Green's identity
\begin{equation}
\label{Green}
  \bigl[S^{[*]}f,g\bigr]_\cK - \bigl[f,S^{[*]} g\bigr]_\cK=
  (\Gamma_1{f})\overline{{(\Gamma_0{g})}}\,
   -(\Gamma_0{f})\overline{{(\Gamma_1{g})}}\,
\end{equation}
holds for all $f$, $g\in \dom(S^{[*]})$.
\end{enumerate}
Then the triple $\bigl(\mathbb{C},\Gamma_0,\Gamma_1\bigr)$ is called
a {\em boundary triple} for ${S}^{+}$,
 see \cite{GG,DM91, D95} for much
more general setting.
\end{definition}

It follows from \eqref{Green} that the extensions ${S}_0$, ${S}_1$ of
${S}$ defined as restrictions of ${S}^{+}$ to the domains
$\dom(S_0) :=  \ker(\Gamma_0)$ and $\dom(S_1) :=\ker(\Gamma_1)$
are self-adjoint extensions of ${S}$.

Given a self-adjoint extension $\wt {S}$ of ${S}$ with nonempty $\rho(\wt {S})$ one can always choose  a boundary triple $\bigl(\mathbb{C},\Gamma_0,\Gamma_1\bigr)$ for ${S}$ such that ${S}_0=\wt {S}$, see~\cite[Proposition~2.2]{D99}.
In this case  for every
$z\in\rho(S_0)$ the decomposition~\eqref{eS*ds} holds with $\wt
{S}={S}_0$ and the mapping $\Gamma_0|_{\mathfrak{N}_z}:\mathfrak{N}_z\to\mathbb{C}$ is invertible for every $z\in\rho(S_0)$. A vector-valued function $z\mapsto\gamma(z)$ defined on $\rho(S_0)$ with values in $\mathfrak{N}_z$  is called the
$\gamma$-field of ${S}$, associated with the boundary triple $\bigl(\mathbb{C},\Gamma_0,\Gamma_1\bigr)$ if
\[
\Gamma_0\gamma(z) = 1  \quad \text{for all} \quad z\in\rho(S_0).
\]
Notice, that $\gamma$ satisfies the equality, see~\cite[Proposition~2.2]{D99},
\begin{equation}\label{eq:gamma}
        \gamma(z)=(S_0-z_0)(S_0-z)^{-1}\gamma(z_0),\quad z,z_0\in\rho(S_0)
\end{equation}
and hence the vector-valued function $\gamma$ is holomorphic on $\rho(S_0)$.
\begin{definition} \label{W00A} The  function $z\mapsto M(z)$ defined by the equality
\begin{equation*}
M(z)\Gamma_0f_{z}=\Gamma_1f_{z},
 \quad
f_z\in\mathfrak{N}_z,\,z\in\rho(S_0),
\end{equation*}
is called the \textit{abstract Weyl function} of $S$,
corresponding to the boundary triple
$\bigl(\mathbb{C},\Gamma_0,\Gamma_1\bigr)$.
\end{definition}
The notion of the  abstract Weyl function was introduced in~\cite{DM91} for a Hilbert space symmetric operator and in~\cite{D95} for a Krein space symmetric operator.

Clearly,
$M(z) = \Gamma_1\gamma(z)$ for $z\in\rho(S_0),$
and hence $M(z)$ is well defined. It follows from \eqref{Green} and \eqref{eq:gamma}  that the Weyl function $M$  satisfies the identity
\begin{equation}\label{eQfun}
        M(z)-\overline{M(w)}=(z-\overline{w})\, [\gamma(z),\gamma(w)]_\cK,
         \quad         z,w \in  \rho(S_0).
\end{equation}
With $w = \overline{z}$ the identity \eqref{eQfun} yields that the Weyl function $M$ satisfies the symmetry condition
\begin{equation} \label{esym}
M(\overline{z}) = \overline{M(z)},  \quad  z \in \rho(S_0).
\end{equation}

In the case when $\bigl(\cH,\ahip_\cH\bigr)$
is a Hilbert space we  will use the notation $B$ for a closed densely
defined symmetric operator in the Hilbert  space $\cH$
with defect numbers $(1,1)$.
Let $\bigl(\mathbb{C},\Gamma_0,\Gamma_1\bigr)$  be a boundary triple for $B^{\langle *\rangle}$. We will use the notations
$m$ and $\gamma_B$ for the abstract Weyl function and for the $\gamma$-field of $B$
corresponding to the boundary triple $\bigl(\mathbb{C},\Gamma_0,\Gamma_1\bigr)$.
It follows from \eqref{eQfun} and  \eqref{esym}
that $m$ is a Nevanlinna function, see~\cite{KaKr74}, i.e.\
$m$ is holomorphic at least on $\mathbb{C}\setminus\mathbb{R}$ and satisfies the following two conditions
\begin{equation*}
m(\overline{z}) = \overline{m(z)} \quad \text{and}
\quad \im m(z) \ge 0, \qquad
z \in \mathbb{C}_+.
\end{equation*}
Since the operator $B$ is densely defined the following two conditions hold  (see~\cite[Theorem~7.36]{DM17})
     \begin{equation}\label{6.6}
\lim_{y\uparrow+\infty}y^{-1}m(i y)=0, \quad \lim_{y\uparrow+\infty}y \im m(iy)= +\infty.
     \end{equation}

Assume that the operators $B$  and its self-adjoint extension $B_0$ with the domain $\dom B_0=\ker \Gamma_0$ are nonnegative.
Then the Weyl function $m$ is holomorphic on $\mathbb{R}_-$.
A Nevanlinna function $m$ with the above property which, in addition, takes nonnegative values for all ${z}\in\mathbb{R}_-$
is called a Stieltjes function. The class of all Stieltjes functions is denoted by $\cS$.

A Stieltjes function $m$ admits the integral representation, \cite{KaKr74},
\begin{equation}\label{eq:int_S}
m(z) = a+\int_{0}^{+\infty}  \frac{d\sigma(t)}{t-z}
\end{equation}
with $a \ge 0$ and with a non-decreasing left-continuous function $\sigma(t)$, such that
$
    \int_{0}^{+\infty} \frac{d\sigma(t)}{1+t}
$ converges.
The following statement is immediate from \eqref{eq:int_S}.

\begin{proposition}\label{p:SiSm}
Let $m\in\cS$ and assume that the support of $d\sigma$ has a nonempty intersection with $\mathbb{R}_+$. Then $\re m(iy)> 0$ for all $y\in\mathbb{R}_+$.
\end{proposition}

\subsection{Real operators}
Recall the notions of real operator and real vector valued function with respect to some conjugation, see~\cite[Section III.5]{EE87} and \cite{DM91,Koch79}.
\begin{definition}\label{def:3.4}
An involution $j_\cK$ on a Krein space $\bigl(\cK,\kip_{\cK} \bigr)$ is called a {\em conjugation} on $\cK$ if
\begin{equation}\label{eq:Conj}
\left[ j_\cK f,j_\cK g \right]_\cK = \left[ g,f \right]_\cK \quad \text{for all} \quad
 f, g \in \cK .
\end{equation}
A closed operator $T$ in a Krein space  $\cK$ is called {\it real}, if
\begin{equation*}
  j_\cK\dom(T)=\dom(T)\quad\text{and}\quad j_\cK T=Tj_\cK.
\end{equation*}
\end{definition}

Every conjugation is an anti-linear operator, see~\cite[Section IX.2]{St32}, i.e.
\[
j_\cK(\lambda f + \mu g)=\overline{\lambda} j_\cK f + \overline{\mu} j_\cK g \quad \text{for
    all} \quad f, g \in\cK, \  \lambda,\mu \in \mathbb{C}.
\]
If $T$ is real and densely defined then its adjoint $T^{[*]}$ is also a real  operator in $\cK$.

A vector $f$ in $\cK$ is called {\it real} with respect to the conjugation $j_\cK$, if $j_\cK f=f$. An arbitrary vector $f\in\cK$ can be decomposed into the sum
\begin{equation}\label{eq:Real_Decom}
  f=f^R+i f^I,\quad\text{where} \quad f^R =\frac{1}{2}(f + j_\cK f) \quad \text{and} \quad  f^I = \frac{1}{2i}(f - j_\cK f) \quad\text{are real}.
\end{equation}

Let $j$ be the standard conjugation in $\mathbb{C}$, $j z=\overline{z}$ for all $z \in\mathbb{C}$. A scalar function $z\mapsto m(z)$ is called {\it real}, if its domain is symmetric with respect to $\mathbb{R}$ and $m(\overline{z})=\overline{m(z)}$ for all $z$ in the  domain of $m$. Similarly, a vector valued function $z\mapsto\gamma(z)$ with the values in $\cK$ is called {\it real} if its domain is symmetric with respect to $\mathbb{R}$ and
\begin{equation}\label{eq:3.10A}
\gamma(\overline{z})=j_\cK \gamma(z)
\end{equation}
for all $z$ in the domain.

Let a symmetric operator $S$ be real in $\cK$ with the conjugation $j_\cK$. A boundary triple
$\bigl(\mathbb{C},\Gamma_0,\Gamma_1\bigr)$ for $S^{[*]}$ is called {\it real}, if
\[
j\Gamma_0=\Gamma_0j_\cK \quad\text{and}\quad
j\Gamma_1=\Gamma_1j_\cK.
\]
Every real symmetric operator $S$ admits a real boundary triple $\bigl(\mathbb{C},\Gamma_0,\Gamma_1\bigr)$ and the corresponding Weyl function $M$ and the $\gamma$-field $\gamma$ are real, see~\cite{Koch79} for the case of a Hilbert space $\cK$.

\section{Regularity of critical points of couplings in Krein spaces} \label{Sec3}

\subsection{Couplings of symmetric operators in Krein spaces} \label{SecCoup}
In this section we consider two Krein spaces $\bigl( \cK_+,
\kip_{\cK_+} \bigr)$ and $\bigl( \cK_-, \kip_{\cK_-} \bigr)$.
Let their direct sum
\begin{equation*} 
    \cK = \cK_+[\dot{+}] \cK_-
\end{equation*}
be endowed with the natural inner product
\begin{equation}\label{eq:Fund_inner}
[f_++f_-,g_++g_-]_{\cK}:=[f_+,g_+]_{\cK_+}+[f_-,g_-]_{\cK_-},\quad
f_\pm,g_\pm\in\cK_\pm.
\end{equation}
Consider two  closed symmetric densely defined operators $A_+$ and $A_-$  with defect numbers $(1,1)$ acting in the Krein spaces $\bigl( \cK_+,
\kip_{\cK_+} \bigr)$ and $\bigl( \cK_-, \kip_{\cK_-} \bigr)$.
Let $\bigl(\mathbb{C},\Gamma_{0}^\pm, \Gamma_{1}^\pm\bigr)$ be a boundary
triple for $A_\pm^{[*]}$. Let $M_\pm$ and $\gamma_{A_\pm}$ be the
corresponding Weyl function and  the $\gamma$-field. By $A_{\pm,0}$ we
denote the self-adjoint extension of $A_\pm$ which is defined on
\[
\dom(A_{\pm,0}) = \ker(\Gamma_{0}^\pm)\quad\text{by}\quad A_{\pm,0}
=A_\pm^{[*]}\bigl|_{\ker(\Gamma_{0}^\pm)}\bigr..
\]
Then the functions $M_\pm$ are defined and holomorphic on $\rho(A_{\pm,0})$. Assume that \begin{equation}\label{eq:rho_Apm}
\rho(A_{+,0})\cap\rho(A_{-,0})\ne\emptyset.
\end{equation}

The following theorem is an indefinite version of results from~\cite[Proposition~4.3]{DHMS00} 
which is, in this form, presented in~\cite[Theorem~4.7]{CD15} and~\cite{DT17}.

\begin{theorem} \label{tHsc}
 Let $A_\pm$  be closed
symmetric densely defined operators  with defect numbers $(1,1)$ in the Krein spaces $\cK_\pm$.
Let $\bigl(\mathbb{C},\Gamma_{0}^\pm, \Gamma_{1}^\pm\bigr)$  be boundary triples for $A_\pm^{[*]}$ which satisfy \eqref{eq:rho_Apm}. Let $M_\pm$ and $\gamma_{A_\pm}$ be the Weyl functions and the $\gamma$-fields of $A_\pm$ corresponding to the boundary triples $\bigl(\mathbb{C},\Gamma_{0}^\pm, \Gamma_{1}^\pm\bigr)$, and let $S$ and $A$  be the restrictions  of $A_+^{[*]} [+]
A_-^{[*]}$ to the domains
\begin{equation} \label{eq:domSHC}
\dom(S) = \left\{\begin{pmatrix} f_+ \\ f_- \end{pmatrix}:
  \begin{array}{l}
 \Gamma_{0}^+(f_+) = \Gamma_{0}^- (f_-)=0,  \\
 \Gamma_{1}^+(f_+) + \Gamma_{1}^- (f_-)=0,
 \end{array} \     f_\pm\in \dom\bigl(A_\pm^{[*]}\bigr)
 \right\},
\end{equation}
\begin{equation} \label{eq:domwAHC}
\dom(A) = \left\{\begin{pmatrix} f_+ \\ f_- \end{pmatrix}:
  \begin{array}{l}
 \Gamma_{0}^+(f_+) = \Gamma_{0}^- (f_-),  \\
 \Gamma_{1}^+(f_+) + \Gamma_{1}^- (f_-)=0,
 \end{array} \   f_\pm\in \dom\bigl(A_\pm^{[*]}\bigr)
 \right\}.
\end{equation}
Then the following statements hold:
\begin{enumerate}
\renewcommand*\theenumi{\alph{enumi}}
\renewcommand*\labelenumi{{\rm (\theenumi)}}
\item \label{ietHsc}
The operator $S$ is  symmetric with defect numbers $(1,1)$ and $A$ is a self-adjoint extension of $S$.
\item 
The adjoint ${S}^+$ of ${S}$ is the restriction of $A_+^+ [+] A_-^+$ to the domain
\begin{equation*}
\dom(S^+) = \left\{ \begin{pmatrix} f_+ \\ f_-
\end{pmatrix}:\,
 \Gamma_{0}^+(f_+) = \Gamma_{0}^- (f_-), \  f_\pm\in \dom\bigl(A_\pm^{[*]}\bigr) \right\}.
\end{equation*}
and a boundary triple $\bigl(\mathbb{C},\Gamma_0,\Gamma_1\bigr)$ for $S^{[*]}$ is
given by
\begin{equation} \label{Htripl}
\Gamma_0 f = \Gamma_0^+ f_+, \quad
 \Gamma_1 f = \Gamma_1^+ f_+ + \Gamma_1^-f_-, \quad
  f = \begin{pmatrix} f_+ \\ f_- \end{pmatrix}  \in \dom\bigl(S^{[*]}\bigr).
\end{equation}
\item 
The corresponding Weyl function  and the $\gamma$-field of $S$ are
\begin{equation} \label{mH2}
 M(z) = M_+(z) + M_-(z), \quad \gamma(z)=\left(%
\begin{array}{c}
 \mkern-6mu \gamma_{A_+}(z) \mkern-6mu \\
 \mkern-6mu \gamma_{A_-}(z) \mkern-6mu
\end{array}%
\right), \quad z \in \mathbb{C} \setminus \mathbb{R}.
\end{equation}
\item \label{irho-i}
If $z\in\rho(A_{+,0})\cap\rho(A_{-,0})$ then $z\in\rho(A)  $  if and only if
$  M_+(z) + M_-(z) \ne 0$.
\item
The resolvent of the operator $A$ is given by
\begin{equation} \label{ereswA1}
\bigl(A - z \bigr)^{-1} f = \bigl(A_0 - z \bigr)^{-1} f
  - \frac{[f,\gamma(\overline{z})]_\cK}
  {M_+(z)+M_-(z)}\gamma(z),\ \
\ z \in \rho(A) \cap\rho(A_0),
 \end{equation}
where   $A_0=A_{+,0}[+]A_{-,0}$ and $f\in\cK$.
\end{enumerate}
\end{theorem}

\begin{definition} \label{def:coup}
The operator $A$ defined in Theorem~\ref{tHsc}(\ref{ietHsc}) is called the {\it coupling of the operators $A_+$ and $A_-$ in the Krein space $\cK$ 
relative to the triples $\bigl(\mathbb{C},\Gamma_{0}^+, \Gamma_{1}^+\bigr)$ and $\bigl(\mathbb{C},\Gamma_{0}^-, \Gamma_{1}^-\bigr)$} and  $A_0=A_{+,0}[+]A_{-,0}$  is called the {\it decoupled} operator.
\end{definition}

The following statement was proved in~\cite[Lemma~5.4]{Kost13}. For the reader's convenience, we present a proof based on Theorem~\ref{tHsc}.
\begin{lemma}\label{lem:Kost1}
Let $(\cH,\ahip_\cH)$ be a Hilbert space with a conjugation $j_{\cH}$, let $B$ be a closed densely defined real symmetric operator in $\cH$ with defect numbers $(1,1)$,  let
$\bigl(\mathbb{C},\Gamma_0,\Gamma_1\bigr)$ be a real  boundary triple for $B^{\langle *\rangle}$,
let $m$ and $\gamma_{B}$ be the corresponding Weyl function and  the $\gamma$-field for $B$ and define
\begin{equation}\label{eq:f}
 \wh h(z) = \bigl\langle h,\gamma_{B}(\overline{z}) \bigr\rangle_{\cH}, \quad h \in \cH, \quad z \in \mathbb{C}\setminus\mathbb{R}.
\end{equation}
Then the following inequality holds for all real $h \in\cH$:
\begin{equation}\label{eq:Ineq1}
\int_{0}^\infty\frac{\bigl|\im \wh h(iy)^2\bigr|}{\im m(iy)}dy \le 2\pi \|h\|^2_{\cH},
\end{equation}
\end{lemma}
\begin{proof}
Let $\cK_+$ and $\cK_-$ be two copies of  the Hilbert space $\cH$
and let us  set $A_+:=B$ and $A_-:=-B$. Notice that $\bigl(\mathbb{C},\Gamma_0,\Gamma_1\bigr)$ is a boundary triple for $A_+^{[*]}$, $\bigl(\mathbb{C},\Gamma_0,-\Gamma_1\bigr)$ is a boundary triple for $A_-^+$  and the corresponding Weyl functions $M_+$, $M_-$ and the $\gamma$-fields $\gamma_{A_+}$ and $\gamma_{A_-}$  take the form
\[
M_+(z)=m(z),\quad M_-(z)=-m(-z),\quad \gamma_{A_+}(z)=\gamma_{B}(z),\quad\gamma_{A_-}(z)=\gamma_{B}(-z).
\]

Let $A$ be the coupling of $A_+$, $A_-$ acting in the Hilbert space  $\cK=\cK_+\oplus\cK_-=\cH\oplus\cH$, let $A_0=B_{0}\oplus (-B_{0})$ be the decoupled operator as defined in Definition~\ref{def:coup}, $B_{0}$ being the restriction of $B^{\langle *\rangle}$ to $\ker\Gamma_0$ and let us denote by $\ahip_\cK$ the scalar product in $\cK=\cK_+\oplus\cK_-=\cH\oplus\cH$. Applying Theorem~\ref{tHsc} to the operators  $A_+$ and $A_-$ one obtains from~\eqref{ereswA1} for vector $f=h\oplus 0$, $h\in\cH$ the equality
\begin{equation} \label{eq:resA2}
\bigl\langle (A - iy )^{-1} f,f \bigr\rangle_{\cK}
 = \bigl\langle (A_0 - iy)^{-1} f,f \bigr\rangle_{\cK}
  - \frac{\wh h(iy) \overline{\wh h(- iy)}}{m(iy)-m(-iy)}
 \end{equation}
Since $A$ and $A_0$ are self-adjoint operators in the Hilbert space $\cK =\cH\oplus\cH$, an application of the functional calculus yields
\begin{equation}\label{eq:Res_SA}
\int_0^\infty \left|\re\bigl\langle (A-iy)^{-1}f,f \bigr\rangle_\cK \right|dy
\le  \frac{\pi}{2}\|f\|_{\cK}^2,
\end{equation}
\begin{equation}\label{eq:Res_SA0}
\int_0^\infty\left|\re \bigl\langle({ A_0}-iy)^{-1}f,f\bigr\rangle_\cK \right|dy \le \frac{\pi}{2}\|f\|_{\cK}^2.
\end{equation}
for all $f\in\cK$.
Since the boundary triple $\bigl(\mathbb{C},\Gamma_0,\Gamma_1\bigr)$ is real, $\gamma_{B}$ is real as well. If, in addition, $h$ is real, then
\[
j_{\cH}h=h,\quad j_{\cH}\gamma_{B}(iy)=\gamma_{B}(-iy)\quad\text{for all}\quad y\in\mathbb{R}_+
\]
and by  \eqref{eq:f} and Definition~\ref{def:3.4}
\begin{equation}\label{eq:Wh_bar}
\overline{\wh h(- iy)}  = \bigl\langle\gamma_{B}( iy),h \bigr\rangle_{\cH}
 = \bigl\langle j_{\cH}h,j_{\cH}\gamma_{B}(iy) \bigr\rangle_{\cH}
 = \bigl\langle h,\gamma_{B}(- iy) \bigr\rangle_{\cH}
 = \wh h(iy).
\end{equation}
By~\eqref{eq:resA2}, \eqref{eq:Res_SA}, \eqref{eq:Res_SA0} and \eqref{eq:Wh_bar}
\begin{equation*}
  \int_0^\infty\left|\re\frac{\wh h( iy)^2}
  {m(iy)-m(-iy)}\right|dy\le {\pi}\|f\|_{\cH}^2={\pi}\|h\|_{\cH}^2.
\end{equation*}
Using the equality $m(iy)-m(-iy)=2i \im m(iy)$ one obtains for all real $h \in\cH$:
\begin{equation*} 
\int_{0}^\infty\frac{\bigl|\im\wh h(iy)^2\bigr|}{2\im m(iy)}dy =   \int_0^\infty\left|\re\frac{\wh h( iy)^2}
  {m(iy)-m(-iy)}\right|dy
  \le {\pi}\|h\|_{\cH}^2.
\end{equation*}
This proves \eqref{eq:Ineq1}.
\end{proof}

In the following lemma we apply Theorem~\ref{tHsc} to two real symmetric operators $B_+$ and $B_-$ acting in Hilbert spaces $\cH_+$ and $\cH_-$ and obtain  estimates for a family of weighted $L^2$-norms of ``generalized Fourier transforms''
\begin{equation}\label{eq:fpm}
 \wh f_\pm(z) = \bigl\langle f_\pm,\gamma_{B_\pm}(\overline{z}) \bigr\rangle_{\cH_\pm}, \quad f_\pm \in \cH_\pm, \quad z \in \mathbb{C}\setminus\mathbb{R}.
\end{equation}

\begin{lemma}\label{lem:Kost2}
Let $\cH_\pm$ be Hilbert spaces with conjugations $j_{\cH_\pm}$, let $B_\pm$ be closed densely defined real symmetric operators in $\cH_\pm$ with defect numbers $(1,1)$,  let $\bigl(\mathbb{C},\Gamma_0^\pm,\Gamma_1^\pm\bigr)$ be real  boundary triples for $B_\pm^{\langle *\rangle}$, and let $m_\pm$ and $\gamma_{B_\pm}$ be the corresponding Weyl functions and  the $\gamma$-fields for $B_\pm$. Then the following inequalities hold for all real $f_\pm \in\cH_\pm$:
\begin{equation}\label{eq:Ineq2}
    \int_{0}^{+\infty}\left|\wh
    f_\pm(iy)\right|^2\frac{\left|\re\bigl(m_+(iy)+{m_-(iy)}\bigr)\right|}
    {\left|m_+(iy)+{m_-(iy)}\right|^2}dy\
     \le 5\pi \|f_\pm\|^2_{\cH_\pm},
\end{equation}
\end{lemma}
\begin{proof}
\noindent{\bf 1.}
In this step we prove that for all real $f_\pm\in\cH_\pm$ we have
\begin{equation}\label{eq:Ineq2Re}
    \int_{0}^{+\infty}\left|\re \bigl( \wh
    f_\pm(iy)^2\bigr)\right|\frac{\left|\re(m_+(iy)+{m_-(iy)})\right|}
    {\left|m_+(iy)+{m_-(iy)}\right|^2}dy\
     \le 3\pi \|f_\pm\|^2_{\cH_\pm}.
\end{equation}
Applying Theorem~\ref{tHsc} to the operators  $A_+:=B_+$ and $A_-:=B_-$  in Hilbert spaces $\cK_+=\cH_+$, $\cK_-:=\cH_-$ and taking $f=f_+\oplus f_-$, $f_\pm\in\cH_\pm$, one obtains the equality
\[ 
\bigl\langle (A - iy )^{-1} f,f \bigr\rangle
 = \bigl\langle (A_0 - iy )^{-1} f,f \bigr\rangle
  - \frac{\bigl(\wh f_{+}(iy)+\wh f_{-}(iy)\bigr)\overline{\bigl(\wh f_{+}(-iy) +\wh f_{-}(-iy)\bigr)}}{m_+(iy)+m_-(iy)}
\]
where $A$ is the coupling of $A_+$ and $A_-$ defined by~\eqref{eq:domwAHC}, ${A_0}$
is the decoupled operator, as defined in Definition~\ref{def:coup}.
Since $A$ and ${ A_0}$ are self-adjoint operators in the Hilbert space $\cH:=\cH_+\oplus\cH_-$ one obtains from~\eqref{eq:Res_SA}, \eqref{eq:Res_SA0} and \eqref{eq:Wh_bar} for all real $f_\pm\in\cH_\pm$
\begin{equation}\label{eq:Re2}
  \int_0^{+\infty}\left|\re\frac{\bigl(\wh f_+(iy)+\wh f_-(iy)\bigr)^2}
  {m_+(iy)+m_-(iy)}\right|dy\le\pi\|f\|_{\cH}^2.
\end{equation}
Set
\begin{equation}\label{eq:ReIm_m}
  u_\pm(iy):=\re m_\pm(iy),\quad v_\pm(iy):=\im m_\pm(iy),
\end{equation}
\begin{equation*}
  U(iy):=\re \bigl((\wh f_+(iy)+\wh f_-(iy))^2\bigr),\quad V(iy):=\im \bigl((\wh f_+(iy)+\wh f_-(iy))^2\bigr).
\end{equation*}
Then inequality~\eqref{eq:Re2} can be rewritten as
\begin{equation*} 
    \int_{0}^{+\infty}\frac{\left|U(iy) \bigl(u_+(iy)+u_-(iy)\bigr)+ V(iy) \bigl(v_+(iy)+{v_-(iy)}\bigr)\right|}
    {\left|m_+(iy)+ m_-(iy)\right|^2}dy  \le \pi\|f\|_{\cH}^2.
\end{equation*}
In particular, setting subsequently $f_-=0$ or $f_+=0$, one obtains
\begin{multline}\label{eq:Ineq3pm}
    \int_{0}^{+\infty}\frac{\left|\re \bigl( \wh f_\pm(iy)^2 \bigr)
    \bigl( u_+(iy)+{u_-(iy)} \bigr)
    + \im \bigl( \wh f_\pm( iy)^2 \bigr) \bigl( v_+(iy)+{v_-(iy)} \bigr) \right|}
    {\left|m_+(iy)+{m_-(iy)}\right|^2}dy  \\
    \le\pi\|f_\pm\|_{\cH_\pm}^2.
\end{multline}
By \eqref{eq:Ineq1} for every real $f_\pm\in\cH_\pm$
\begin{equation}\label{eq:Ineq3A}
\int_{0}^{+\infty} \frac{\left|\im \bigl(\wh f_\pm(iy)^2\bigr)%
\bigl(v_+(iy)+v_-(iy)\bigr)\right|}{\left|m_+(iy)+m_-(iy)\right|^2}dy
 \le
\int_{0}^{+\infty}\frac{\left|\im \bigl( \wh f_\pm( iy)^2\bigr) \right|}{\im m_\pm(iy)}dy\
    \le 2\pi\|f_\pm\|_{\cH_\pm}^2
\end{equation}
and thus \eqref{eq:Ineq3pm} and \eqref{eq:Ineq3A} imply for every real $f_\pm\in\cH_\pm$
\begin{equation*} 
\int_{0}^{+\infty}\left|\re \bigl( \wh f_\pm( iy)^2 \bigr) \right|
\frac{\left|u_+(iy)+u_-(iy)\right|}{\left|m_+(iy)+{m_-(iy)}\right|^2}dy
     \le 3\pi \|f_\pm\|^2_{\cH_\pm},
\end{equation*}
which proves~\eqref{eq:Ineq2Re}.

\noindent{\bf 2.} To prove \eqref{eq:Ineq2} we notice that from
\begin{equation*}
\frac{ |u_+(iy)+{u_-(iy)}|}{\left|m_+(iy)+{m_-(iy)}\right|}\le 1
\end{equation*}
and \eqref{eq:Ineq1} we obtain
\begin{equation}\label{eq:Whf_H1}
\begin{split}
\int_{0}^{+\infty} \left|\im \bigl( \wh f_\pm (iy)^2 \bigr) \right|
\frac{|u_+(iy)+{u_-(iy)}|}{\left|m_+(iy)+{m_-(iy)}\right|^2}dy
 & \le
 \int_{0}^{+\infty}
     \frac{\left|\im \bigl(\wh f_\pm(iy)^2\bigr)\right|}{\left|m_+(iy)+{m_-(iy)}\right|}dy \\
 & \le
    \int_{0}^{+\infty}
\frac{\left|\im \bigl( \wh f_\pm(iy)^2\bigr) \right|}{\im m_\pm(iy)}dy \\
 & \le 2\pi\|f_\pm\|_{\cH_\pm}^2,
\end{split}
\end{equation}
for all real $f_\pm\in\cH_\pm$. Now \eqref{eq:Ineq2} follows from \eqref{eq:Ineq2Re} and~\eqref{eq:Whf_H1}.
\end{proof}

\begin{remark}
\begin{enumerate}
\renewcommand*\theenumi{\alph{enumi}}
\renewcommand*\labelenumi{{\rm (\theenumi)}}
  \item
It follows from \eqref{eq:Wh_bar} that the inequalities \eqref{eq:Ineq1} in Lemma~\ref{lem:Kost1} and \eqref{eq:Ineq2} in Lemma~\ref{lem:Kost2} remain in force when we substitute $\wh f_\pm(iy)$ by $\wh f_\pm(-iy)$ for all real $f_\pm\in\cH_\pm$:
  \begin{equation*}
\int_{0}^{+\infty}
    \frac{\left| \im \bigl( \wh f_\pm(-iy)^2 \bigr)\right|}{ \im m_\pm(iy)}dy \
    \le \ 2 \pi \|f_\pm\|^2_{\cH_\pm},
\end{equation*}
\begin{equation*}
    \int_{0}^{+\infty}\left|\wh
    f_\pm(-iy)\right|^2\frac{\left|\re (m_+(iy)+{m_-(iy)})\right|}
    {\left|m_+(iy)+{m_-(iy)}\right|^2}dy\
    \le 5\pi \|f_\pm\|^2_{\cH_\pm}.
\end{equation*}
\item
Notice that the statement of Lemma 3.3  is essentially contained in~\cite{Kost13} for the case when $A$ is a coupling of Sturm-Liouville operators. The result of Lemma~3.4 is new and in the symmetric case, when $m_+=m_-$, implies the statement of Corollary~5.6 in~\cite{Kost13}.\\
\end{enumerate}
\end{remark}

\subsection{Veseli\'{c} condition and coupling}
In the rest of this section we make the following general assumptions:

\begin{enumerate}
\renewcommand*\theenumi{\arabic{enumi}}
\renewcommand*\labelenumi{{\rm (A\theenumi)}}
\item \label{assume1}
 $\bigl(\cH_\pm, \ahip_{\cH_\pm}\bigr)$ are Hilbert spaces with conjugations $j_{\cH_\pm}$ and $\bigl( \cK_\pm, \kip_{\cK_\pm} \bigr)$ are Krein spaces defined by
\begin{equation*}
\cK_+ = \cH_+, \quad \kip_{\cK_+} = \ahip_{\cH_+}, \quad
\cK_- = \cH_-, \quad \kip_{\cK_-} = -\ahip_{\cH_-}.
\end{equation*}
\item
$\bigl( \cK, \kip_{\cK} \bigr)$ is a  Krein space with the fundamental decomposition $\cK:= \cK_+[+]\cK_-$ and the inner product~\eqref{eq:Fund_inner}.
\item \label{assume3}
$B_\pm$ are real closed nonnegative symmetric  densely defined operators with defect indices $(1,1)$ in the Hilbert spaces $\bigl(\cH_\pm, \ahip_{\cH_\pm}\bigr)$  and let $A_+$ and $A_-$ be symmetric operators  in the Krein spaces $\cK_+$ and $\cK_-$, respectively:
\begin{equation*}
    A_+:=B_+,\quad A_-:=-B_-.
\end{equation*}
  \item
$\bigl(\mathbb{C},\Gamma_0^\pm,\Gamma_1^\pm\bigr)$ are real boundary triples for
$B_\pm^{\langle *\rangle}$ and $m_\pm$ and $\gamma_{B_\pm}$ are the corresponding Weyl functions and the $\gamma$-fields.
\item \label{assumelast}
$A$ is the coupling of the operators $A_+$ and $A_-$ in the Krein space
$\bigl( {\cK}, \kip_{{\cK}} \bigr)$ relative to the triples
$\bigl(\mathbb{C},\Gamma_{0}^+, \Gamma_{1}^+\bigr)$ and $\bigl(\mathbb{C},\Gamma_{0}^-,\Gamma_{1}^-\bigr)$.
\end{enumerate}

By $B_{\pm,0}$ we denote the self-adjoint extension of $B_\pm$ which is defined on
\[
\dom(B_{\pm,0}) = \ker(\Gamma_{0}^\pm)\quad\text{by}\quad B_{\pm,0}
= B_\pm^{\langle *\rangle}\bigl|_{\ker(\Gamma_{0}^\pm)}\bigr..
\]
Clearly, $\bigl(\mathbb{C},\Gamma_0^\pm,\Gamma_1^\pm\bigr)$ are also boundary triples for $A_\pm^{[*]}$. The Weyl functions $M_\pm$  and  the $\gamma$-fields $\gamma_{A_\pm}$ of the operators $A_{\pm}$ corresponding to $\bigl(\mathbb{C},\Gamma_0^\pm,\Gamma_1^\pm\bigr)$ are connected with the Weyl functions $m_\pm$ and  the $\gamma$-fields $\gamma_{B_\pm}$ of the operators $B_{\pm}$ by the equalities
\begin{equation*} 
    M_\pm(z)=m_\pm(\pm z),\quad
    \gamma_{A_\pm}(z)=\gamma_{B_\pm}(\pm z), \quad z\in\rho(B_{\pm,0}).
\end{equation*}

In the next lemma we reformulate the Veseli\'{c} condition from Theorem~\ref{Veselic} for the coupling $A$ of two nonnegative operators as defined in Definition~\ref{def:coup}, cf.~\cite{DT17} and~\cite{CD15}.

\begin{lemma}\label{thm:Reg_infty}
Let conditions {\rm (A\ref{assume1})} through {\rm (A\ref{assumelast})} be satisfied. Then the following statements hold:
\begin{enumerate}
\renewcommand*\theenumi{\roman{enumi}}
\renewcommand*\labelenumi{{\rm (\theenumi)}}
  \item  \label{Coupl}
The coupling $A$ is definitizable and $\infty\in c(A)$.
  \item  \label{iVeselic-1}
$\infty\in c_r(A)$ if and only if there exists $\eta_0>0$ such that for all real $f_\pm\in\cH_\pm$
\begin{equation}\label{eq:Crit_Inf1}
    \int_{\eta_0}^{\eta}\re \frac{ \bigl(\wh f_{+}(iy)+\wh f_{-}(-iy)\bigr)^2}{m_+(iy)+m_-(-iy)} dy = O(1)\quad\text{as}\quad \eta\to+\infty.
\end{equation}

\item \label{iVeselic-2}
 $0\not\in c_s(A)$ and
  $\ker A=\ker A^2$ if and only if   there is $\eta_0>0$  such that  for all real $f_\pm\in\cH_\pm$
\begin{equation}\label{eq:Crit_Inf2}
 \int_{\eta}^{\eta_0}\re \frac{\bigl(\wh f_{+}(iy)+\wh f_{-}(-iy)\bigr)^2}{m_+(iy)+m_-(-iy)} dy
 =O(1)\quad\text{as}\quad \eta\downarrow 0.
\end{equation}
\end{enumerate}
\end{lemma}
\begin{proof}
\noindent
(\ref{Coupl})
Let us show that $\rho(A)\ne\emptyset$ for the operator $A$ from Theorem~\ref{tHsc}. Assume $\rho(A)=\emptyset$. Then, by Theorem~\ref{tHsc}(\ref{irho-i})
\begin{equation}\label{eq:m_pm_poles}
  m_+(z)+m_-(-z) = 0 \quad \text{for all} \quad z \in \rho(A_{+,0})\cap\rho(A_{-,0}).
\end{equation}
Since $B_\pm$ is nonnegative, its self-adjoint extension $B_{\pm,0} = \pm A_{\pm,0}$ has at most one isolated negative eigenvalue. Hence $m_\pm$ has at most one pole $a_\pm$ in $\dR_-$. Now, equality~\eqref{eq:m_pm_poles} implies that $m_+$ has at most one pole $-a_-$ in $\dR_+$ and possibly a pole at $0$. Therefore,
\[
m_+(z)=\frac{\sigma_-}{-a_--z}+\frac{\sigma_0}{-z}+\frac{\sigma_+}{a_+-z}
\]
for some $a_-,a_+ < 0$, $\sigma_0,\sigma_-,\sigma_+\ge 0$. Since $B_+$ is densely defined, see (A\ref{assume3}), the last displayed formula contradicts \eqref{6.6}.

The operator  $ A_0=B_{+,0}\oplus (-B_{-,0})$ is  definitizable, and $\infty\in c(A_0)$ as, by assumptions, $B_{\pm,0}$ are unbounded.
Since the resolvent $(A-z)^{-1}$ of $A$ is a one-dimensional perturbation of  the resolvent $(A_0-z)^{-1}$, see~\eqref{ereswA1},
the claim (\ref{Coupl}) follows from~\cite[Theorem 1]{JonLan79}.

(\ref{iVeselic-1})
 Applying Theorem~\ref{tHsc} to the operators $A_+=B_+$ and $A_-=-B_-$  in the inner product spaces $(\cK_\pm,\kip_{\cK_\pm})=(\cH_\pm,\pm\ahip_{\cH_\pm})$  one obtains the equality
\begin{equation*} 
\bigl[\bigl(A - z \bigr)^{-1} f,f \bigr]_\cK = \bigl[\bigl( A_0 - z \bigr)^{-1} f,f\bigr]_\cK
  - \frac{\bigl[f,\gamma(\overline{z})\bigr]_\cK
  \bigl[\gamma({z}),f\bigr]_\cK}{m_+(z)+m_+(-z)},
 \end{equation*}
where $ A_0=B_{+,0}\oplus (-B_{-,0})$, $z \in \rho(A) \cap\rho(A_0)$. Since $A_0$ is a self-adjoint operator in $\cH$ Theorem~\ref{Veselic} and~\eqref{eq:Res_SA0} imply that the condition $\infty\not\in c_s(A)$ is equivalent to
\begin{equation}\label{eq:RegCrit1}
\int_{\eta_0}^{\eta}\re \frac{\bigl[f,\gamma(-iy)\bigr]_\cK
    \bigl[\gamma(iy),f\bigr]_\cK}{m_+(iy)+m_-(-iy)} dy=O(1), \quad
\eta\to +\infty\quad
\text{for all} \quad  f\in\cK.
\end{equation}
Decompose $f\in\cK$ into its ``real'' and ``imaginary'' part, $f=f^R+if^I$, where $f^R$ and   $f^I$ are real, see~\eqref{eq:Real_Decom}.
Since the vector valued functions $\gamma_{B_\pm}(z)$ are real, it follows from~\eqref{eq:Conj} and  \eqref{eq:3.10A} that
\[
\bigl[\gamma(iy),f^R\bigr]_\cK=\bigl[f^R,\gamma(-iy)\bigr]_\cK,\quad
\bigl[\gamma(iy),f^I\bigr]_\cK=\bigl[f^I,\gamma(-iy)\bigr]_\cK
\]
and hence (see also analogous identity in~\cite[Proof of Theorem~4.5]{Kost13})
\[
\begin{split}
\bigl[f,\gamma(-iy)\bigr]_\cK  \bigl[\gamma(iy),f\bigr]_\cK
  &=\bigl[f^R+if^I,\gamma(-iy)\bigr]_\cK  \bigl[\gamma(iy),f^R+if^I\bigr]_\cK  \\
  &=\bigl[f^R,\gamma(-iy)\bigr]_\cK^2+\bigl[f^I,\gamma(-iy)\bigr]_\cK^2.
\end{split}
\]
Therefore,  $\infty\not\in c_s(A)$ if and only if ~\eqref{eq:RegCrit1} hold for all real $f\in\cK$. Since
\begin{equation} \label{eq:cG}
\bigl[f,\gamma(-iy)\bigr]_\cK = \bigl\langle f_+,\gamma_{B_+}(-iy) \bigr\rangle_{\cH_+}
 - \bigl\langle f_-,\gamma_{B_-}(iy) \bigr\rangle_{\cH_-} = \wh f_+(iy) - \wh f_-(-iy),
\end{equation}
condition~\eqref{eq:RegCrit1} takes the form
\[
    \int_{\eta_0}^{\eta}\re \frac{ \bigl(\wh f_{+}(iy)-\wh f_{-}(-iy)\bigr)^2}{m_+(iy)+m_-(-iy)} dy = O(1)\quad\text{as}\quad \eta\to+\infty,
\]
which reduces to~\eqref{eq:Crit_Inf1} when we substitute $f_-$ with $- f_-$.

(\ref{iVeselic-2})
By Theorem~\ref{Veselic}  and~\eqref{eq:Res_SA0} we have that the conjunction $0\not\in c_s(A)$ and $\ker A=\ker A^2$
is equivalent to
\begin{equation}\label{eq:RegCrit10}
\int_{\eta}^{\eta_0}\re \frac{\bigl[f,\gamma(-iy)\bigr]_\cK
    \bigl[\gamma(iy),f\bigr]_\cK}{m_+(iy)+m_-(-iy)} d{y}=O(1) \quad \text{as} \quad
\eta\downarrow 0 \quad  \text{for all} \quad f\in\cK.
\end{equation}
The reasoning in the proof of item (\ref{iVeselic-1}) shows that the preceding equivalence is preserved if $f$ in~\eqref{eq:RegCrit10} is restricted to be real, which in view of~\eqref{eq:cG} yields \eqref{eq:Crit_Inf2}.
\end{proof}

\begin{remark}
  The coupling $A$ in Lemma~\ref{thm:Reg_infty}  is not necessary non-negative, so it may have other  critical points $\alpha$ distinct from $0$. A criterion for $\alpha\not\in c_s(A)$ can be derived from Lemma~\ref{thm:Reg_infty} (\ref{iVeselic-2}) by using the shift $A\rightarrow A-\alpha I$. However, we
   do not pursue this case here, since the main object of this paper, the operator $A$ in Section~5 associated with the indefinite Sturm-Liouville expression \eqref{eq2ndq69A} is  non-negative.
\end{remark}

\subsection{$D$-properties and conditions for regularity} 
\begin{definition}
A pair of Nevanlinna functions $m_+$ and $m_-$ is said to have the $D_\infty$-{\em property} (resp.~$D_0$-{\em property}) if
\begin{equation} \label{eqDinf}
\frac{\max\bigl\{\im m_+(iy),\im m_-(iy)
\bigr\}}{\bigl|m_+(iy)+m_-(-iy)\bigr|}=O(1)\quad
\text{as} \quad y\to+\infty \quad(\text{resp.} \ y\downarrow 0).
\end{equation}
\end{definition}

\begin{lemma}\label{lem:3.11}
Assume that a pair $m_+$ and $m_-$ has the $D_\infty$-property (resp.\ $D_0$-property). Then
\begin{equation}\label{eq:ratio2}
  \frac{m_+(iy)+{m_-(iy)}}{m_+(iy)+\overline{m_-(iy)}}=O(1)\quad
\text{as}\quad y\to+\infty \quad(\text{resp.}\quad y\downarrow 0).
\end{equation}
If, in addition, there exists $y_0>0$ such that
\begin{equation}\label{eq:Rem_pm}
     \re m_+(iy) \re m_-(iy) > 0 \quad \text{for all} \quad y > y_0 \quad (\text{resp.} \quad 0<y<y_0),
\end{equation}
then
\begin{equation}\label{eq:ratio1}
  \frac{m_+(iy)-\overline{m_-(iy)}}{m_+(iy)+\overline{m_-(iy)}}=O(1)\quad
\text{as}\quad y\to+\infty \quad(\text{resp.}\quad y\downarrow 0).
\end{equation}
\end{lemma}
\begin{proof}
Assume that the pair $m_+$ and $m_-$ has the $D_\infty$-property.  To prove \eqref{eq:ratio2} we use the notation $u_\pm(iy)$ and $v_\pm(iy)$ introduced in~\eqref{eq:ReIm_m}. With this notation we have
\begin{equation*}
  \bigl|m_+(iy)+m_-(iy)\bigr|^2
  = \bigl(u_+(iy)+u_-(iy)\bigr)^2 + \bigl(v_+(iy)+v_-(iy) \bigr)^2
\end{equation*}
and
\begin{equation*}
\frac{|u_+(iy)+u_-(iy)|}{|m_+(iy)+\overline{m_-(iy)}|}<  1\quad\text{for all} \quad y>0 .
\end{equation*}
By the $D_\infty$-property
\[
\frac{v_+(iy)+v_-(iy)}{|m_+(iy)+\overline{m_-(iy)}|}=O(1)\quad\text{as} \quad y \to +\infty.
\]
Hence~\eqref{eq:ratio2} holds.

To prove~\eqref{eq:ratio1}, assume further that there exists $y_0 > 0$ such that \eqref{eq:Rem_pm} holds for all $y > y_0$. Then~\eqref{eq:Rem_pm} yields
\begin{equation*} 
\frac{|\re m_\pm(iy)|}{|m_+(iy)+\overline{m_-(iy)}|}< \frac{|\re m_\pm(iy)|}{|\re m_+(iy)+ \re  m_-(iy)|} \le 1 \quad \text{for all} \quad y > y_0,
\end{equation*}
which, together with~\eqref{eqDinf}, imply~\eqref{eq:ratio1}.

To prove the claims involving the $D_0$-property we notice that if the pair $m_+(z)$ and $m_-(z)$ has the $D_0$-property, then the pair $m_+(-1/z)$ and $m_-(-1/z)$ has the $D_\infty$-property and we apply already proven statements to the functions $m_+(-1/z)$ and $m_-(-1/z)$.
\end{proof}

\begin{remark}
  As was shown in 
\cite{KaKost08}  the $D_\infty$-property ($D_0$-property, respectively) is necessary for the condition $\infty\not\in c_s( A)$
 ($0\not\in c_s( A)$, respectively).
A weaker form of condition \eqref{eqDinf} for the Sturm-Liouville operator \eqref{eq2ndq69A} with $w(x)=\sgn x$, and $r\equiv 1$ was presented in~\cite{KaMal07}.
\end{remark}

In the next theorem we show that the $D_\infty$-property becomes also sufficient for  $\infty\not\in c_s(A)$ if it is supplemented by the assumption \eqref{eq:Rem_pm}.

\begin{theorem}\label{prop:Semibound2}
Let conditions {\rm (A\ref{assume1})} through {\rm (A\ref{assumelast})} be satisfied and assume that there exists $y_0>0$ such that
\begin{equation}\label{eq:Rem_pm2}
     \re m_+(iy)\re m_-(iy)>0 \quad \text{for all} \quad y>y_0.
\end{equation}
Then the coupling $A$ is definitizable in the Krein
space $\bigl( \cK, \kip_{\cK} \bigr)$, $\infty\in c(A)$ and the following equivalence holds:
\begin{equation*}
    \infty\in c_r(A) \quad \Leftrightarrow
\quad \text{the pair $m_+$ and $m_-$ has the $D_\infty$-property.}
\end{equation*}
\end{theorem}
\begin{proof}
The definitizability of $A$ and $\infty\in c(A)$ follow from item (\ref{Coupl}) in Lemma~\ref{thm:Reg_infty}.

The necessity of the condition that the pair $m_+$ and $m_-$ has the $D_\infty$-property for  $\infty\not\in c_s(A)$ was proved in~
\cite{KaKost08}.

To prove sufficiency, assume that the pair $m_+$ and $m_-$ has the $D_\infty$-property. We use Lemma~\ref{thm:Reg_infty}(\ref{iVeselic-1}) to prove $\infty\in c_r(A)$.
The integral in \eqref{eq:Crit_Inf1} can be rewritten as the  sum $I_1(f_+,f_-)+I_2(f_+,f_-)$ of two integrals
\begin{align*}
 I_1(f_+,f_-) & = \int_{\eta_0}^\eta
 \frac{\re \bigl( (\wh f_+(iy)+\wh f_-(-iy))^2 \bigr) \bigl(u_+(iy)+ u_-(iy)\bigr)}%
    {\left|m_+(iy)+\overline{m_-(iy)}\right|^2}dy, \\
  I_2(f_+,f_-) & = \int_{\eta_0}^\eta
  \frac{\im \bigl( (\wh f_+(iy)+\wh f_-(-iy))^2\bigr)\bigl( v_+(iy) - v_-(iy) \bigr)}%
    {\left|m_+(iy)+\overline{m_-(iy)}\right|^2}dy,
\end{align*}
 where we use the notation introduced in \eqref{eq:ReIm_m}.

 We will prove that both of these integrals are bounded as $\eta\to+\infty$. By Lemma~\ref{lem:3.11} there exist $y_1,C_1>0$ such that
\begin{equation}\label{eq:y1C1}
  \left|\frac{m_+(iy)+ m_-(iy)}{m_+(iy)+\overline{m_-(iy)}}\right|\le C_1\quad\text{for all}\quad y>y_1.
\end{equation}
It follows from~\eqref{eq:y1C1} and \eqref{eq:Ineq2} that for all $\eta \geq \eta_0> y_1$ we have
\begin{align*}
\int_{\eta_0}^\eta\left|\wh f_\pm(iy)\right|^2  \frac{|u_+(iy)+u_-(iy)|}
    {\left|m_+(iy)+\overline{m_-(iy)}\right|^2}dy
     & \le C_1^2
    \int_{\eta_0}^\eta\left|\wh
    f_\pm(iy)\right|^2\frac{|u_+(iy)+u_-(iy)|}
    {\left|m_+(iy)+m_-(iy)\right|^2}dy \\
    &< 5\pi C_1^2\|f_\pm\|^2.
\end{align*}
This proves that
\begin{equation}\label{eq:IneqI1}
\bigl| I_1(f_+,f_-) \bigr| <10\pi C_1^2 \bigl( \|f_+\|_{\cH_+}^2+\|f_-\|_{\cH_-}^2 \bigr)
\end{equation}
for all real $f_\pm\in\cH_\pm$ and for all $\eta \geq \eta_0 > y_1$.

It follows from  \eqref{eq:ratio1}  that there exist $y_2, C_2>0$ such that
\begin{equation*} 
  \left|\frac{m_+(iy)-\overline{m_-(iy)}}{m_+(iy)+\overline{m_-(iy)}}\right|\le C_2
 \quad\text{ for all}\quad y>y_2.
\end{equation*}
 The preceding inequality yields that for all real $f_\pm\in\cH_\pm$ and  for all $\eta \geq \eta_0>y_2$  we have
\begin{align} \nonumber
\bigl|I_2(f_+,f_-)\bigr|
& \le C_2^2 \int_{\eta_0}^\eta
 \frac{\left| \im \bigl( (\wh f_+(iy)+\wh f_-(-iy))^2 \bigr)\right| |v_+(iy)-v_-(iy)|}
    {\bigl|m_+(iy)-\overline{m_-(iy)}\bigr|^2}dy \\
    \label{eq:Ineq5}
 & \le C_2^2 \int_{\eta_0}^\eta\frac{\left|\im \bigl((\wh f_+(iy)+\wh f_-(-iy))^2\bigr)\right|}
    {\im\bigl( m_+(iy)-\overline{m_-(iy)} \bigr)}dy,
\end{align}
 where for the second inequality we used that
\[
\left|m_+(iy)-\overline{m_-(iy)}\right|\ge v_+(iy)+v_-(iy)\ge |v_+(iy)-v_-(iy)|\quad\text{for all}\quad y>0.
\]

Next we prove the inequality
\begin{equation} \label{eq:ref_I2a}
\int_{\eta_0}^\eta
 \frac{\bigl|\im \bigl((\wh f_+(iy)+\wh f_-(-iy))^2\bigr)\bigr|}%
    {\im \bigl( m_+(iy)-\overline{m_-(iy)} \bigr)}dy
 \le 2\pi \bigl(\|f_+\|_{\cH_+}^2+\|f_-\|_{\cH_-}^2\bigr)
\end{equation}
for all real $f_\pm\in\cH_\pm$ and for all $\eta > \eta_0 > y_2$. This inequality will be deduced from  \eqref{eq:Ineq1} in Lemma~\ref{lem:Kost1} applied to the symmetric operator $S$ defined in Theorem~\ref{tHsc} with the following specific setting.

The Krein space $\bigl(\cK_+, \kip_{\cK_+} \bigr)$ is the Hilbert spaces $\bigl(\cH_+, \ahip_{\cH_+} \bigr)$, the Krein space $\bigl(\cK_-, \kip_{\cK_-} \bigr)$ is the Hilbert spaces $\bigl(\cH_-, \ahip_{\cH_-} \bigr)$, the symmetric operator $A_+$ is given by $A_+ = B_+$ and the symmetric operator $A_-$ is given by $A_- = - B_-$. Notice that in this part of the proof the spaces $\cK_+$ and $\cK_-$ differ from those in (A\ref{assume1}). The operator $A_\pm$ is a closed symmetric densely defined operator with defect numbers $(1,1)$ in the Hilbert space $\bigl(\cK_\pm, \kip_{\cK_\pm} \bigr)$. Furthermore, $\bigl(\mathbb{C},\Gamma_0^+,\Gamma_1^+\bigr)$ is a boundary triple for $A_+^{\langle *\rangle}$ (the adjoint in the Hilbert space $\cK_+ =\cH_+$) with the corresponding Weyl function and the $\gamma$-field given by
\[
z \mapsto m_+(z),  \quad z \mapsto \gamma_{B_+}(z), \quad z \in \mathbb{C}\setminus\mathbb{R},
\]
while $\bigl(\mathbb{C},\Gamma_0^-,-\Gamma_1^-\bigr)$ is a boundary triple for $A_-^{\langle *\rangle}$ (the adjoint in the Hilbert space $\cK_-=\cH_-$) with the corresponding Weyl function and the $\gamma$-field given by
\[
z \mapsto -m_-(-z),  \quad z \mapsto \gamma_{B_-}(-z), \quad z \in \mathbb{C}\setminus\mathbb{R}.
\]
The operator $S$ defined in Theorem~\ref{tHsc} on the domain~\eqref{eq:domSHC} is a real densely defined symmetric operator with defect numbers $(1,1)$ acting in a Hilbert space $\cK,\kip_{\cK}$ which is the direct sum, $\cK = \cH_+\oplus \cH_-$, of the Hilbert spaces $\bigl(\cH_+, \ahip_{\cH_+} \bigr)$ and $\bigl(\cH_-, \ahip_{\cH_-} \bigr)$.  A real boundary triple for $S^{\langle *\rangle}$ is the boundary triple  $\bigl(\mathbb{C},\Gamma_0,\Gamma_1\bigr)$ given in~\eqref{Htripl}. The corresponding Weyl function $M$ and the $\gamma$-field $\gamma$ are given by
\begin{equation*}
M(z)= m_+(z)-m_-(-z), \quad
  \gamma(z) = \begin{pmatrix}
                   \gamma_{B_+}(z) \\
                   \gamma_{B_-}(-z)
                 \end{pmatrix},\quad z \in \mathbb{C} \setminus \mathbb{R},
\end{equation*}
see~\eqref{mH2}.
 For $f=\begin{pmatrix}
                   f_+ \\ f_-
                 \end{pmatrix}\in \cK = \cH_+\oplus \cH_-$, from \eqref{eq:fpm}, we obtain
\[
\wh f(iy) = \bigl[ f, \gamma(-iy) \bigr]_{\cK}
 = \bigl\langle f_+,\gamma_{B_+}(-iy) \bigr\rangle_{\cH_+}
 + \bigl\langle f_-,\gamma_{B_-}(iy) \bigr\rangle_{\cH_-}
 = \wh f_+(iy)+\wh f_-(-iy).
\]
Now \eqref{eq:ref_I2a} follows by applying \eqref{eq:Ineq1} in Lemma~\ref{lem:Kost1} to the real symmetric operator $S$ acting in the Hilbert space $\cK$. Inequalities \eqref{eq:Ineq5} and \eqref{eq:ref_I2a} yield
\begin{equation} \label{eq:ref_I2}
\bigl|I_2(f_+,f_-)\bigr| \le C_2^2 2\pi \bigl(\|f_+\|_{\cH_+}^2+\|f_-\|_{\cH_-}^2\bigr)
\end{equation}
for all real $f_\pm\in\cH_\pm$ and for all $\eta > \eta_0 > y_2$. From \eqref{eq:IneqI1} and \eqref{eq:ref_I2} it follows that \eqref{eq:Crit_Inf1} holds. Hence Lemma~\ref{thm:Reg_infty}(\ref{iVeselic-1}) implies $\infty\in c_r(A)$.
\end{proof}

In the next theorem we give a criterion  for  $0\not\in c_s( A)$ formulated in terms of the $D_0$-property.
\begin{theorem} \label{thm:a_cs}
Let conditions {\rm (A\ref{assume1})} through {\rm (A\ref{assumelast})} be satisfied and assume that there exists $y_0>0$ such that the Weyl functions $m_+$ and $m_-$ satisfy the condition \begin{equation*}
     \re m_+(iy)\re m_-(iy)>0 \quad \text{for all} \quad  0 < y < y_0.
\end{equation*}
Then
\begin{equation*} 
0 \not\in c_s(A) \ \text{and} \ \ker A=\ker A^2
  \quad \Leftrightarrow \quad
  \text{the pair $m_+$ and $m_-$ has the $D_0$-property.}
\end{equation*}
\end{theorem}
\begin{proof}
The necessity of the $D_0$-property for $0\not\in c_s(A)$ was proved in 
\cite{KaKost08}.

To prove the sufficiency we will employ Lemma~\ref{thm:Reg_infty} and decompose the integral in
  \eqref{eq:Crit_Inf1} into a sum $I_1(f_+,f_-)+I_2(f_+,f_-)$ of two integrals
\begin{align*}
 I_1(f_+,f_-)
 & = \int_\eta^{\eta_0}
 \frac{\re \bigl( (\wh f_+(iy)+\wh f_-(-iy))^2 \bigr) \bigl(u_+(iy) + u_-(iy) \bigr)}%
    {\left|m_+(iy)+\overline{m_-(iy)}\right|^2}dy, \\
 I_2(f_+,f_-)
 & = \int_\eta^{\eta_0}
 \frac{\im \bigl( (\wh f_+(iy)+\wh f_-(-iy) )^2\bigr)\bigl(v_+(iy)- v_-(iy) \bigr)}%
    {\left|m_+(iy)+\overline{m_-(iy)}\right|^2}dy
\end{align*}
The estimates for $ I_1(f_+,f_-)$ and $ I_2(f_+,f_-)$ for every $f_\pm\in\cH_\pm$ similar to those in~\eqref{eq:IneqI1} and \eqref{eq:ref_I2} follow in the same way as  in the proof of Theorem~\ref{prop:Semibound2}.
\end{proof}

\begin{remark}
Notice that the condition~\eqref{eq:Rem_pm2} is not necessary for the non-negativity of the coupling $A$. For example,
let $B_\pm$ be minimal operators generated by the differential expression $-\frac{d^2}{dx^2}$ in $L^2(\dR_\pm)$, let $d_\pm\in\dR_\pm$ be such that  $d_++d_->0$, and let boundary triples $(\dC,\Gamma_0^\pm,\Gamma_1^\pm)$ for $B_\pm^*$ be given by
\[
\Gamma_0^+f_+=f'_+(0),\quad\Gamma_1^+f_+=-f_+(0)+d_+f'_+(0),\quad f_+\in\dom(B_+^*),
\]
\[
\Gamma_0^-f_-=f'_-(0),\quad\Gamma_1^-f_-=f_-(0)+d_-f'_-(0),\quad f_-\in\dom(B_-^*).
\]
Then the operator $A$ defined by~\eqref{eq:domwAHC} as the restriction of $-\sgn x\frac{d^2}{dx^2}$ to the domain
\[
\dom(A) =\left\{ f \in \dom(B_+^{\langle*\rangle})\oplus\dom(B_-^{\langle*\rangle}):
\begin{array}{c}
  f'_+(0)=f'_-(0) \\
  f_+(0)-f_-(0)=(d_-+d_+)f'_+(0)
\end{array}
 \right\}
\]
is nonnegative in the Krein space $\bigl(L_{w}^2(\dR),\kip_{w} \bigr)$, where $w(x) = \sgn x, x \in \mathbb{R}$. Indeed, for $f = f_+\oplus f_-\in\dom A$ we obtain
\begin{align*}
[Af,f]_{w}
& = -\int_{-\infty}^0f_-''\overline{f_-}-\int_0^{+\infty}f_+''\overline{f_+} \\
&= -f'_-(0)\overline{f_-(0)}+f'_+(0)\overline{f_+(0)} +\int_{\dR}|f'|^2 \\
& =(d_-+d_+)|f'_+(0)|^2+\int_{\dR}|f'|^2\ge 0.
\end{align*}
The Weyl functions $m_\pm$ of the operators $B_\pm$ corresponding to the boundary triples $(\dC,\Gamma_0^\pm,\Gamma_1^\pm)$ have the form
\[
m_+(z)=\frac{i}{\sqrt{z}}+d_+,\quad m_-(z)=\frac{i}{\sqrt{z}}+d_-,\quad d_\pm\in\dR_\pm
\]
and hence there exists $y_0>0$ such that $\re m_+(iy)\re m_-(iy)<0$ for all $ y>y_0.$
\end{remark}

\subsection{One-sided sufficient conditions for regularity}

In the next theorem we give a one-sided condition
which is sufficient for   $\infty\not\in c_s(A)$.

\begin{theorem}\label{prop:Semibound2C}
Let conditions {\rm (A\ref{assume1})} through {\rm (A\ref{assumelast})} be satisfied and assume that:
\begin{enumerate}
\renewcommand*\theenumi{\roman{enumi}}
\renewcommand*\labelenumi{{\rm (\theenumi)}}
\item \label{prop:Semibound2C-i-1}
there exists $y_0>0$ such that~\eqref{eq:Rem_pm2} holds for all $ y>y_0$;

\item \label{prop:Semibound2C-i-2}
either $ \im m_+(iy)\!=\!O(\re m_+(iy))$\! or $\im m_-(iy)=O(\re m_-(iy))$\! as $y\to +\infty$.
\end{enumerate}
Then the coupling $A$ of $A_+$ and $A_-$ is definitizable in the Krein space $\bigl( \cK, \kip_{\cK} \bigr)$,  $\infty\in c(A)$ and
\[
    \infty\in c_r(A).
\]
\end{theorem}
\begin{proof}
The definitizability of $A$ and $\infty\in c(A)$ follow from item (\ref{Coupl}) in Lemma~\ref{thm:Reg_infty}.
Let us assume that $ \im m_+(iy)\!=\!O(\re m_+(iy))$ as $y\to +\infty$ and show that the pair
$m_+$, $m_-$ has the $D_\infty$-property. Indeed, in view of (\ref{prop:Semibound2C-i-1})
\[
\left|\frac{\im m_+(iy)}{m_+(iy)+\overline{m_-(iy)}}\right|\le \frac{\im m_+(iy)}{\left|\re m_+(iy)\right|}
\]
and by assumption (\ref{prop:Semibound2C-i-2}) there exist $C,y_0>0$, such that
\begin{equation}\label{eq:mpm_2}
\left|\frac{\im m_+(iy)}{m_+(iy)+\overline{m_-(iy)}}\right|\le C,\quad y>y_0.
\end{equation}

Next, if ${\im m_-(iy)} > 2 \im m_+(iy)$, then
\[
|{\im m_-(iy)}-{\im m_+(iy)}|\ge|{\im m_-(iy)}|-|{\im m_+(iy)}|>\frac12 |{\im m_-(iy)}|
\]
and hence
\begin{equation*}
\left|\frac{\im m_-(iy)}{m_+(iy)+\overline{m_-(iy)}}\right|\le
2,\quad y>y_0.
\end{equation*}
Now,  if ${\im m_-(iy)} \le 2 \im m_+(iy)$, then
\begin{equation*}
\left|\frac{\im m_-(iy)}{m_+(iy)+\overline{m_-(iy)}}\right|\le
2\left|\frac{\im m_+(iy)}{m_+(iy)+\overline{m_-(iy)}}\right|\le 2C, \quad y > y_0.
\end{equation*}
Thus  the pair $m_+$, $m_-$ has the $D_\infty$-property, and the statement of Theorem~\ref{prop:Semibound2C} follows from Theorem~\ref{prop:Semibound2}.
\end{proof}

In the next theorem we formulate a one-sided condition which is sufficient for $0\not\in c_s(A)$.

\begin{theorem} \label{thm:a_cs02}
Let the conditions {\rm (A\ref{assume1})} through {\rm (A\ref{assumelast})} be satisfied and assume that
\begin{enumerate}
\renewcommand*\theenumi{\roman{enumi}}
\renewcommand*\labelenumi{{\rm (\theenumi)}}
  \item \label{thm:a_cs02-i-1}
there exist $y_0>0$ such that~\eqref{eq:Rem_pm} holds;
  \item \label{thm:a_cs02-i-2}
either $ \im m_+(iy)=O\bigl(\re m_+(iy)\bigr)$ or $\im m_-(iy)=O\bigl(\re m_-(iy)\bigr)$  as $y\downarrow 0$.
\end{enumerate}
Then $  0\not\in c_s(A)$ and $\ker A=\ker A^2.$
\end{theorem}
\begin{proof}
  Let us assume that $\im m_+(iy)=O(\re m_+(iy))$ as $y\downarrow 0$.
Then in view of (\ref{thm:a_cs02-i-1}) and (\ref{thm:a_cs02-i-2}) the inequality~\eqref{eq:mpm_2} holds for $0<y<y_0$
and, hence,  $\im m_+(iy)= O\bigl(m_+(iy)+\overline{m_-(iy)}\bigr)$ as $y\downarrow 0$.

The proof of the relation $\im m_-(iy)=O\bigl(m_+(iy)+\overline{m_-(iy)}\bigr)$ as $y\downarrow 0$ is similar to that in Theorem~\ref{prop:Semibound2C}.
Therefore, the pair $m_+$, $m_-$ has the $D_0$-property, and the statement of Theorem~\ref{thm:a_cs02} follows from Theorem~\ref{thm:a_cs}.
\end{proof}

\section{Sturm-Liouville operator with indefinite weight} \label{Sec:S-L}

\subsection{Indefinite Sturm-Liouville operator as a coupling} 

Let $I=(b_-,b_+)$ be a finite or infinite interval such that $-\infty \leq b_- <0 < b_+ \leq + \infty$ and let $\mathfrak{a}$ be the differential expression \eqref{eq2ndq69A} subject to the assumptions~\eqref{eq:rw}. In this section we study a nonnegative self-adjoint operator $A$ associated with $\mathfrak{a}$ in the Krein space $\bigl(L_{w}^2(I), \kip_w \bigr)$. In the definition of $A$ given in  \eqref{eq:A} we use nonnegative symmetric operators $B_\pm$  generated by the differential expressions  $\mathfrak{b}_{\pm}$ in the Hilbert spaces
$L^2_{w_\pm}(I_\pm)$ 
with the inner products
\[
\langle f,g\rangle_{w_\pm}=\int_{I_\pm} f(x)\overline{g(x)} w_{\pm}(x) dx.
\]
Let $B_{\pm,\max}$ be the maximal differential operator generated
in $L_{w_\pm}^2(I_\pm)$ by the differential expression $\mathfrak{b}_\pm$
(see \eqref{eq:Bpm}),
with the domain
\begin{equation*}
  \dom\bigl(B_{\pm,\max}\bigr) = \bigl\{ f \in L_{w_\pm}^2(I_\pm): f, f^{[1]} \in
AC_{\loc}(I_\pm), \ \mathfrak{b}_\pm(f) \in L_{w_\pm}^2(I_\pm) \bigr\},
\end{equation*}
where $f^{[1]}(x):=r(x)^{-1}f'(x), x\in I$. Let $B_{\pm,\min}( = (B_{\pm,\max})^{\langle*\rangle})$ be the minimal differential operator generated by $\mathfrak{b}_\pm$
in $L_{w_\pm}^2(I_\pm)$.

Let $z \in \mathbb{C}\setminus\mathbb{R}$ and denote by ${s}_\pm(\cdot,z)$ and ${c}_\pm(\cdot,z)$ the solutions on $I_\pm$ of the equation
\begin{equation} \label{eq:2.13a}
  \mathfrak{b}_\pm(f)=z f,
\end{equation}
satisfying the boundary conditions
 \[
{c}_\pm(0,z) = 1, \ c_\pm^{[1]}(0,z) = 0  \quad \text{and} \quad
{s}_\pm(0,z) = 0,
\ s_\pm^{[1]}(0,z) = 1.
 \]

If $\mathfrak{b}_\pm$ is in the limit point case at $b_\pm$ then neither ${s}_\pm(\cdot,z)$ nor ${c}_\pm(\cdot,z)$ belongs to  $L_{w_\pm}^2(I_\pm)$, however there exists a coefficient $m_\pm(z)$ such that the solution
\begin{equation} \label{eTWeyl}
\psi_\pm(t,z) = {s}_\pm(t,z) \mp m_\pm(z) {c}_\pm(t,z), \quad t
\in I_\pm,
\end{equation}
of the equation \eqref{eq:2.13a}  belongs to $L_{w_\pm}^2(I_\pm)$.

In the limit point case the  operator  $B_\pm:=B_{\pm,\min}$ is a symmetric operator in $L_{w_\pm}^2(I_\pm)$ with defect numbers $(1,1)$ and with the domain
\begin{equation}  \label{eq:domSpm*}
 \dom(B_\pm)  = \bigl\{ f \in \dom(B_{\pm,\max}) : f(0) = f^{[1]}(0) = 0\bigr\}.
 \end{equation}

In the limit circle case, by~\cite[Section 10.7]{KaKr74B}, for every
$f\in\dom(B_{\pm,\max})$ the following one-sided limit exists
\[
f^{[1]}(b_\pm):=\lim_{x\to b_\pm\mp 0}r_\pm(x)^{-1}f'(x).
\]
Let $m_\pm(z)$ be a coefficient  such that the solution $\psi_\pm(x,z)$ in~\eqref{eTWeyl} satisfies the condition
\begin{equation}\label{eq:m_function}
\psi_\pm^{[1]}(b_\pm,z)=0 \quad \text{for all} \quad z \in \mathbb{C}\!\setminus\!\mathbb{R}.
\end{equation}
Clearly, $m_\pm(z)$ is calculated as
$m_\pm(z)=\pm s^{[1]}(b_\pm,z)/c^{[1]}(b_\pm,z)$.
In the limit circle case the  operator  $B_{\pm,\min}$ is a symmetric operator in
$L_{w_\pm}^2(I_\pm)$ with defect numbers $(2,2)$ and we define its symmetric extension $B_\pm$ with defect numbers $(1,1)$
as the restriction of $\mathfrak{b}_\pm$ to the domain
\begin{equation}  \label{eq:domSpmLC}
 \dom(B_\pm)  = \bigl\{ f \in \dom(B_{\pm,\max}) : f(0) = f^{[1]}(0) =  f^{[1]}(b_\pm) =0\bigr\}.
 \end{equation}
 The adjoint operator $B_\pm^{\langle*\rangle}$ is the restriction of $\mathfrak{b}_\pm$ to the domain
\begin{equation*}
 \dom(B_\pm^{\langle*\rangle})  = \bigl\{ f \in \dom(B_{\pm,\max}) :  f^{[1]}(b_\pm) =0\bigr\}.
 \end{equation*}

In the following definition (see~\cite{Kost13}) the notion of  Neumann $m$-function is introduced both for the limit point case and the limit circle case.
\begin{definition}\label{def:5.1}
The function $m_\pm$ for which the solution $\psi_\pm(x,z)$ in~\eqref{eTWeyl} satisfies the condition
\begin{equation}\label{eq:m_functionLPLC}
  \left.\begin{array}{cl}
          \psi_\pm^{[1]}(b_{\pm},z)=0  & \text{if} \  \mathfrak{b}_\pm \ \text{is in the limit circle case at} \ b_\pm \\[6pt]
           \psi_\pm(\cdot,z)\in L_{w_\pm}^2(I_\pm) & \text{if} \ \mathfrak{b}_\pm \ \text{is in the limit point case at} \ b_\pm
            \end{array}\!\!\right\}
\end{equation}
is called the {\em  Neumann $m$-function of  $\mathfrak{b}_\pm$ on $I_\pm$ subject to~\eqref{eq:m_functionLPLC}}.
\end{definition}

The following proposition collects some facts from~\cite{DM17} about boundary triples for the operator $B_\pm^{\langle*\rangle}$.

\begin{proposition}\label{prop:5.1}
Assume that $\mathfrak{b}_\pm$ satisfies~\eqref{eq:rw}, let $B_\pm$
be defined as in~\eqref{eq:domSpm*} or in~\eqref{eq:domSpmLC}, respectively, (depending on limit point or limit circle case) and let $m_\pm$ be the {Neumann \em $m$-function} of $\mathfrak{b}_\pm$ on $I_\pm$, subject to~\eqref{eq:m_functionLPLC}.
Then:
  \begin{enumerate}
\renewcommand*\theenumi{\alph{enumi}}
\renewcommand*\labelenumi{{\rm (\theenumi)}}
    \item \label{prop:5.1i1}
$B_\pm$ is a symmetric nonnegative operator in the Hilbert space $L_{w_\pm}^2(I_\pm)$ with defect numbers $(1,1)$.
    \item \label{prop:5.1i2}
The triple $\bigl(\mathbb{C},\Gamma_0^\pm,\Gamma_1^\pm\bigr)$, where
\begin{equation}\label{eq:BTpm}
  \Gamma_0^\pm f_\pm=f_\pm^{[1]}(0),\quad \Gamma_1^\pm f_\pm=\mp f_\pm(0), \quad f \in \dom(B_\pm^{\langle*\rangle}),
\end{equation}
is a real boundary triple for $B_\pm^{\langle*\rangle}$.
\item \label{prop:5.1i3}
The Weyl function of $B_\pm$ corresponding to the boundary triple $\bigl(\mathbb{C},\Gamma_0^\pm,\Gamma_1^\pm\bigr)$ coincides with
the Neumann $m$-function $m_\pm$. That is
\begin{equation}\label{eq:WF_pm}
  m_\pm(z)=\mp \frac{\psi_\pm(0,z)}{\psi^{[1]}_\pm(0,z)}, \quad z\in \mathbb{C}\setminus\mathbb{R}.
\end{equation}
If $\mathfrak{b}_\pm$ is in the limit circle case at $b_\pm$, then, in addition to \eqref{eq:WF_pm}, the following formula holds
\begin{equation}\label{eq:WF_pmLC}
m_\pm(z) = \pm \frac{s_\pm^{[1]}(b_\pm,z)}{c_\pm^{[1]}(b_\pm,z)}, \quad z\in \mathbb{C}\setminus\mathbb{R}.
\end{equation}
\item \label{prop:5.1i4}
The Weyl function $m_\pm$ of $B_\pm$ belongs to the Stieltjes class $\cS$ and satisfies the condition
$\lim_{x\to -\infty}m_\pm(x)=0$. In particular,
    \begin{equation*}
  \re m_\pm(iy)\ge 0  \quad \text{for all } \quad y>0.
    \end{equation*}
  \end{enumerate}
\end{proposition}
\begin{proof}
Since
\[
\lim_{x \to b_{\pm}\mp 0}f^{[1]}(x)\overline{f(x)}=0 \quad
\text{for all} \quad f\in\dom(B_\pm^{\langle*\rangle})
\]
both in the limit point case~\cite[Corollary, p.\ 199]{Kal74} and in the limit circle case~\cite[Lemma 2.1]{EvEv91} the following formula holds
\begin{equation}\label{eq:GreenSL17}
  \int_{I_\pm}\mathfrak{b}_\pm(f_\pm)\overline{f_\pm} w_\pm dx = \pm f_\pm^{[1]}(0)\overline{f_\pm(0)}+\int_{I_\pm}\frac{1}{r_\pm}|f_\pm'|^2dx,\quad
  f_\pm\in\dom(B_\pm^{\langle*\rangle}).
\end{equation}
 By \eqref{eq:rw} and Definition~\ref{def:BTriple} this proves statements (\ref{prop:5.1i1}) and (\ref{prop:5.1i2}),
see also~\cite[Proposition 9.51, Theorem 9.69]{DM17}.

The statement (\ref{prop:5.1i3}) is implied by Definition~\ref{W00A} and the equalities
\[
\Gamma_0^\pm\psi_\pm(\cdot,z)=\psi_\pm^{[1]}(0,z) = 1, \quad
\Gamma_1^\pm\psi_\pm(\cdot,z)=\mp\psi_\pm(0,z)= m_\pm(z)\quad z\in \mathbb{C}\setminus\mathbb{R}.
\]
The formula~\eqref{eq:WF_pmLC} follows from \eqref{eq:m_function} and the equality
\[
0=\psi_\pm^{[1]}(b_\pm,z)=s_\pm^{[1]}(b_\pm,z)\mp m_\pm(z)c_\pm^{[1]}(b_\pm,z)\quad z\in \mathbb{C}\setminus\mathbb{R}.
\]

The extension $B_{\pm,0}$ of $B_\pm$ defined by
\begin{equation}\label{eq:B_pm0}
  B_{\pm,0}f=B_\pm^{\langle*\rangle}f,\quad f\in\dom (B_{\pm,0}):=\ker \Gamma_0^\pm
\end{equation}
is the von Neumann extension of $B_\pm$. Hence  $B_{\pm,0}\ge 0$, see also~\eqref{eq:GreenSL17}, and thus the function $m_\pm$ is holomorphic on $\mathbb{R}_-$. Moreover, as it follows from~\cite[Theorem~3.1]{KaKr74B}, see also~\cite[Proposition~3.6]{DSW20},
\begin{equation*}
  \lim_{x\to-\infty} m_\pm(x)=0
\end{equation*}
 and hence $m_\pm\in\cS$. This proves (\ref{prop:5.1i4}).
\end{proof}

With the differential expression $\mathfrak{a}$ we associate the following operator $A$ in the Krein $\bigl(L_{w}^2(I), \kip_w\bigr)$:
\begin{equation} \label{eq:domA}
\dom(A) = \bigl\{ f \in \dom(B_+^{\langle*\rangle})\oplus\dom(B_-^{\langle*\rangle}):
f, r^{-1}f' \in AC_{\loc}(I) \bigr\}
\end{equation}
and
\begin{equation}\label{eq:A}
Af = \mathfrak{a}(f), \qquad f \in \dom(A).
\end{equation}

\begin{lemma}\label{lem:Ker}
For every $\lambda\in\mathbb{R}$ the subspace $\ker (A-\lambda I)$  is at most one-dimensional.
\end{lemma}
\begin{proof}
Let $\lambda\in\mathbb{R}$. If $\mathfrak{b}_\pm$ is limit circle at $b_\pm$, then by Weyl's alternative the equation $\mathfrak{b}_\pm(f) = \pm \lambda f$ has two linearly independent solutions $c_\pm(x,\pm \lambda)$ and $s_\pm(x,\pm \lambda)$ in $L_{w_\pm}^2(I_\pm)$. Since the Wronskian of these solutions is not zero, it is not possible that both of these solutions satisfy $f^{[1]}(b_\pm) = 0$. Therefore, $\ker (B_\pm^{\langle *\rangle}\mp\lambda I)$ is one-dimensional.

If $\mathfrak{b}_\pm$ is limit point at $b_\pm$, then by Weyl's alternative the equation $\mathfrak{b}_\pm(f) = \pm \lambda f$ has at most one solution in $L_{w_\pm}^2(I_\pm)$. Consequently, $\ker (B_\pm^{\langle *\rangle}\mp\lambda I)$ is at most one-dimensional.

By the uniqueness theorem for linear initial value problems, the only solution of the problem $\mathfrak{b}_\pm(f) = \pm \lambda f$,  $f_\pm(0) = f_\pm^{[1]}(0) = 0$ is the zero function. Therefore, the subspace
\[
\ker (A-\lambda I) = \left\{ f =f_+\oplus f_-:\,
\begin{array}{cc}
  f_+\in \ker (B_+^{\langle *\rangle}-\lambda I), & f_+(0) = f_-(0)\\
  f_-\in \ker (B_-^{\langle *\rangle}+\lambda I), & f_+^{[1]}(0) = f_-^{[1]}(0)
\end{array} \right\}
\]
is also at most one-dimensional.
\end{proof}

\begin{theorem}\label{thm:SL}
Let the differential expression $\mathfrak{b}$ satisfy \eqref{eq:rw} and
let $m_\pm$ be the Neumann $m$-function of $\mathfrak{b}_\pm$ subject to~\eqref{eq:m_functionLPLC} on $I_\pm$.
Then the operator $A$ associated with the expression $\mathfrak{a}$ is the coupling of the operators  $A_+:=B_+$ and $A_-:=-B_-$ in the sense of Theorem~{\rm\ref{tHsc}}. The operator $A$ is a nonnegative self-adjoint operator in the Krein space $\bigl(L_{w}^2(I), \kip_w\bigr)$
with $\rho(A)\ne\emptyset$ and $\infty \in c(A)$. We have
\begin{enumerate}
 \renewcommand*\theenumi{\roman{enumi}}
\renewcommand*\labelenumi{\rm{(\theenumi)}}
\item \label{thm:SL-i-1}
$\infty \in c_r(A)\ \  \Leftrightarrow \ \
\text{the pair $m_+$ and $m_-$ has the $D_\infty$-property}.$
\item \label{thm:SL-i-2}
$ 0 \not\in c_s(A) \, \text{and} \, \ker A = \ker A^2
\, \Leftrightarrow
\, \text{the pair $m_+$ and $m_-$ has the $D_0$-property}.$
\item \label{thm:SL-i-3}
$\im m_+(iy)=O\bigl(\re m_+(iy)\bigr)$ as $y\to+\infty
\ \
  \Rightarrow
\ \
\infty \in c_r(A).$
\item  \label{thm:SL-i-4}
$\im m_-(iy)=O\bigl(\re m_-(iy)\bigr)$ as $y\to+\infty
\ \
  \Rightarrow
\ \
\infty \in c_r(A).$
\item \label{thm:SL-i-5}
$\im m_+(iy)=O\bigl(\re m_+(iy)\bigr)$ as $y\downarrow 0
\ \
  \Rightarrow
\ \
0\not\in c_s(A) \, \text{and} \, \ker A=\ker A^2$.
\item \label{thm:SL-i-6}
$\im m_-(iy)=O\bigl(\re m_-(iy)\bigr)$ as $y\downarrow 0
\ \
  \Rightarrow
\ \
0\not\in c_s(A)\, \text{and} \,\ker A=\ker A^2$.
\end{enumerate}
\end{theorem}
\begin{proof}
The boundary triples $\bigl(\mathbb{C},\Gamma_0^\pm,\Gamma_1^\pm\bigr)$ from Proposition~\ref{prop:5.1} are also boundary triples for $A_\pm^{+}$.
The coupling of the  operators $A_\pm$ in Theorem~\ref{tHsc} is characterized by the conditions
\[
 \Gamma_{0}^+(f_+) - \Gamma_{0}^- (f_-)=0,  \quad
 \Gamma_{1}^+(f_+) + \Gamma_{1}^- (f_-)=0, \quad f_\pm\in\dom(B_\pm^{\langle*\rangle})
\]
which in view of \eqref{eq:BTpm} can be rewritten as
\begin{equation}\label{eq:A_Coupl}
f^{[1]}_+(0)=f^{[1]}_-(0),\quad f_+(0)=f_-(0), \quad f_\pm\in\dom(B_\pm^{\langle*\rangle}).
\end{equation}
Therefore, the differential operator $A$ associated with the expression $\mathfrak{a}$ is the coupling of the operators
$A_\pm:=\pm B_\pm$ relative to the boundary triples  $\bigl(\mathbb{C},\Gamma_0^\pm,\Gamma_1^\pm\bigr)$.
It follows from~\eqref{eq:GreenSL17} and~\eqref{eq:A_Coupl} that for $f=f_++f_-\in\dom A$, $f_\pm\in\dom(B_\pm^{\langle*\rangle})$, we have
\begin{equation*}
\begin{split}
   [Af,f]_{w} & = \bigl\langle B_+^{\langle*\rangle}f_+,f_+ \bigr\rangle_{w_+}
   + \bigl\langle B_-^{\langle*\rangle}f_-,f_- \bigr\rangle_{w_-} \\
     & = f_+^{[1]}(0)\overline{f_+(0)}-f_-^{[1]}(0)\overline{f_-(0)}
     +\int_{I}\frac{1}{r}|f'|^2dt
     =\int_{I}\frac{1}{r}|f'|^2dt\ge 0.
\end{split}
\end{equation*}
Hence the operator $A$ is nonnegative.

The Weyl function $M_\pm$  of the operator $A_{\pm}$ corresponding to
$(\mathbb{C},\Gamma_0^\pm,\Gamma_1^\pm)$ and the Weyl function $m_\pm$ of the operator $B_{\pm}$ satisfy
\[
M_{\pm}(z) = m_{\pm}(\pm z), \qquad z \in \mathbb{C}\setminus\mathbb{R}.
\]
Since $m_+, m_- \in \cS$, by Proposition~\ref{p:SiSm} we have $\re\bigl(m_+(iy)+m_-(-iy)\bigr)>0$ for all $y\in\mathbb{R}_+$. Consequently,   Theorem~\ref{tHsc}(\ref{irho-i}) yields $\rho(A)\ne\emptyset$.
Therefore the operator $A$ is definitizable and $\infty\in c(A)$, see
Lemma~\ref{thm:Reg_infty}.
As $m_+,m_-\in\cS$, the assumptions of Theorems~\ref{prop:Semibound2},
\ref{thm:a_cs}, \ref{prop:Semibound2C}, and \ref{thm:a_cs02} are satisfied, and, thus, the remaining claims  follow.
\end{proof}

\begin{remark} \label{rBenn}
In the limit circle case Bennewitz, see~\cite{Ben89}, considered a more general class of Neumann $m$-functions than introduced in Definition~\ref{def:5.1}. We restate Bennewitz's definition here.

Denote the Wronskian of two functions $f,g\in\dom{B_{\pm,\max}}$ by
\begin{equation*}
{\mathsf W}_t(f,g) := f(t)g^{[1]}(t)-f^{[1]}(t)g(t), \quad t \in I_\pm.
\end{equation*}
The one-sided limit
\begin{equation*}
{\mathsf W}_{b_\pm}\!(f,g) := \lim_{x\to b_\pm \mp 0} \bigl(f(t)g^{[1]}(t)-f^{[1]}(t)g(t) \bigr)
\end{equation*}
exists for all $f,g\in\dom(B_{\pm,\max})$. Furthermore, according to Titchmarsh~\cite{Tit62} (see also~\cite[Theorem 9.69]{DM17}, every symmetric boundary condition at $b_\pm$ for arbitrary $f\in\dom(B_{\pm,\max})$ can be written as
\begin{equation*}
{\mathsf W}_{b_\pm}\!\bigl(f, (\cos\alpha) {s}_\pm(\cdot,z_0) + (\sin\alpha) {c}_\pm(\cdot,z_0) \bigr) = 0
\end{equation*}
for some $\alpha \in (-\pi/2,\pi/2]$ and some $z_0 \in \mathbb{C}\!\setminus\!\mathbb{R}$.

If $m_\pm(z)$ is a coefficient for which the solution $\psi_\pm(t,z)$ in~\eqref{eTWeyl} satisfies the condition
\begin{equation}\label{eq:m_func_Ben}
 {\mathsf W}_{b_\pm}\!\bigl(\psi_\pm(\cdot,z), (\cos\alpha) {s}_\pm(\cdot,z_0) + (\sin\alpha) {c}_\pm(\cdot,z_0) \bigr) = 0, \quad z \in \mathbb{C}\!\setminus\!\mathbb{R},
\end{equation}
for some $\alpha \in (-\pi/2,\pi/2]$, then $m_\pm$ is called the {\em  Neumann $m$-function of $\mathfrak{b}_\pm$ on $I_\pm$}.
Clearly, $m_\pm(z)$ can be expressed as
\begin{equation*}
m_\pm(z)=
\frac{%
(\cos\alpha){\mathsf W}_{b_\pm}\!\bigl(s_\pm(\cdot,z),s_\pm(\cdot,z_0)\bigr)+(\sin\alpha){\mathsf W}_{b_\pm}\!\bigl(s_\pm(\cdot,z),c_\pm(\cdot,z_0)\bigr)%
}
{%
(\cos\alpha){\mathsf W}_{b_\pm}\!\bigl(c_\pm(\cdot,z),s_\pm(\cdot,z_0)\bigr)+(\sin\alpha){\mathsf W}_{b_\pm}\!\bigl(c_\pm(\cdot,z),c_\pm(\cdot,z_0)\bigr)%
},
\  z \in \mathbb{C}\setminus\mathbb{R}.
\end{equation*}
Since all the symmetric boundary conditions at $b_\pm$ are included in the boundary condition~\eqref{eq:m_func_Ben}, the boundary condition~\eqref{eq:m_function} is included as well. Therefore, the class of Neumann $m$-functions introduced in this remark contains the Neumann $m$-functions introduced in Definition~\ref{def:5.1}.
\end{remark}

\subsection{Asymptotic properties of $m$-functions}\label{sec:5.2}

V.A.~Mar\v{c}enko~\cite{Mar52} (for Sturm-Liouville operator $-\frac{d^2}{dx^2}+q$), and I.S.~Kac \cite{Kac73} and Y.~Kasahara~\cite{Kas75} (for  weighted Sturm-Liuoville operator) showed that the asymptotic behaviour of the Weyl function $m$ along the imaginary axes at $+\infty$ is closely related to the behaviour of the coefficients of the differential expression at $0$. In this section we present some results in this direction from~\cite{Ben87,Ben89} and their recent developments in~\cite{Kost13}.

Recall the definition~\eqref{eq:WR} of functions $W_\pm$ and $R_\pm$:
\begin{equation}\label{eq:WR_pm}
    W_\pm( x):=\int_0^x w_\pm(\xi)d\xi,\quad R_\pm(x):=\int_0^x r_\pm (\xi)d\xi,\quad x\in I_\pm,
\end{equation}
where $W_+$ and $R_+$ are positive and increasing on $I_+$, while
$W_-$ and $R_-$ are negative and increasing functions on $I_-$. Define the function $F_\pm: R_\pm\bigl(I_\pm\bigr) \to \mathbb{R}_+$ as follows
\begin{equation}\label{eq:F_pm}
F_\pm(x) :=  \frac{1}{x W_\pm\bigl(R_\pm^{-1}(x)\bigr)}, \quad x \in R_\pm\bigl(I_\pm\bigr).
\end{equation}
Here
\begin{equation}\label{Graz}
R_-(I_-) = (c_-,0), \quad  R_+(I_+) = (0,c_+) \quad \text{with} \quad -\infty \leq c_- < 0 < c_+ \leq +\infty.
\end{equation}
The function $F_+$ is decreasing and unbounded, and $F_-$ is an unbounded increasing function.
Denote by $f_\pm$ the inverse of $F_\pm$. Notice that both $f_-$ and $f_+$ are defined in a neighbourhood of $+\infty$, the function $f_+$ is positive and decreasing, the function
$f_-$ is negative and increasing, and
\[
\lim_{x\to +\infty} f_-(x) = 0 \quad \text{and} \quad \lim_{x\to +\infty} f_+(x) = 0.
\]

The following result was proved by F. Atkinson \cite{Atk85}, see also Bennewitz \cite[Theorem~3.4]{Ben89} for an improved version which we use here. The concept of the Neumann $m$-function of $\mathfrak{b}_\pm$ on $I_\pm$ is used in the sense defined in Remark~\ref{rBenn}.  For the concept of a  slowly varying function at $0_\pm$ we refer to Definition~\ref{def:sv} in  Appendix~\ref{appen}.

\begin{theorem}\label{thm:Asym_m2}
Let $W_\pm$ and $R_\pm$ be the functions defined in \eqref{eq:WR_pm}, let $f_\pm$ be the inverse of the function defined in \eqref{eq:F_pm} and let $m_\pm$ be the Neumann $m$-function of $\mathfrak{b}_\pm$ on $I_\pm$. If $W_\pm\circ R_\pm^{-1}$ is a slowly varying function at $0_\pm$, then
\begin{equation*} 
   m_\pm(iy)\sim \pm i f_\pm(y)\quad \text{as}\quad y\to+\infty.
\end{equation*}
\end{theorem}
\begin{proof}
Assume that $W_\pm\circ R_\pm^{-1}$ is a slowly varying function at $0_\pm$.
By Corollary~\ref{cor:AtkCon} this condition is equivalent to
  \begin{equation}\label{eq:Sub_RW}
    \int_0^xR_\pm(\xi) \, d W_\pm(\xi)= o(R_\pm(x)W_\pm(x)) \quad \text{as} \quad |x| \downarrow 0 \quad \text{with} \quad x \in I_\pm.
  \end{equation}
The claim about the function $m_+$ was proved in~\cite[Theorem~3.4]{Ben89}.  We use this result to prove the claim about $m_-$.  Let us set $\wh w_+(x)=w_-(-x)$, $\wh r_+(x)=r_-(-x)$, $x\in\wh I_+=(0,-b_-)$. Then the Hilbert space $L^2_{\wh w_+}(\wh I_+)$ consists  of functions
\[
\wh y(x):=y(-x),\quad y\in L^2_{w_-}( I_- ).
\]
Let $\wh B_+$ be the minimal operator generated in $L^2_{\wh w_+}(\wh I_+)$ by the differential expression
\[
(\widehat{\mathfrak{b}}_+ (\wh{ f}))(x):=-(\mathfrak{b}_- (f))(-x),\quad x\in \wh I_+.
\]
Then the Neumann $m$-function $\wh m_+$ of $\widehat{\mathfrak{b}}_+$ on $\wh I_+$ is connected with $m_-$ by
\begin{equation}\label{eq:whm}
  \wh m_+(z)=-m_-(-z).
\end{equation}
Next the functions
\[
   \wh W_+( x):=\int_0^x \wh w_+(\xi)d\xi,\quad \wh R_+( x):=\int_0^x\wh r_+ (\xi)d\xi,\quad x\in \wh I_+
\]
are connected with $W_-$ and $R_-$ by the equalities
\begin{equation}\label{eq:wh_WR}
  \wh W_+( x)=-W_-(-x),\quad \wh R_+( x)=-R_-(-x),\quad x\in \wh I_+
\end{equation}
and the inverse $\wh f_+$ of $\wh F_+:(0,\epsilon) \ni x \mapsto 1/\bigl(x (\wh W_+\circ \wh R_+^{-1})(x) \bigr)\in\mathbb{R}_+$ is connected with
 the  inverse $f_-$ of $F_-: (-\epsilon, 0) \ni x \mapsto 1/\bigl(x (W_-\circ R_-^{-1})(-x)\bigr)\in\mathbb{R}_+$  by the equality
\begin{equation}\label{eq:wh_f}
  \wh f_+( y)=-f_-(y),\quad y\in\mathbb{R}_+.
\end{equation}
It is easy to see, that $\wh W_+$ and $\wh R_+$ satisfy the condition~\eqref{eq:Sub_RW}. Therefore, by Theorem~\ref{thm:Asym_m2}
$\wh m_+(iy)\sim i \wh f_+(y)$.
Hence one obtains by~\eqref{eq:whm}, \eqref{eq:wh_WR}, \eqref{eq:wh_f}
\begin{equation*}
m_-(iy)=-\overline{\wh m_+(iy)}\sim-\overline{i\wh f_+(y)}={i\wh f_+(y)}=-if_-(y). \qedhere
\end{equation*}
\end{proof}

The sufficiency  part of the following lemma was proved by Bennewitz~\cite{Ben87}. The condition that appears in \cite{Ben87} is equivalent to the definition of a positively increasing function, see Definition~\ref{CieloDiablo} in Appendix~\ref{appen}.
The necessity of condition~\eqref{eq:RI_infty} below was proved by Kostenko in~\cite{Kost14}.

\begin{lemma}\label{lem:I_infty}
Let $m_\pm$ be the Neumann $m$-function of $\mathfrak{b}_\pm$ on $I_\pm$.
Then
\begin{equation}\label{eq:RI_infty}
\re m_\pm(iy)=O(\im m_\pm(iy))\quad \text{as}\quad y\to\pm\infty
\end{equation}
if and only if the function $R_\pm \circ W_\pm^{-1}$ is positively increasing
at $0_\pm$.
\end{lemma}

Notice that the concept of the Neumann $m$-function of $\mathfrak{b}_\pm$ on $I_\pm$ in Lemma~\ref{lem:I_infty} is used in the sense defined in Remark~\ref{rBenn}, while in the rest of the paper we use Definition~\ref{def:5.1}.
The following analog of Lemma~\ref{lem:I_infty} was proved in~\cite[Corollary~2.7]{Kost13}.
\begin{lemma}\label{cor:R_inf}
Let $m_\pm$ be the  Neumann $m$-function  of $\mathfrak{b}_\pm$ on $I_\pm$, subject to~\eqref{eq:m_functionLPLC}.
Then
\begin{equation}\label{eq:IR_inftyA}
\im m_\pm(iy)=O(\re m_\pm(iy))\quad \text{as} \quad y\to \pm\infty
\end{equation}
if and only if the function $W_\pm\circ R_\pm^{-1}$ is positively increasing at $0_\pm$.
\end{lemma}

Similar criteria for estimates~\eqref{eq:RI_infty}
and~\eqref{eq:IR_inftyA} at 0 were proved by Kostenko
in~\cite[Theorem~2.11 and Corollary~2.15]{Kost13}.

\begin{lemma}
Let $w_\pm,r_\pm\not\in L^1(I_\pm)$ and let $m_\pm$ be the  Neumann $m$-function  of $\mathfrak{b}_\pm$ on $I_\pm$, subject to~\eqref{eq:m_functionLPLC}.
Then
\begin{equation*}
\re m_\pm(iy)=O(\im m_\pm(iy))\quad \text{as}\quad y \to 0_\pm
\end{equation*}
if and only if the function $R_\pm\circ W_\pm^{-1}$ is positively increasing
at $\pm\infty$.
\end{lemma}
\begin{lemma}\label{cor:R_0}
Let $w_\pm,r_\pm\not\in L^1(I_\pm)$ and let $m_\pm$ be the  Neumann $m$-function  of $\mathfrak{b}_\pm$ on $I_\pm$, subject to~\eqref{eq:m_functionLPLC}.
Then
\begin{equation*}
\im m_\pm(iy)=O(\re m_\pm(iy)) \quad \text{as} \quad y \to 0_\pm
\end{equation*}
if and only if the function $W_\pm\circ R_\pm^{-1}$ is positively increasing at $\pm\infty$.
\end{lemma}

In the following lemma we consider the cases in which the conditions $w_\pm,r_\pm\not\in L^1(I_\pm)$ are not satisfied.

 \begin{lemma}\label{lem:PolF}
Let $m_\pm$ be the {\em Neumann $m$-function} of $\mathfrak{b}_\pm$ on $I_\pm$, subject to~\eqref{eq:m_functionLPLC}.
    \begin{enumerate}
 \renewcommand*\theenumi{\roman{enumi}}
\renewcommand*\labelenumi{\rm{(\theenumi)}}
\item \label{lem:PolFi1}
Let $a_\pm =\pm \lim_{x\to b_\pm} 1/W_\pm(x)$. Then $a_\pm \geq 0$ and the function
        \begin{equation*}
\wt m_\pm({z}):=m_\pm({z})+\frac{a_\pm}{{z}}, \quad {z} \in \mathbb{C}_+,  \end{equation*}
belongs to $\cS$  and
$
       \lim_{y\downarrow 0}y\mkern 1mu \wt m_\pm(iy)=0.
$ 
In particular, if $w_\pm\in L^1(I_\pm)$, then $a_{\pm}>0$ ,
$y\mkern 1mu m_\pm(iy)\sim i {a_\pm}$ at $0_+$ and
\begin{equation*} 
\re m_\pm(iy)=o(\im m_\pm(iy))\quad \text{as}\quad y\downarrow 0.
\end{equation*}
\item \label{lem:PolFi2}
If $r_\pm \in L^1(I_\pm)$ and $w_\pm \not\in L^1(I_\pm)$, then
        \begin{equation}\label{eq:Pol_m2}
    \im m_\pm(iy)=o(\re m_\pm(iy))\quad \text{as}\quad
y\downarrow 0.
    \end{equation}
\end{enumerate}
\end{lemma}
\begin{proof}
The claims (\ref{lem:PolFi1}) and (\ref{lem:PolFi2}) appear in~\cite[Lemma 2.10]{Kost13}. For the proof of (\ref{lem:PolFi1}) see also~\cite[Propositions 3.6, 4.6]{DSW20}.
\end{proof}

\subsection{Regularity of the critical point $\infty$} 

Statements (\ref{thm:SL-i-3}), (\ref{thm:SL-i-4}) of Theorem~\ref{thm:SL} can be restated as follows.

\begin{theorem}\label{thm:Reg_infty_PIS}
Let the differential expression $\mathfrak{b}_\pm$  satisfy \eqref{eq:rw} and let  the functions $R_\pm$ and $W_\pm$ be defined by~\eqref{eq:WR_pm}.
If either $W_+\circ R_+^{-1}$ is positively increasing at $0_+$ or $W_-\circ R_-^{-1}$ is positively increasing at $0_-$, then $\infty\in c_r(A)$.
\end{theorem}
\begin{proof}
Let $m_\pm$ be the Neumann $m$-function of $\mathfrak{b}_\pm$ on $I_\pm$, subject to~\eqref{eq:m_functionLPLC}.
By Proposition~\ref{prop:5.1}(\ref{prop:5.1i4}) $m_+$ and $m_-$ belong to the Stieltjes class $\cS$.
Thus $m_+$ and $m_-$ satisfy the assumption~\eqref{eq:Rem_pm2} of Theorem~\ref{prop:Semibound2}.
Assume that $W_+\circ R_+^{-1}$ is positively
increasing at $0_+$. Then by Lemma~\ref{cor:R_inf} condition~\eqref{eq:IR_inftyA} holds. Hence by Theorem~\ref{thm:SL} (\ref{thm:SL-i-3}) we have  $\infty\in c_r(A)$. Similar argument proves the theorem if $W_-\circ R_-^{-1}$ is positively increasing at $0_-$.
\end{proof}
\begin{example}\label{ex:Beals}
Let $I=(-1,1)$. Consider differential operators $B_\pm$ generated by $\mathfrak{b}_\pm$ in $L^2(I_\pm)$, where $r_-$, $w_-$ are arbitrary subject to conditions \eqref{eq:rw} and $r_+ = 1$, and $w_+$ satisfies the condition:
  \begin{equation}\label{eq:Beals}
      w_+(x)=x^{\alpha} v_+(x),\quad x\in I_+,\quad   \alpha>-1,
  \end{equation}
where $v_+(x)$ is slowly varying at $0_+$. Then by Karamata's characterization theorem,  Theorem~\ref{th-Kon}, we have
\[
  W_+(x)=\int_0^x t^{\alpha} v_+(t)dt \sim \frac{x^{\alpha+1}}{\alpha+1}v_+(x)\quad
  \text{as}\quad x\downarrow 0,
\]
and hence $W_+(x)$ is regularly varying at $0_+$ of order $\alpha + 1>0$ by Proposition~\ref{pro:rvsv}. Theorem~\ref{thm:Reg_infty_PIS} yields that $\infty\in c_r(A)$.

In the case when both $w_+$ and $w_-$ satisfy the condition \eqref{eq:Beals} with $v_\pm\in C^1(\overline{I_\pm})$ and  $\alpha>-1/2$ (so called Beals conditions) this result was obtained by R. Beals in~\cite{Bea85},
and by B. \'{C}urgus and H. Langer in~\cite{CL89} for $\alpha>-1$. That one-sided condition for the weight $w$ on $I_+$ is enough for  $\infty\in c_r(A)$ was noticed by A.~Fleige in~\cite{Fl95}.
\end{example}

\begin{lemma} \label{lem:simbdd}
Let $a \in \mathbb{R}_+$ and let $\alpha, \beta, f, g: [a,+\infty) \to \mathbb{C}\!\setminus\!\{0\}$ be functions such that $\alpha$ and $\beta$ are bounded,
\begin{equation}\label{eq:simbddA}
\lim_{x\to+\infty} \frac{\alpha(x)}{\beta(x)} = 1 \quad \text{and} \quad \lim_{x\to+\infty} \frac{f(x)}{g(x)} = 1.
\end{equation}
Then
\begin{equation}\label{eq:simbdd}
\frac{1}{\alpha(x) - f(x)} = O(1) \ \text{as} \ x\to+\infty \quad \Leftrightarrow \quad \frac{1}{\beta(x) - g(x)} = O(1) \ \text{as} \ x\to+\infty.
\end{equation}
\end{lemma}
\begin{proof}
We will prove the equivalence of the negations of the statements in \eqref{eq:simbdd}. The negation of the statement on the left-hand side of \eqref{eq:simbdd} is: There exists an increasing sequence $(x_n)$ in $[a,+\infty)$ such that
\[
\lim_{n\to+\infty} x_n = +\infty \qquad \text{and} \qquad \lim_{n\to+\infty} \bigl(\alpha(x_n) - f(x_n)\bigr) = 0.
\]
Since for all $n\in \mathbb{N}$ we have
\begin{equation*}
\beta(x_n) - g(x_n) =\alpha(x_n) \left( \frac{\beta(x_n)}{\alpha(x_n)} - \frac{g(x_n)}{f(x_n)}\right) + \bigl( \alpha(x_n) - f(x_n)\bigr) \frac{g(x_n)}{f(x_n)}
\end{equation*}
and since $\alpha$ is bounded, \eqref{eq:simbddA} and the stated negation imply that the negation of the right-hand side of \eqref{eq:simbdd} holds.  The proof of the converse is similar.
\end{proof}

\begin{lemma} \label{lem:bddrec}
Let $a \in \mathbb{R}_+$ and let $f$ and $g$ be positive functions defined on $[a,+\infty)$. Then
\begin{equation}\label{eq:bddrec}
\begin{split}
 \left(\frac{f(x)}{g(x)} - 1 \right)^{-1} = O(1) \ \text{as} \ x  & \to +\infty \\
 &
 \Leftrightarrow
 \quad  \left(\frac{g(x)}{f(x)} - 1 \right)^{-1} = O(1) \ \text{as} \ x \to +\infty.
\end{split}
\end{equation}
\end{lemma}
\begin{proof}
The equivalence of the negations of the propositions in \eqref{eq:bddrec} is clear.
\end{proof}

Application of Theorem~\ref{thm:SL}(\ref{thm:SL-i-1}) and Theorem~\ref{thm:Asym_m2} leads to the following characterization of regularity of critical point $\infty$
under the assumptions of Theorem~\ref{thm:Asym_m2}.

\begin{theorem}\label{thm:Reg_infty_PI}
Let the differential expression $\mathfrak{a}$ satisfy \eqref{eq:rw}. Let $W_\pm$ and $R_\pm$ be the functions defined in \eqref{eq:WR_pm} and assume that $W_\pm \circ R_\pm^{-1}$ is slowly varying function at $0_\pm$. Then the operator $A$ associated with $\mathfrak{a}$ is nonnegative in the Krein space $\bigl(L_{w}^2(I), \kip_w\bigr)$, $\rho(A)\ne\emptyset$, $\infty$ is a critical point of $A$, and
\begin{equation*}
    \infty\in c_r(A)\quad\Leftrightarrow
\quad
\left( 1 + \frac{W_-\bigl(R_-^{-1}(-x)\bigr)}{W_+\bigl(R_+^{-1}(x)\bigr)} \right)^{\!\!-1} = O(1) \quad \text{as}\quad  x\downarrow 0.
\end{equation*}
\end{theorem}
\begin{proof}
Assume that $W_\pm \circ R_\pm^{-1}$ is slowly varying function at $0_\pm$. An immediate consequence of  the definition in \eqref{eq:F_pm} is the equivalence
\begin{equation*}
\left( 1 + \frac{W_-\bigl(R_-^{-1}(-x)\bigr)}{W_+\bigl(R_+^{-1}(x)\bigr)} \right)^{\!\!-1}\!\!= O(1) \  \text{as}\ x \downarrow 0  \ \Leftrightarrow \
\left( 1 - \frac{F_+(x)}{F_-(-x)} \right)^{\!\!-1}\!\!= O(1) \  \text{as} \ x \downarrow 0.
\end{equation*}
Recall that $F_+$ is unbounded decreasing, and $F_-$ is an unbounded increasing function. Since $W_\pm \circ R_\pm^{-1}$ is slowly varying at $0_\pm$, the function $F_\pm$ is regularly varying at $0_\pm$ with index $-1$, see  the definition in \eqref{eq:F_pm}. As the function $f_\pm$ is the inverse of $F_\pm$, Corollary~\ref{cor:RVinv} yields the following equivalence
\begin{equation*}
\left( 1 - \frac{F_+(x)}{F_-(-x)} \right)^{\!\!-1}\!\!= O(1) \  \text{as} \ x \downarrow 0  \ \ \Leftrightarrow \ \
\left( 1 + \frac{f_+(y)}{f_-(y)} \right)^{\!\!-1}\!\!= O(1) \  \text{as} \ y \to +\infty.
\end{equation*}

Let $m_\pm$ be the Neumann $m$-function of $\mathfrak{b}_\pm$ on $I_\pm$. By Theorem~\ref{thm:Asym_m2} we have
\begin{equation*}
  \mp i \, m_\pm(iy)\sim f_\pm(y)\quad \text{as}\quad y\to+\infty.
\end{equation*}
The preceding asymptotic relation and Lemma~\ref{lem:simbdd} imply
\begin{equation*}
\frac{\im m_+(iy)}{m_+(iy)+m_-(-iy)} = O(1) \ \text{as} \ y \to +\infty
\ \Leftrightarrow \ \left(\!1 + \frac{f_-(y)}{f_+(y)}\! \right)^{\!\!-1}\!\!= O(1) \ \text{as} \ y \to +\infty.
\end{equation*}
To see how Lemma~\ref{lem:simbdd} applies here we write
\[
\frac{i \, \im m_+(iy)}{m_+(iy)+m_-(-iy)} = \frac{1}{\frac{-i \, m_+(iy)}{\im m_+(iy)} - \frac{i \,m_-(-iy)}{\im m_+(iy)}},
\]
set
\[
\alpha(y) = \frac{-i \, m_+(iy)}{\im m_+(iy)}, \quad f(y) = \frac{i \,m_-(-iy)}{\im m_+(iy)}, \quad \beta(y) = 1, \quad g(y) = - \frac{f_-(y)}{f_+(y)},
\]
and observe that the above asymptotic relation from Theorem~\ref{thm:Asym_m2} implies
\[
\lim_{y \to + \infty} \alpha(y) = 1 \quad \text{and} \quad \lim_{y \to + \infty} \frac{f(y)}{g(y)} = 1.
\]
Since by Lemma~\ref{lem:bddrec} we have
\[
\left( 1 + \frac{f_-(y)}{f_+(y)} \right)^{\!\!-1}\!\!= O(1) \  \text{as} \ y \to +\infty \ \ \Leftrightarrow \ \ \left( 1 + \frac{f_+(y)}{f_-(y)} \right)^{\!\!-1}\!\!= O(1) \  \text{as} \ y \to +\infty,
\]
we have proved that
\begin{multline*}
\left( 1 + \frac{W_-\bigl(R_-^{-1}(-x)\bigr)}{W_+\bigl(R_+^{-1}(x)\bigr)} \right)^{\!\!-1}\!\!= O(1) \  \text{as}\ x \downarrow 0  \\ \Leftrightarrow \quad
\frac{\im m_+(iy)}{m_+(iy)+m_-(-iy)} = O(1) \ \text{as} \ y \to +\infty.
\end{multline*}
Similarly, we can prove that
\begin{multline*}
\left( 1 + \frac{W_-\bigl(R_-^{-1}(-x)\bigr)}{W_+\bigl(R_+^{-1}(x)\bigr)} \right)^{\!\!-1}\!\!= O(1) \  \text{as}\ x \downarrow 0  \\ \Leftrightarrow \quad
\frac{\im m_-(iy)}{m_+(iy)+m_-(-iy)} = O(1) \ \text{as} \ y \to +\infty.
\end{multline*}
Therefore,
\begin{multline*}
\left( 1 + \frac{W_-\bigl(R_-^{-1}(-x)\bigr)}{W_+\bigl(R_+^{-1}(x)\bigr)} \right)^{\!\!-1}\!\!= O(1) \  \text{as}\ x \downarrow 0  \\ \Leftrightarrow \quad
\text{the pair} \ m_+ \ \text{and} \ m_- \ \text{has the $D_\infty$-property}.
\end{multline*}
Now the theorem follows from Theorem~\ref{thm:SL}.
\end{proof}

\begin{corollary}
Under the assumptions of Theorem~{\rm\ref{thm:Reg_infty_PI}} the following equivalence holds
\begin{equation} \label{eq-chsincp}
    \infty\in c_s(A)\ \ \Leftrightarrow
\ \
\liminf_{x\downarrow 0} \frac{-W_-\bigl(R_-^{-1}(-x)\bigr)}{W_+\bigl(R_+^{-1}(x)\bigr)} \leq 1 \leq \limsup_{x\downarrow 0} \frac{-W_-\bigl(R_-^{-1}(-x)\bigr)}{W_+\bigl(R_+^{-1}(x)\bigr)}.
\end{equation}
\end{corollary}

\begin{proof}
By Theorem~\ref{thm:Reg_infty_PI} $\infty\in c_s(A)$ is equivalent to the negation of
\begin{equation} \label{con:bdd}
\left( 1 + \frac{W_-\bigl(R_-^{-1}(-x)\bigr)}{W_+\bigl(R_+^{-1}(x)\bigr)} \right)^{\!\!-1}  = O(1) \quad \text{as}\quad  x\downarrow 0.
\end{equation}
Subsection~\ref{subsec:AsEq} of Appendix we give two equivalent negations of \eqref{con:bdd}. One is
\begin{equation*}
W_+\bigl(R_+^{-1}(x)\bigr) \ssim -W_-\bigl(R_-^{-1}(-x)\bigr) \quad \text{at} \quad 0_+,
\end{equation*}
and the other, $1$ is a cluster value at $0_+$ of the function
\begin{equation}\label{eq-fracfun}
x \mapsto \frac{-W_-\bigl(R_-^{-1}(-x)\bigr)}{W_+\bigl(R_+^{-1}(x)\bigr)} \quad
\text{with} \quad x \in (0,c),
\end{equation}
where $c = \min\{c_+,-c_-\}$ with $c_-$ and $c_+$ as defined in \eqref{Graz}. Since the function in \eqref{eq-fracfun} is continuous on $(0,c)$ it is an exercise in elementary analysis, see \cite[5.10.11]{TBB}, that $1$ is a cluster value at $0_+$ of the function in \eqref{eq-fracfun} if and only if the inequalities on the right-hand side of the equivalence in \eqref{eq-chsincp} hold.
\end{proof}

\begin{remark}
The criteria in Theorem~\ref{thm:Reg_infty_PI} nicely complements the result of Kostenko in~\cite[Corollary~4.8(i)]{Kost13}. To see this, we notice that \cite[Corollary~4.8(i)]{Kost13} can be restated as follows: If $W_- \circ R_-^{-1}$ is slowly varying function at $0_-$, $W_+ \circ R_+^{-1}$ is slowly varying function at $0_+$ and $\infty \in c_r(A)$, then $w$ is not odd or $r$ is not even.

The ``only if'' part of Theorem~\ref{thm:Reg_infty_PI} gives \eqref{con:bdd} which is more than the fact that $w$ is not odd or $r$ is not even, that is, \eqref{con:bdd} gives that for all small enough positive $x$ we have $W_+\bigl(R_+^{-1}(x)\bigr) \neq -W_-\bigl(R_-^{-1}(-x)\bigr)$.

In this setting the negation of \eqref{con:bdd}, that is the right-hand side of the equivalence in \eqref{eq-chsincp}, appears to be a natural generalization of the condition that the function $w$ is odd and $r$ is even. In the case when $r=1$, this condition also generalizes the condition of $w$ being odd-dominated which was used in Fleige's criterion for $\infty\in c_r(A)$, see \cite[Definition~3.8 and~Theorem~3.11]{CFK13}.
\end{remark}

For slowly varying functions the following corollary extends the result of \cite[Corollary~3.15]{CFK13}.

\begin{corollary} \label{cor:cfk}
Let $0 < b_+ \leq +\infty$, $I_+ = [0,b_+)$ and $r_+, w_+ \in L^1_{\loc}(I_+)$ be positive functions. Let $\alpha, \beta \in \mathbb{R}_+$, set $b_- = - b_+/\beta$ and define
\begin{equation*}
  r(x) = \begin{cases}
           \phantom{\alpha}r_+(x) & \text{if} \quad  x \in [0,b_+) \\
           \alpha r_+(-\beta x) & \text{if} \quad x \in (b_-,0),
         \end{cases}
         \quad
 w(x) = \begin{cases}
           \phantom{-\alpha}w_+(x) & \text{if} \quad  x \in [0,b_+) \\
           -\alpha w_+(-\beta x) & \text{if} \quad x \in (b_-,0).
         \end{cases}
\end{equation*}
Let $W_+$ and $R_+$ be the functions defined in \eqref{eq:WR_pm} and assume that $W_+ \circ R_+^{-1}$ is slowly varying function at $0_+$. Then the operator $A$ associated with $\mathfrak{a}$ is nonnegative in the Krein space $\bigl(L_{w}^2(I), \kip_w\bigr)$, $\rho(A)\ne\emptyset$, $\infty$ is a critical point of $A$, and $\infty \in c_r(A)$ if and only if $\alpha \neq \beta$.
\end{corollary}

\begin{proof}
To apply Theorem~\ref{thm:Reg_infty_PI} we first calculate
 for $x\in (c_-,0)$ (cf.\ \eqref{Graz})
\[
W_-\bigl(R_-^{-1}(x)\bigr) = -(\alpha/\beta) W_+\bigl(R_+^{-1}\bigl(-(\beta/\alpha) x \bigr)\bigr).
\]
Hence $W_- \circ R_-^{-1}$ is a slowly varying function at $0_-$. Further
\[
\frac{W_-\bigl(R_-^{-1}(-x)\bigr)}{W_+\bigl(R_+^{-1}(x)\bigr)} = - \frac{\alpha}{\beta} \frac{W_+\bigl(R_+^{-1}\bigl( (\beta/\alpha) x \bigr)\bigr)}{W_+\bigl(R_+^{-1}(x)\bigr)},
\]
and since $W_+ \circ R_+^{-1}$ is a slowly varying function at $0_+$ we have
\[
\lim_{x\downarrow 0} \frac{W_-\bigl(R_-^{-1}(-x)\bigr)}{W_+\bigl(R_+^{-1}(x)\bigr)} = -\frac{\alpha}{\beta}.
\]
Therefore
\[
\left( 1 + \frac{W_-\bigl(R_-^{-1}(-x)\bigr)}{W_+\bigl(R_+^{-1}(x)\bigr)} \right)^{\!\!-1} = O(1) \quad \text{as}\quad  x\downarrow 0
\]
holds if and only if $\alpha \neq \beta$. Now the claim follows from Theorem~\ref{thm:Reg_infty_PI}.
\end{proof}

We illustrate Corollary~\ref{cor:cfk} with an example which has appeared in \cite[Example~3.17]{CFK13}. The novelty here is that we can give a characterization of the regularity of the critical point $\infty$ for all positive coefficients $\alpha$ and $\beta$.

\begin{example}
Let $w_+, r_+ :(0,1) \to \mathbb{R}_+$ be given by
\begin{equation*} 
    w_+(x)=\frac{1}{x(\ln x)^2}, \quad  r_+(x)=1, \quad x \in (0,1).
\end{equation*}
Then
\[
  W_+(x) = W_+\bigl( R_+^{-1} (x) \bigr) = -\frac{1}{\ln x },\quad x \in [0,1).
\]
Hence, $W_+ \circ R_+^{-1}$ is a slowly varying function at $0_+$. Therefore the operator $A$ from Corollary~\ref{cor:cfk} is nonnegative in the Krein space
$\bigl(L_{w}^2(-1,1), \kip_w\bigr)$,  $\rho(A) \neq \emptyset$, $\infty$ is its critical point and $\infty \in c_r(A)$ if and only if $\alpha \neq \beta$.
\end{example}

\begin{example}\label{ex:Reg1}
Let $\alpha_\pm > 0$ and $I=(-1,1)$. Let $r = 1$ on $I$ and
\begin{equation*}
w_{-}(x)  =\frac{\alpha_-}{x\bigl(-\ln (-x)\bigr)^{1+\alpha_-}},
\ \ x  \in (-1,0), \quad
w_{\scriptscriptstyle{+}}(x)  =\frac{\alpha_+}{x (-\ln x)^{1+\alpha_+}},  \ \ x  \in (0,1).
\end{equation*}
Then
\begin{equation}\label{ex:Reg1_R}
R_-(x) = x, \ \ x \in [-1,0], \quad
R_+(x) = x, \ \ x \in [0,1],
\end{equation}
\begin{equation}\label{ex:Reg1_W}
W_-(x) = \frac{-1}{\bigl(-\ln(-x)\bigr)^{\alpha_-}}, \ \ x \in (-1,0], \quad
W_+(x) = \frac{1}{\bigl(-\ln x\bigr)^{\alpha_+}}, \ \ x \in [0,1).
\end{equation}
Thus $W_- \circ R_-^{-1} = W_-$ is slowly varying at $0_-$, $W_+ \circ R_+^{-1} = W_+$ is slowly varying at $0_+$ and
\begin{equation*}
\left( 1 + \frac{W_-(-x)}{W_+(x)} \right)^{\!\!-1}
= \left(1 - \bigl(- \ln(x) \bigr)^{\alpha_+-\alpha_-} \right)^{-1} = O(1) \ \ \text{as} \ \ x\downarrow 0
\end{equation*}
holds if and only if $\alpha_+ \neq \alpha_-$.

By Theorem~\ref{thm:Reg_infty_PI} the operator $A$ associated with the differential expression $\mathfrak{a}$ with the above defined $w$ and $r$ is nonnegative in the Krein space $\bigl(L_{w}^2(I), \kip_w\bigr)$, $\rho(A) \neq \emptyset$, $\infty$ is its critical point and $\infty \in c_r(A)$ if and only if $\alpha_+ \neq \alpha_-$. That is, $\infty$ is a singular critical point of $A$ if and only if $\alpha_+ = \alpha_-$. Notice that the implication
\begin{equation*}
 \alpha_+ = \alpha_-  \quad \Rightarrow \quad \infty \in c_s(A)
\end{equation*}
follows from a result of Parfenov \cite[Theorem~6]{Par03}, as with $\alpha_+ = \alpha_-$ the weight function $w(x), x\in I$, is odd on $I$.

The converse of the last displayed implication does not follow from neither of the following sufficient conditions for regularity: Volkmer's condition, see \cite[Corollary~2.7]{Vol} or \cite[Theorem~3.14]{CFK13}, Fleige's condition for odd-dominated weights, see~\cite{CFK13}, Parfenov's condition~\cite[Corollary 8]{Par05} for non-odd weights.
\end{example}

\subsection{Discreteness}\label{subsec:disc} For a closed operator $T$ its discrete spectrum consists of its isolated eigenvalues of finite algebraic multiplicity. The complement of the discrete spectrum is called the essential spectrum of $T$; it is denoted by $\sigma_{\!\operatorname{ess}}(T)$.
The differential expression $\mathfrak{b}_{+}$ is said to be quasi-regular at the end-point $b_+$ if $w_+,r_+\in L^1(I_+)$.
As is known, see~\cite{KacKr58}, in the quasi-regular case the spectrum of the operator $B_{+,0}$ is discrete.
The following statement for a non-quasi-regular case is also based on a result from~\cite{KacKr58}.

\begin{theorem}\label{lem:KacKr}
Let $0< b_+ \leq +\infty$, let  $B_{+,0}$ be defined by~\eqref{eq:B_pm0} and let either
$w_+ \not\in L^1(0,b_+)$ or $r_+ \notin L^1(0,b_+)$.
 Then  $0\not\in\sigma_{\!\operatorname{ess}}(B_{+,0})$ if and only if:

 \noindent Either
\begin{enumerate}
\renewcommand*\theenumi{\Roman{enumi}}
\renewcommand*\labelenumi{\rm{(\theenumi)}}
  \item  \label{lem:KacKr-i1}
  $w_+ \in L^1(0,b_+)$, $r_+ \notin L^1(0,b_+)$ and
\begin{equation}\label{eq:Ess_Sp} 
 \sup_{x\in(0,b_+)} R_+(x) \bigl(W_+(b_+)-W_+(x)\bigr)  <+\infty,
\end{equation}
\end{enumerate}
or
\begin{enumerate}
\renewcommand*\theenumi{\Roman{enumi}}
\renewcommand*\labelenumi{\rm{(\theenumi)}}
\setcounter{enumi}{1}
  \item \label{lem:KacKr-i2}
 $w_+ \not\in L^1(0,b_+)$, $r_+ \in L^1(0,b_+)$ and
\begin{equation*} 
 \sup_{x\in(0,b_+)} W_+(x) \bigl(R_+(b_+)-R_+(x) \bigr) < +\infty.
\end{equation*}
\end{enumerate}
Moreover, the spectrum of  $B_{+,0}$ is discrete if and only if:

\noindent Either

\begin{enumerate}
\renewcommand*\theenumi{\Roman{enumi}}
\renewcommand*\labelenumi{\rm{(\theenumi)}}
\setcounter{enumi}{2}
  \item
$w_+ \in L^1(0,b_+)$, $r_+ \notin L^1(0,b_+)$ and
\begin{equation}
\label{eq:Disc_Sp0}
 \lim_{x\to b_+}  R_+(x) \bigl(W_+(b_+)-W_+(x) \bigr)   = 0,
\end{equation}
\end{enumerate}
or
\begin{enumerate}
\renewcommand*\theenumi{\Roman{enumi}}
\renewcommand*\labelenumi{\rm{(\theenumi)}}
\setcounter{enumi}{3}
  \item
$w_+ \not\in L^1(0,b_+)$, $r_+ \in L^1(0,b_+)$ holds and
\begin{equation*} 
 \lim_{x\to b_+}  W_+(x) \bigl(R_+(b_+)-R_+(x) \bigr)   = 0.
\end{equation*}
\end{enumerate}

\end{theorem}
\begin{proof}
By using the change of variable $\xi=R_+(x)$, $x\in (0,b_+)$ the statements of Theorem~\ref{lem:KacKr} are easily reduced to~\cite{KacKr58}, see~\cite{CDT21} for the details.
\end{proof}

The statements of Theorem~\ref{lem:KacKr} remain in force for $B_{-,0}$ with $b_+$, $w_+$, $r_+$ replaced by $b_-$, $w_-$, $r_-$, respectively. In particular, if $w_- \in L^1(b_-,0)$, then $0\not\in\sigma_{\!\operatorname{ess}}(B_{{-},0})$ if and only if
\begin{equation}\label{eq:Ess_Sp-}
 \sup_{x\in(b_{-},0)} R_{-}(x) \bigl(W_-(b_-)- W_-(x)\bigr) < +\infty.
\end{equation}

\begin{remark}
It seems  Kac and Krein~\cite{KacKr58} were the first to introduce condition~\eqref{eq:Disc_Sp0}, with $r_+\equiv 1$,
as a discreteness criterion for the  Krein string.
A condition similar to~\eqref{eq:Ess_Sp} was used by Chisholm and Everitt~\cite{ChisEve71}
as a criterion for the boundedness of the integral operator
$(Tf)(x) =  v(x)\int_0^x u(t)f(t)dt$ with $f  \in L^2(\dR_+)$. Here $u,v\in L^2(\dR_+)$. Stuart~\cite{Stu72} proved that the compactness of the operator $T$ is characterized by a condition of type~\eqref{eq:Disc_Sp0}
and this allowed  him to characterize the discreteness of a general Sturm-Liouville operator. See also \cite{OinOtel99} and~\cite{CRead} where
discreteness criteria were formulated in terms of the coefficients of the Sturm-Liouville operator. Conditions similar to~\eqref{eq:Ess_Sp} appeared also in papers by Muckenhaupt~\cite{Muck72},  \cite[Section~1.3.1]{Maz85} as criteria for some Hardy-type inequalities in weighted spaces. Criteria for the discreteness of the spectra of canonical systems, that contain the Krein string as a special case, were found recently in~\cite{RomWor}, see  also~\cite{RemSc} for a class of semibounded canonical systems.
\end{remark}

The next proposition shows that the spectrum of the operator $B_{\pm,0}$ defined in~\eqref{eq:B_pm0} can be  discrete even in the limit point case.

\begin{proposition}\label{lem:2intd}
Assume $w_\pm\in L^1(I_\pm)$. Then
\begin{equation}\label{eq-2int}
 \int_{I_\pm} \bigl| R_\pm(\xi) \bigr| w_{\pm}(\xi) d\xi  =  \int_{I_\pm}\bigl| W_\pm(b_{\pm})- W_\pm(\xi)\bigr| r_\pm(\xi) d\xi,
\end{equation}
meaning that either the two integrals diverge simultaneously, or, if one converges, then the other one converges as well and the integrals are equal. Further, if
 $R_\pm\in L^1_{w_\pm}(I_\pm)$, then the spectrum of $B_{\pm,0}$ is discrete and $0\in\wh\rho(B_\pm)$.
\end{proposition}

\begin{proof}
Equality~\eqref{eq-2int} is verified using integration by parts in  $\int_{I_\pm} \bigl| R_\pm(\xi) \bigr| w_{\pm}(\xi) d\xi$.

Assume now that $R_\pm\in L^1_{w_\pm}(I_\pm)$.
Applying again integration by parts to the integral
$\int_0^x \bigl(W_\pm(b_\pm)- W_\pm(\xi)\bigr) dR_\pm(\xi)$ we obtain for all $x \in  I_\pm$
\begin{equation}\label{eq:Est10}
\begin{split}
\int_0^x \bigl(W_\pm(b_\pm)- W_\pm(\xi)\bigr) & dR_\pm(\xi) \\
&=
  R_\pm(x)\bigl(W_\pm(b_\pm)-W_\pm(x)\bigr)+\int_0^x R_\pm(\xi)dW_\pm(\xi).
  \end{split}
\end{equation}
Taking the limit as $x\to b_\pm$ in \eqref{eq:Est10} and using ~\eqref{eq-2int} yields
\begin{equation}\label{eq:LimRW}
  \lim_{x\to b_\pm} R_\pm(x)\bigl(W_\pm(b_\pm)-W_\pm(x)\bigr) = 0.
\end{equation}
Hence, by Theorem~\ref{lem:KacKr}, the spectrum of $B_{\pm,0}$ is discrete.
\end{proof}

The next theorem combines the results of Theorem~\ref{lem:KacKr} 
and Theorem~\ref{thm:Reg_infty_PI} to provide a necessary and sufficient condition for the existence of a Riesz basis consisting of eigenfunctions of the differential operator $A$.

\begin{theorem} \label{thm:Riesz}
Let the differential expression $\mathfrak{a}$ satisfy \eqref{eq:rw} and let $W_\pm$ and $R_\pm$ be the functions defined in \eqref{eq:WR_pm}. Assume

\begin{enumerate}
\renewcommand*\theenumi{\alph{enumi}}
\renewcommand*\labelenumi{\rm{(\theenumi)}}
\item \label{thm:Riesz-i+d}
The functions $w_+$ and $r_+$ satisfy one of the following three conditions:
\begin{enumerate}
\renewcommand*\theenumii{\roman{enumii}}
\renewcommand*\labelenumii{\rm{(\theenumii)}}
\item \label{thm:Riesz-i+d-1}
$w_+\in L^1(I_+)$ and $r_+\in L^1(I_+).$
\item \label{thm:Riesz-i+d-2}
$w_+\in L^1(I_+)$, $r_+ \not\in L^1(I_+)$ and
$\lim_{x\uparrow b_+} R_+(x) \bigl(W_+(b_+) - W_+(x)\bigr) = 0.$
\item \label{thm:Riesz-i+d-3}
$w_+ \not\in L^1(I_+)$, $r_+\in L^1(I_+)$ and
$\lim_{x\uparrow b_+} W_+(x) \bigl(R_+(b_+) - R_+(x) \bigr)= 0.$
\end{enumerate}
\item \label{thm:Riesz-i-d}
The functions $w_-$ and $r_-$ satisfy one of the following three conditions
\begin{enumerate}
\renewcommand*\theenumii{\roman{enumii}}
\renewcommand*\labelenumii{\rm{(\theenumii)}}
\item \label{thm:Riesz-i-d-1}
$w_-\in L^1(I_-)$ and $r_-\in L^1(I_-)$.
\item \label{thm:Riesz-i-d-2}
$w_-\in L^1(I_-)$, $r_-\not\in L^1(I_-)$ and
$\lim_{x\downarrow b_-} R_-(x) \bigl(W_-(b_-) - W_-(x)\bigr) = 0.$
\item \label{thm:Riesz-i-d-3}
$w_- \not\in L^1(I_-)$, $r_-\in L^1(I_-)$ and
$\lim_{x\downarrow b_-} W_-(x) \bigl(R_-(b_-) - R_-(x) \bigr) = 0.$
\end{enumerate}
\end{enumerate}

Then the spectrum of the operator $A$ associated with the differential expression $\mathfrak{a}$ in the Hilbert space $L_{|w|}^2(I)$ is real and discrete, its eigenvalues accumulate on both sides of $\infty$, all nonzero eigenvalues are simple and Jordan chain at $0$ is of length at most $2$. The following statements hold.
\begin{enumerate}
\renewcommand*\theenumi{\Alph{enumi}}
\renewcommand*\labelenumi{\rm{(\theenumi)}}
\item \label{thm:Riesz-B}
If either $W_+\circ R_+^{-1}$ is positively increasing at $0_+$ or $W_-\circ R_-^{-1}$ is positively increasing at $0_-$, then
$A$ has the Riesz basis property {\rm (Ri)}.

\item \label{thm:Riesz-A}
If $W_+ \circ R_+^{-1}$ is slowly varying at $0_+$ and $W_- \circ R_-^{-1}$ is slowly varying at $0_-$,
then $A$ has the Riesz basis property {\rm (Ri)} if and only if
\begin{equation}\label{eq:thm4.26}
  \left( 1 + \frac{W_-\bigl(R_-^{-1}(-x)\bigr)}{W_+\bigl(R_+^{-1}(x)\bigr)} \right)^{\!\!-1} = O(1) \quad \text{as}\quad  x\downarrow 0.
\end{equation}
\end{enumerate}
\end{theorem}

\begin{proof}
In either of the three cases in (\ref{thm:Riesz-i+d}), the spectrum of the operator $B_{+,0}$ is discrete and its eigenvalues accumulate at $+\infty$. This follows from the fact that in case (\ref{thm:Riesz-i+d-1}) in (\ref{thm:Riesz-i+d}) the operator $B_{+,0}$ is either regular or in the limit-circle case at $b_+$. In the remaining two cases in (\ref{thm:Riesz-i+d}) this follows from Theorem~\ref{lem:KacKr}. Similarly, in either of the three cases in (\ref{thm:Riesz-i-d}), the spectrum of the operator $B_{-,0}$ is discrete and its eigenvalues accumulate at $+\infty$.
Since $A$ is a rank-one perturbation of the operator $B_{+,0}\oplus (-B_{-,0})$,
by Weyl's theorem the spectrum of the operator $A$ is also discrete (see \cite[Theorem~XIII.14]{RS}). By Lemma~\ref{thm:Reg_infty} the eigenvalues of $A$ accumulate on both sides of $\infty$. Since the operator $A$ is nonnegative in the Krein space $\cK$ all nonzero eigenvalues of $A$ are semi-simple and the length of the Jordan chain at $0$ is at most $2$. Moreover, by Lemma~\ref{lem:Ker} all nonzero eigenvalues of $A$ are simple.

Let $\Delta$ be an arbitrary finite open interval such that $0\in\Delta$ and let $E$ be the spectral function of $A$ in the sense of~\cite{La}. By the properties
of this spectral function~\cite{La}, $\infty\in c_r(A)$ if and only if  there exists a Riesz basis of ${(I-E(\Delta))\cK}$ which consists of eigenfunctions and the generalized eigenfunctions of the restriction of $A$ on $(I-E(\Delta))\cK$. Since $E(\Delta)\cK$ is a finite-dimensional space, the eigenfunctions and the generalized eigenfunctions of the restriction of $A$ on ${E(\Delta)\cK}$ form a Riesz basis of $E(\Delta)\cK$. Therefore,  the Riesz basis property (Ri) 
is equivalent to $\infty\in c_r(A)$. By Theorem~\ref{thm:Reg_infty_PIS}, if either $W_+\circ R_+^{-1}$ is positively increasing at $0_+$ or $W_-\circ R_-^{-1}$ is positively increasing at $0_-$, then  $\infty\in c_r(A)$ and hence the claim in (\ref{thm:Riesz-B}) holds.

If $W_+ \circ R_+^{-1}$ is slowly varying at $0_+$ and $W_- \circ R_-^{-1}$ is slowly varying at $0_-$,
then by Theorem~\ref{thm:Reg_infty_PI} condition~\eqref{eq:thm4.26} 
is equivalent to $\infty\in c_r(A)$. Since we already proved that $\infty\in c_r(A)$ is equivalent to (Ri), 
the equivalence in (\ref{thm:Riesz-A}) is proved.
\end{proof}

\begin{remark}
For the differential expression $\mathfrak{a}$ introduced in Example~\ref{ex:Reg1} we have
\[
W_+(x)\bigl( R_+(1)-R_+(x) \bigr) = \frac{1-x}{(-\ln x)^{\alpha_+}}\sim (1-x)^{1-\alpha_+} \quad \text{as} \quad  x\uparrow 1
\]
and
\[
W_-(x)(R_-(-1)-R_-(x)) = \frac{1+x}{\bigl(-\ln(-x)\bigr)^{\alpha_-}}\sim (1+x)^{1-\alpha_-} \quad \text{as} \quad  x\downarrow -1.
\]
Therefore $\mathfrak{a}$ satisfies conditions (\ref{thm:Riesz-i+d})(\ref{thm:Riesz-i+d-3}) and (\ref{thm:Riesz-i-d})(\ref{thm:Riesz-i-d-3}) in Theorem~\ref{thm:Riesz} if and only if $\alpha_-\in(0,1)$ and $\alpha_+\in(0,1)$.
By (\ref{thm:Riesz-A}) in Theorem~\ref{thm:Riesz} the operator $A$ in Example~\ref{ex:Reg1} with $\alpha_-,\alpha_+\in(0,1)$ has the Riesz basis property if and only if $\alpha_-\ne\alpha_+$.
\end{remark}

\subsection{Regularity at $0$}\label{sec:5.4}

Since the operator $A$ associated with the differential expression $\mathfrak{a}$ is nonnegative it may have another critical point at $0$. In this subsection we consider the problem of regularity of the critical point $0$ of the operator $A$. Let  $W_\pm$ and $R_\pm$ be defined by~\eqref{eq:WR_pm}.

\begin{theorem}\label{thm:Reg0}
 Let  $W_\pm$ and $R_\pm$ be defined by~\eqref{eq:WR_pm},
    and let $A$ be the differential operator associated with the expression $\mathfrak{a}$ with the domain defined by~\eqref{eq:domA}.
Assume that one of the following cases is  in force:
\begin{enumerate}
\renewcommand*\theenumi{\roman{enumi}}
\renewcommand*\labelenumi{\rm{(\theenumi)}}
 \item \label{thm:Reg0i1}
$w_-$, $r_-\not\in L^1(I_-)$, $w_+$, $r_+\not\in L^1(I_+)$ and
     either $W_-\circ R_-^{-1}$ is positively increasing at $-\infty$ or $W_+\circ R_+^{-1}$ is positively increasing at $+\infty$;
\item  \label{thm:Reg0i2}
$w_-\not\in L^1(I_-)$ and $w_+\not\in L^1(I_+)$ and either $r_+\in L^1(I_+)$,   or $r_-\in L^1(I_-)$;
 \item  \label{thm:Reg0i3}
either $w_+\in L^1(I_+)$, $w_-\not\in L^1(I_-)$
     or $w_-\in L^1(I_-)$, $w_+\not\in L^1(I_+)$.
    \item  \label{thm:Reg0i4}
$w_-\in L^1(I_-)$ and $w_+\in L^1(I_+)$, and $  W_+(b_+)+W_-(b_-)\ne 0.$
\end{enumerate}
Then
\begin{equation}\label{eq:Reg0_SL}
  0\not\in c_s(A)\quad\text{and}\quad \ker A=\ker A^2.
\end{equation}
Moreover, the following statements hold.
\begin{enumerate}
\renewcommand*\theenumi{\alph{enumi}}
\renewcommand*\labelenumi{\rm{(\theenumi)}}
 \item \label{thm:Reg0i5}
If $w_-\in L^1(I_-)$ and $w_+\in L^1(I_+)$, then \eqref{eq:Reg0_SL} holds if and only if $W_+(b_+)+W_-(b_-)\ne 0.$

\item \label{thm:Reg0i6}
If $w_-\in L^1(I_-)$, $w_+\in L^1(I_+)$ and~\eqref{eq:Ess_Sp}, \eqref{eq:Ess_Sp-} hold,
then $0\not\in\sigma_{\!\operatorname{ess}}(A)$ and the following three statements are equivalent
\begin{equation}\label{eq:KerAA2}
W_+(b_+)+W_-(b_-)\ne 0 \quad \Leftrightarrow \quad  \ker A = \ker A^2 \quad \Leftrightarrow \quad 0 \notin c(A).
\end{equation}
\end{enumerate}
\end{theorem}

\begin{proof}
{\bf 1.} {\it Proof of} \eqref{eq:Reg0_SL} {\it under assumption} (\ref{thm:Reg0i1}).
 Due to Lemma~\ref{cor:R_0} the assumption that $W_+\circ R_+^{-1}$  is positively increasing at $+\infty$ is equivalent to the condition
 \begin{equation*} 
    \im m_+(iy)=O\bigl(\re m_+(iy)\bigr) \quad \text{as} \quad y\downarrow 0.
\end{equation*}
By Theorem~\ref{thm:SL} (\ref{thm:SL-i-5})
        this implies $0\not\in c_{\mathrm{s}}(A)$ and $\ker A=\ker A^2$.

\medskip

\noindent
{\bf 2.} {\it Proof of} \eqref{eq:Reg0_SL} {\it under assumption} (\ref{thm:Reg0i2}).
If $r_+\in L^1(\mathbb{R}_+)$ and  $w_+\not\in L^1(\mathbb{R}_-)$, then
by Lemma \ref{lem:PolF} \eqref{eq:Pol_m2} holds and, hence, by Theorem~\ref{thm:SL} (\ref{thm:SL-i-5}) we have $0\not\in c_{\mathrm{s}}(A)$ and $\ker A=\ker A^2$.

\medskip

\noindent
{\bf 3.}  {\it Proof of} \eqref{eq:Reg0_SL} {\it under assumption} (\ref{thm:Reg0i3}).
If $w_+\in L^1(I_+)$ and $w_-\not\in L^1(I_-)$, then  by  Lemma \ref{lem:PolF}
\begin{equation*}
  m_+(iy)=i\frac{a_+}{y}+\wt m_+(iy),\quad   m_-(iy)=o(1/y) \quad \text{as} \quad y\downarrow 0
    \end{equation*}
    for  $a_+=\dfrac{1}{W_+(b_+)} > 0$, $\wt m_+(iy)=o(1/y)$.
 Then
    \[
     m_+(iy)+m_-(-iy)\sim -\frac{a_+}{iy} \quad \text{as} \quad y\downarrow 0.
    \]
  and
\[
 \im m_+(iy)\sim \frac{a_+}{y},\quad \im m_-(iy)\to 0 \quad \text{as} \quad y\downarrow 0.
\]
Hence
\begin{equation*}
    \im  m_\pm(iy)=O(m_+(iy)+m_-(-iy)) \quad \text{as} \quad y\downarrow 0
\end{equation*}
and by Theorem~\ref{thm:SL}(\ref{thm:SL-i-2}) $0\not\in c_s(A)$ and $\ker A=\ker A^2$.

\medskip

\noindent
{\bf 4.}  {\it Proof of} \eqref{eq:Reg0_SL} {\it under assumption} (\ref{thm:Reg0i4}).
If $w_+\in L^1(I_+)$ and $w_-\in L^1(I_-)$ then  by  Lemma \ref{lem:PolF}
\[
  m_+(iy)\sim i\frac{a_+}{y},\quad   m_-(iy)\sim i\frac{a_-}{y} \quad \text{as} \quad y\downarrow 0
\]
for  $a_\pm=\pm 1/W_\pm(b_\pm)$.
Since
$W_+(b_+)+W_-(b_-)\ne 0$ then $a_+\ne a_-$,
\[
     m_+(iy)+m_-(-iy)\sim i\frac{a_+-a_-}{y} \quad \text{as} \quad y\downarrow 0
\]
and
\begin{equation}\label{eq:Im_m5}
     \im m_\pm(iy)\sim \frac{a_\pm}{y} \quad \text{as} \quad y\downarrow 0.
\end{equation}
Hence
\begin{equation}\label{eq:Pol_m5}
    \im m_\pm(iy)=O\bigl(  m_+(iy)+m_-(-iy) \bigr) \quad \text{as} \quad y\downarrow 0.
\end{equation}
By Theorem~\ref{thm:SL}(\ref{thm:SL-i-2}) \eqref{eq:Pol_m5} is equivalent to $0\not\in c_s(A)$ and $\ker A=\ker A^2$.

\medskip

\noindent
{\bf 5.} {\it Proof of} (\ref{thm:Reg0i5}).
Assume now that $W_+(b_+)+W_-(b_-)=0$. Then by  Lemma \ref{lem:PolF} $a_+= a_-$  and hence
\[
m_+(iy)+m_-(-iy)=o(1/y)  \quad \text{as} \quad y\downarrow 0.
\]
In view of~\eqref{eq:Im_m5} the relation \eqref{eq:Pol_m5} is not fulfilled and     by Theorem~\ref{thm:SL} the relations~\eqref{eq:Reg0_SL} fail to hold, i.e either $0\in c_s(A)$ or $\ker A\subsetneq\ker A^2$.

\medskip

\noindent
{\bf 6.}  
{\it Proof of} (\ref{thm:Reg0i6}).
If $w_\pm\in L^1(I_\pm)$ and~\eqref{eq:Ess_Sp}, \eqref{eq:Ess_Sp-} hold, then by Theorem~\ref{lem:KacKr} $0\not\in\sigma_{\!\operatorname{ess}}(B_{\pm,0})$. Since $A$ is a rank-one perturbation of the operator $B_{+,0}\oplus (-B_{-,0})$ we have
$0\not\in \sigma_{\!\operatorname{ess}}(A)$.

Since $w \in L^1(I)$, all constant functions on $I$ belong to $\dom A$ defined in \eqref{eq:domA}. Consequently, all constant functions on $I$ belong to $\ker A$. As by Lemma~\ref{lem:Ker} $\ker A$ is at most one-dimenional, we deduce that $\ker A$ consists of all the constant functions on $I$. Denote by $\mathbf{1}$ the constant function on $I$ equal to $1$. Notice that
\begin{equation}\label{eq:h^0}
  [\mathbf{1},\mathbf{1}]_w=W_+(b_+)+W_-(b_-).
\end{equation}

If $W_+(b_+)+W_-(b_-)\ne 0$ then the subspace $\ker A$ is nondegenerate. Moreover, we have $\ker A = \ker A^2$, since the existence of an associated vector $f \in \dom A$ such that $Af = \mathbf{1}$ implies $[\mathbf{1},\mathbf{1}]_w = [Af,\mathbf{1}]_w = [f,A\mathbf{1}]_w = 0$. This proves the implication
\[
W_+(b_+)+W_-(b_-)\ne 0 \quad \Rightarrow \quad  \ker A = \ker A^2.
\]

Now assume $\ker A = \ker A^2$. Then $0$ is a simple eigenvalue of $A$ and since $0\not\in \sigma_{\!\operatorname{ess}}(A)$ it is an isolated eigenvalue. By \cite{La} $\ker A$ is nondegenerate and thus $0\not\in c(A)$.

And finally, if $0\not\in c(A)$, then by \cite{La} $\ker A^2 = \ker A$ and $\ker A$ is definite. Hence $[\mathbf{1},\mathbf{1}]_w \ne 0$ and by~\eqref{eq:h^0} we have $W_+(b_+)+W_-(b_-)\ne 0$. This proves the implication $0\not\in c(A)\Rightarrow W_+(b_+)+W_-(b_-)\ne 0$ and hence the
equivalence~(\ref{eq:KerAA2}).
\end{proof}

\begin{remark}\label{rem:Karabash09}
The first equivalence in~\eqref{eq:KerAA2} can also be  derived from~\cite[Theorem 3.1]{Kar10}. Indeed, let $\sigma_\pm$ be measures from the integral representations $m_\pm(z)=\int_\dR (t-z)^{-1}d\sigma_\pm(t)$ of the $m$-functions $m_\pm$. If $w_-\in L^1(I_-)$, $w_+\in L^1(I_+)$ and~\eqref{eq:Ess_Sp}, \eqref{eq:Ess_Sp-} hold,
then $0\not \in \sigma_{\operatorname{ess}}(B_-) \cap \sigma_{\operatorname{ess}}(B_+)$. Hence the conditions
\[
\int_{\dR\setminus\{0\}} t^{-2} d\sigma_-(t)<+\infty,\quad \int_{\dR\setminus\{0\}}t^{-2}d\sigma_+(t)<+\infty
\]
are automatically fulfilled.
By~\cite[Theorem 3.1, 2(ii)]{Kar10},   the condition  $\ker A^2 = \ker A$ is equivalent to $d\sigma_-(\{0\})\ne d\sigma_+(\{0\})$, which, by Lemma~\ref{lem:PolF},  is equivalent to $  W_+(b_+)+W_-(b_-)\ne 0$.
\end{remark}

In Theorem~\ref{thm:Reg0} it is not claimed that $0\in c_r(A)$, since it may happen that $0$ is not a critical point  of $A$ at all. In the next corollary we specify some cases when $0$ is indeed a regular critical point of $A$.

\begin{corollary}\label{cor:Reg0}
Assume that $w_+$ and $r_+$ satisfy one of the following assumptions:
  \begin{enumerate}
\renewcommand*\theenumi{\alph{enumi}}
\renewcommand*\labelenumi{\rm{(\theenumi)}}
    \item  \label{cor:Crit-a}
$w_+ \in L^1(I_+)$, $r_+ \notin L^1(I_+)$ and $\sup_{x\in I_+} R_+(x) \bigl(W_+(b_+) - W_+(x)\bigr)  =+\infty$,
    \item  \label{cor:Crit-b}
$w_+ \not\in L^1(I_+)$, $r_+ \in L^1(I_+)$ and $\sup_{x\in I_+}   W_+(x) \bigl(R_+(b_+) - R_+(x)\bigr)  =+\infty$.
  \end{enumerate}
Assume that $w_-$ and $r_-$ satisfy one of the following assumptions:
  \begin{enumerate}
\renewcommand*\theenumi{\alph{enumi}}
\renewcommand*\labelenumi{\rm{(\theenumi)}}
\setcounter{enumi}{2}
    \item  \label{cor:Crit-c}
  $w_-\in L^1(I_-)$, $r_-\not\in L^1(I_-)$ and
$\sup_{x\in I_-}  R_-(x) \bigl( W_-(b_-) - W_-(x) \bigr) = +\infty$,
    \item  \label{cor:Crit-d}
$w_- \not\in L^1(I_-)$, $r_-\in L^1(I_-)$ and
$\sup_{x\in I_-}  W_-(x) \bigl( R_-(b_-) - R_-(x) \bigr) = + \infty$.
\end{enumerate}
In cases {\rm (\ref{cor:Crit-a})} and {\rm (\ref{cor:Crit-c})} assume $W_+(b_+)+W_-(b_-)\ne 0$. Then $0\in c_r(A)$ and the spectrum of the operator $A$  accumulates on both sides of $0$.
\end{corollary}
\begin{proof}
In either of the cases (\ref{cor:Crit-a}) and (\ref{cor:Crit-b}) ((\ref{cor:Crit-c}) and (\ref{cor:Crit-d}), respecitively), $0$ is an accumulation point for the spectrum of the operator $B_{+,0}$ ($B_{-,0}$, respectively) from the right. Therefore, $0$ is an accumulation point for the spectrum of the decoupled operator $A_0=A_{+,0}\oplus (A_{-,0})$ from both sides.
Since the resolvent $(A-z)^{-1}$ of $A$ is a one-dimensional perturbation of  the resolvent $(A_0-z)^{-1}$, see~\eqref{ereswA1}, it follows from~\cite[Theorem 1]{JonLan79} that $0\in c(A)$.

 The statement $0\in c_r(A)$ follows from Theorem~\ref{thm:Reg0}.
\end{proof}

\begin{remark}
  The list of assumptions of Theorem~\ref{thm:Reg0} covers all possible cases except the following:
   \begin{enumerate}
     \item [(v)] $w_-$, $r_-\not\in L^1(I_-)$, $w_+$, $r_+\not\in L^1(I_+)$ and
     both $W_-\circ R_-^{-1}$ is not positively increasing at $-\infty$ and $W_+\circ R_+^{-1}$ is not positively increasing at $+\infty$.
   \end{enumerate}
In this case we cannot apply our abstract results from Theorem~\ref{prop:Semibound2} because the asymptotic behaviour of the Weyl functions at finite points is insufficiently studied.

Notice that if $r_\pm = 1$, then \cite[Theorem~5.2]{Vol} and \cite[Theorem~7.4]{Kost14} yield that the set all not positively increasing functions is dense in some metric subspace of $L^1(I_\pm)$.

\end{remark}

\begin{remark}
 If $w_+\in L^1(I_+)$, $w_-\in L^1(I_-)$, $R_+\in L^1_{w_+}(I_+)$ and $R_-\in L^1_{w_-}(I_-)$
then by Proposition~\ref{lem:2intd} the spectrum of $A$ is discrete and in the case
$W_+(b_+)+W_-(b_-)=0$ the root subspace $\ker A^2$  can be found explicitly.
As was mentioned above $\ker A = \lspan\{\mathbf{1}\}$.
Let  us find a generalized eigenvector $f\in\dom(A)$ such that  $Af=\mathbf{1}$, i.e. $f=f_+\oplus f_-$, where $f_\pm\in \dom(B_{\pm}^{\langle *\rangle})$ are solutions of the equations
\begin{equation}\label{eq:Jord}
  \mathfrak{b}_+f_+=1,\quad
     -\mathfrak{b}_-f_-=1,
\end{equation}
such that
\begin{equation}\label{eq:A_Coupl2}
f_+(0)=f_-(0), \quad  f^{[1]}_+(0)=f^{[1]}_-(0).
\end{equation}
holds. Straightforward calculations show that the functions
\begin{equation}\label{eq:h1}
f_\pm(x)=\pm \int_0^{x}R_\pm(\xi)w_\pm(\xi)d\xi \pm \int_x^{b_\pm}R_\pm(x)w_\pm(\xi)d\xi
\end{equation}
satisfy \eqref{eq:Jord} and the first  boundary condition in \eqref{eq:A_Coupl2}.
The second boundary condition in \eqref{eq:A_Coupl2} holds since  $W_+(b_+)+W_-(b_-)=0$.

It follows from~\eqref{eq:LimRW} that the second term in~\eqref{eq:h1}
\[
\int_x^{b_{\pm}}R_\pm(x)w_\pm(\xi)d\xi = R_\pm(x)(W_\pm(b_{\pm})-W_\pm(x))
\]
is bounded. The first term in the right hand part of~\eqref{eq:h1} is also bounded since $R_\pm\in L^1_{w_\pm}(I_\pm)$ and hence $f_\pm\in L^2_{w_\pm}(I_\pm)$.
Therefore, $f_\pm\in \dom(B_{\pm,\max})$ and hence $f_\pm\in\dom(B_\pm^{\langle*\rangle})$ in the limit point case.

In the limit circle case we also get $f_\pm\in\dom(B_\pm^{\langle*\rangle})$, since $f_\pm^{[1]}(b_\pm)=0$. Therefore, $f= f_+\oplus f_-\in\dom A$ and the equation $Af=\mathbf{1}$ has a solution $f\in \dom(A)$. Thus $\ker A\ne\ker A^2$.
\end{remark}

\begin{remark}
Let $b_+$, $\alpha, \beta$, $b_-$ and the function $w$ be defined as in Corollary~\ref{cor:cfk} and arbitrary $r \in L^1_{\operatorname{loc}}(I)$. Assume that $w_+\in L^1(I_+)$. Then $w_-\in L^1(I_-)$. Let $W_+$ be the function defined in \eqref{eq:WR_pm}. Then, for all $x \in [b_-,0]$ we have $W_-(x)=-(\alpha/\beta) W_+(-\beta x)$. Consequently,
\begin{equation*}
W_+(b_+) + W_-(b_-) = \left(1 - \alpha/\beta\right) W_+(b_+).
\end{equation*}
By Theorem~\ref{thm:Reg0}~\!(\ref{thm:Reg0i5}) we have that \eqref{eq:Reg0_SL} holds if and only if $\alpha \neq \beta$. In particular, if $\alpha = \beta = 1$ the weight function $ w(t)$ is odd and condition \eqref{eq:Reg0_SL} does not hold. This has been proved in \cite[Theorem~4.7]{Kost13} under additional conditions that $r$ is even and $r \notin L^1(I)$.
\end{remark}

\begin{example}\label{ex:Reg1_0}
We consider the differential expression studied in Example~\ref{ex:Reg1} on an interval $I=(b_-,b_+)$ with $-1 \leq b_- < 0 < b_+ \leq 1$.

First assume that $b_- = -1$ and $b_+ =1$, as in Example~\ref{ex:Reg1}. Then $R_\pm$ and $W_\pm$ are given by the formulas~\eqref{ex:Reg1_R} and ~\eqref{ex:Reg1_W}.
 Due to Theorem~\ref{thm:Reg0}(\ref{thm:Reg0i2})  $0\not\in c_s(A)$ and $\ker A=\ker A^2$.

Moreover, $\ker(A)=\{0\}$ since a function $f\in\ker(A)$
should have a form   $f= f_+\oplus f_-\in\dom A$, where $f_+\in\dom(B_+^{\langle*\rangle})$, $f_-\in\dom(B_-^{\langle*\rangle})$ and  the coupling conditions \eqref{eq:A_Coupl2} hold.
The conditions $f_+\in\dom(B_+^{\langle*\rangle})$, $f_-\in\dom(B_-^{\langle*\rangle})$ yield that $f_+$ and $f_-$ are proportional to $1-x$ and $1+x$, respectively.
But then the coupling conditions~\eqref{eq:A_Coupl2}  yield $f_+=f_-=0$.

Further,
\begin{equation*}
\lim_{x\uparrow 1} \bigl(R_+(1)-R_+(x)\bigr)W_+(x)=\lim_{x\uparrow 1}\frac{1-x}{(-\ln x)^{\alpha_+}}
=\left\{ \begin{array}{cl}
           0 & \text{if} \quad 0 < \alpha_+ < 1, \\
           1 & \text{if} \quad \alpha_+ = 1, \\
           +\infty& \text{if} \quad \alpha_+ > 1,
         \end{array} \right.
\end{equation*}
and
\begin{equation*}
\lim_{x\downarrow -1} \bigl(R_-(-1)-R_-(x)\bigr)W_-(x)=\lim_{x\downarrow -1}\frac{1+x}{\bigl(-\ln(- x)\bigr)^{\alpha_-}}
=\left\{ \begin{array}{cl}
           0 & \text{if} \quad 0 < \alpha_- < 1, \\
           1 & \text{if} \quad \alpha_- = 1, \\
           +\infty& \text{if} \quad \alpha_- > 1.
         \end{array} \right.
\end{equation*}
Hence Theorem~\ref{lem:KacKr} yields
\begin{equation*} 
0\in \sigma_{\!\operatorname{ess}}(B_{-,0})\cap \sigma_{\!\operatorname{ess}}(B_{+,0}) \quad \Leftrightarrow  \quad  \alpha_+> 1  \quad
\text{and} \quad \alpha_- > 1.
\end{equation*}
Since $0$ is not an eigenvalue of $A$, it follows from the preceding equivalence that $0 \in c(A)$ if and only if $\alpha_+> 1$ and $\alpha_- > 1$. Theorem~\ref{thm:Reg0}~\!(\ref{thm:Reg0i2}) yields that $0 \in c_r(A)$, whenever $\alpha_+> 1$ and $\alpha_- > 1$. Conversely, if $\alpha_+ \in (0,1]$ or $\alpha_- \in (0,1]$, then $0 \notin c(A)$. Consequently, $0 \in c_r(A)$ if and only if $\alpha_+> 1$ and $\alpha_- > 1$.

Next assume that $b_+=1$ and $b_- \in (-1,0)$.
Due to Theorem~\ref{thm:Reg0}(\ref{thm:Reg0i3})  $0\not\in c_s(A)$. In this case $\ker(A)=\{0\}$,
since a function $f\in\ker(A)$
should have a form   $f= f_+\oplus f_-\in\dom A$, where $f_+\in\dom(B_+^{\langle*\rangle})$, $f_-\in\dom(B_-^{\langle*\rangle})$ and
satisfy the conditions
\begin{equation}\label{eq:coupl5}
  f_+(0)=f_-(0),\quad f'_+(0)=f'_-(0),\quad f'_-(b_-)=0.
\end{equation}
The conditions $f_+\in\dom(B_+^{\langle*\rangle})$ and $f'_-(b_-)=0$ yield that $f_+$ is proportional to $1-x$ and $f_-$ is constant. Then~\eqref{eq:coupl5} implies $ f_+ = 0$ and hence, also $f_- = 0$.

Since the spectrum of $B_{-,0}$ is discrete, $0$ is not an accumulation point of the negative spectrum of $A$ and consequently $0 \notin c(A)$. The same conclusion holds if
$b_-=-1$ and $b_+ \in (0,1)$.

Finally we assume that $b_- \in (-1,0)$ and $b_+ \in (0,1)$. In this case the differential expression $\mathfrak{a}$ is regular, so the spectrum of $A$ is discrete. Therefore the root space at $0$ is nondegenerate. Consequently, $0 \notin c_s(A)$.
Since
\[
W_-(b_-)+W_+(b_+)=\frac{1}{(-\ln b_+)^{\alpha_+}}-\frac{1}{(-\ln |b_-|)^{\alpha_-}},
\]
statement (\ref{thm:Reg0i6}) from Theorem~\ref{thm:Reg0} takes the form:
\begin{equation*}
 \frac{\alpha_+}{\alpha_-} \neq \frac{\ln\bigl|\ln |b_-|\bigr|}{\ln|\ln b_+|} \quad \Leftrightarrow \quad \ker A = \ker A^2 \quad \Leftrightarrow \quad 0 \notin c(A).
\end{equation*}
It is interesting to write the preceding equivalences in the following form:
\begin{equation*}
 \frac{\alpha_+}{\alpha_-} = \frac{\ln\bigl|\ln |b_-|\bigr|}{\ln|\ln b_+|} \quad \Leftrightarrow \quad \ker A \subsetneq \ker A^2 \quad \Leftrightarrow \quad 0 \in c_r(A).
\end{equation*}
\end{example}

\subsection{Similarity} 

The coupling operator $A$ in the Krein space $\cK$ is nonnegative and $\rho(A) \ne \emptyset$. Hence, it has at most two critical points $0$ and $\infty$. Thus, $A$ is similar to a self-adjoint operator in a Hilbert space if and only if its critical points are regular and $\ker A=\ker A^2$, see Theorem \ref{WellKnown}.
Combining Theorems~\ref{thm:SL}, \ref{thm:Reg0}, \ref{thm:Reg_infty_PIS}, and~\ref{thm:Reg_infty_PI}
 we obtain the following list of sufficient conditions for similarity of $A$ to a self-adjoint operator in a Hilbert space, which equals Property (Si) from the introduction.

\begin{theorem}\label{thm:PolF}
 Let   $A$ be the differential operator associated with the expression $\mathfrak{a}$ with the domain defined by~\eqref{eq:domA} and let
 $W_\pm$ and $R_\pm$ be defined by~\eqref{eq:WR_pm}. Let at least one of the conditions {\rm(\ref{thm:Reg0i1})-(\ref{thm:Reg0i4})} in Theorem~{\rm\ref{thm:Reg0}} be in force.
Then the following statements hold.
\begin{enumerate}
\renewcommand*\theenumi{\alph{enumi}}
\renewcommand*\labelenumi{\rm{(\theenumi)}}
  \item \label{thm:PolFia}
If either $W_+\circ R_+^{-1}$ is positively increasing at $0_+$ or $W_-\circ R_-^{-1}$ is positively increasing at $0_-$, then the operator $A$ is similar to a self-adjoint operator in a Hilbert space.
  \item \label{thm:PolFib}
  Let $W_+\circ R_+^{-1}$ be slowly varying at $0_+$ and let $W_-\circ R_-^{-1}$ be slowly varying at $0_-$. Then the operator $A$ is similar to a self-adjoint operator in a Hilbert space if and only if
\begin{equation*} 
\left( 1 + \frac{W_-\bigl(R_-^{-1}(-x)\bigr)}{W_+\bigl(R_+^{-1}(x)\bigr)} \right)^{\!\!-1} = O(1) \quad \text{as}\quad  x\downarrow 0.
\end{equation*}
\end{enumerate}
\end{theorem}

\begin{example}
Let us consider Example~\ref{ex:Reg1} on an interval $I=(b_-,b_+)$ with $-1 \leq b_- < 0 < b_+ \leq 1$.
Combining the conclusions made in Example~\ref{ex:Reg1} and Example~\ref{ex:Reg1_0} we obtain the following equivalence:\\
The operator $A$ is similar to a self-adjoint operator in a Hilbert space if and only if \begin{enumerate}
  \item either  $\max\{b_+,|b_-|\}=1$ and $\frac{\alpha_+}{\alpha_-}\ne 1$;
  \item or $\max\{b_+,|b_-|\}<1$ and $\frac{\alpha_+}{\alpha_-}\not\in\left\{1, \frac{\ln|\ln |b_-||}{\ln|\ln b_+|}\right\}$.
\end{enumerate}

\end{example}

\appendix

\section{Some results from Karamata's theory} \label{appen}

In Appendix we present the definitions and the results from Karamata's theory of regularly varying functions that we use in the paper. Standard references for Karamata's theory are \cite{BGT} and~\cite{Sen}. For completeness we include a few standard results from Karamata's theory and some of them are reformulated to fit our needs. In addition, we present  Theorem~\ref{th-acf} and Corollary~\ref{cor:AtkCon}
that seem to be new.

\subsection{Definitions and basic results} First we give definitions of regularly varying functions.

\begin{definition} \label{def:sv}
Let $a,\alpha \in \mathbb{R}$ with $a > 0$.  A measurable function $f:(0,a]\to\mathbb{R}_+$ is called {\em regularly varying at $0$ from the right with index}  $\alpha$ if the following condition is satisfied:
\[
\text{for all} \quad \lambda \in  \mathbb{R}_+  \qquad \text{we have} \qquad \lim_{x\downarrow 0} \frac{f(\lambda x)}{f(x)} = \lambda^{\alpha}.
\]
When $\alpha = 0$ the function $f$ is called {\em slowly varying at $0$ from the right}.

A measurable function $g:[a,+\infty) \to\mathbb{R}_+$ is called {\em regularly varying at $+\infty$ with index} $\alpha$ if the following condition is satisfied:
\[
\text{for all} \quad \lambda \in  \mathbb{R}_+  \qquad \text{we have} \qquad \lim_{x\to +\infty} \frac{g(\lambda x)}{g(x)} = \lambda^{\alpha}.
\]
When $\alpha = 0$ the function $g$ is called {\em slowly varying at $+\infty$}.

A measurable function $g:[-a,0) \to\mathbb{R}_-$ is called {\em regularly varying at $0$ from the left with index} $\alpha$ if the function $f(x) = -g(-x)$ where $x \in (0,a]$ is regularly varying at $0$ from the right with index $\alpha$. When $\alpha = 0$ the function $g$ is called {\em slowly varying at $0$ from the left}.
\end{definition}

We will often use ``at $0_+$'' as an abbreviation for the phrase ``at $0$ from the right'' and ``at $0_-$'' as an abbreviation for the phrase ``at $0$ from the left.''

The Karamata's theory of regular variation is commonly presented for functions regularly varying at $+\infty$. The results for functions regularly varying at $0_+$ follow from the following equivalence. Let $f$ and $g$ be measurable functions such that $g(x) = f(1/x)$ for all $x$ in the domain of $g$ for which $1/x$ is in the domain of $f$. Then $g$ is regularly varying at $+\infty$ with index $\alpha$ if and only if $f$ is regularly varying at $0_+$ with index $-\alpha$.

In this section some results will be presented at $0_+$ and some at $+\infty$. This choice is sometimes made based on our needs in this paper and sometimes on convenience.

Slow variation plays the central role in the theory of regular variation. That centrality is expressed in the following proposition that follows immediately from the definition.

\begin{proposition} \label{pro:rvsv}
Let $a, \alpha \in \mathbb{R}$ with $a > 0$ and let $f, g:(0,a]\to\mathbb{R}_+$ be measurable functions such that $g(x) = x^\alpha f(x)$ for all $x\in(0,a]$. The function $g$ is regularly varying at $0_+$ with index $\alpha$ if and only if $f$ is slowly varying at $0_+$.
\end{proposition}

The next theorem is Karamata's Representation Theorem, see \cite[Theorem~1.3.1]{BGT} or \cite{Kar} for Karamata's original paper.

\begin{theorem} \label{th:SVRT}
Let $a \in \mathbb{R}$. A function $f:[a,+\infty)\to\mathbb{R}_+$ is slowly varying at $+\infty$ if and only if there exist $b \in [a,+\infty)$, a measurable function $m : [b,+\infty) \to \mathbb{R}_+$ and a continuous function $\varepsilon : [b,+\infty) \to \mathbb{R}$ such that
\[
\lim_{x\to+\infty} m(x) = M \in \mathbb{R}_+, \qquad \lim_{x\to+\infty} \varepsilon(x) = 0,
\]
and for all $x \geq b$ we have
\begin{equation*} 
f(x) = m(x) \exp\left(\int_{b}^{x} \frac{\varepsilon(t)}{t} dt \right).
\end{equation*}
\end{theorem}

The following property of regularly varying functions follows from Proposition~\ref{pro:rvsv} and  Theorem~\ref{th:SVRT}, see \cite[1\degree~on page~18]{Sen}.

\begin{corollary} \label{co-rvl}
If $g$ is a regularly varying function at $+\infty$ with a positive {\rm (}negative, respectively{\rm )} index, then
\[
\lim_{x\to+\infty} g(x) = +\infty \qquad \text{\rm (} \lim_{x\to+\infty} g(x) = 0, \
\text{respectively}\text{\rm )}.
\]
If $f$ is a regularly varying function at $0_+$ with a positive {\rm (}negative, respectively{\rm )} index, then
\[
\lim_{x\downarrow 0} f(x) = 0 \qquad \text{\rm (} \lim_{x\downarrow 0} f(x) = +\infty, \ \text{respectively} \text{\rm )}.
\]
\end{corollary}

\subsection{Karamata's characterization and consequences}

The following theorem is our restatement of Karamata's characterization of regular variation as it appears in \cite[Theorem~1.2.1]{dH}, \cite[Theorems~IV.5.2 and~IV.5.3]{Kor}, \cite[Theorems~1.5.11 and~1.6.1]{BGT} and \cite{BoPa}. In \cite{BGT, BoPa, dH, Kor} regular variation at $+\infty$ is considered. Here we characterize regular variation at $0_+$.

\begin{theorem} \label{th-Kon}
Let $a\in \mathbb{R}_+$ and let $f:(0,a] \to \mathbb{R}_+$ be a locally integrable function on $(0,a]$. Let $\alpha, \gamma \in \mathbb{R}$ be such that $\gamma + \alpha \neq 0$ and consider the following two  conditions:
\begin{gather} \label{eq-KorCon1}
\int_{0}^{a} s^{\gamma-1} f(s) ds \quad \text{exists as an improper integral at $0$},
\\
\label{eq-KorCon2}
\lim_{v \downarrow 0} \frac{1}{v^\gamma f(v)} \int_{0}^{v} s^{\gamma-1} f(s) ds
= \frac{1}{\gamma + \alpha}.
\end{gather}
The following statements are equivalent:
\begin{enumerate}
\renewcommand*\theenumi{\alph{enumi}}
\renewcommand*\labelenumi{{\rm (\theenumi)}}
\item
$f$ is regularly varying at $0_+$ with index $\alpha$.
\item \label{th-Kon-2}
For all $\gamma \in \mathbb{R}$ such that $\gamma + \alpha > 0$ conditions \eqref{eq-KorCon1} and \eqref{eq-KorCon2} hold.
\item \label{th-Kon-3}
There exists $\gamma \in \mathbb{R}$ such that $\gamma + \alpha > 0$ and 
\eqref{eq-KorCon1} and \eqref{eq-KorCon2} hold.
\end{enumerate}
\end{theorem}

The next theorem is a reformulation of the preceding one in terms of the differential of the function under consideration.

\begin{theorem} \label{th-acf}
Let $a, \alpha, \gamma \in \mathbb{R}$ be such that $a > 0$, $\gamma \neq 0$ and $\gamma + \alpha \neq 0$. Let $f:(0,a] \to \mathbb{R}_+$ be a measurable function which is of bounded variation on each closed interval contained in $(0,a]$. Consider the following three conditions:
\begin{gather} \label{eq-Conn1}
\int_{0}^{a} s^\gamma df(s) \quad \text{exists as an improper Riemann-Stieltjes integral at $0$},
\\
\label{eq-Conn2}
\lim_{v \downarrow 0} v^\gamma f(v) = 0,
\\
\label{eq-Conn3}
\lim_{v \downarrow 0} \frac{1}{v^\gamma f(v)} %
 \displaystyle\int_{0}^{v}\! s^\gamma df(s) = \frac{\alpha}{\gamma + \alpha}.
\end{gather}
The following statements are equivalent:
\begin{enumerate}
\renewcommand*\theenumi{\roman{enumi}}
\renewcommand*\labelenumi{{\rm (\theenumi)}}
\item \label{th-acf-1}
$f$ is regularly varying at $0_+$ with index $\alpha$.
\item \label{th-acf-2}
For all $\gamma \in \mathbb{R}\!\setminus\!\{0\}$ such that $\gamma + \alpha > 0$ conditions \eqref{eq-Conn1}, \eqref{eq-Conn2} and \eqref{eq-Conn3} hold.
\item \label{th-acf-3}
There exists $\gamma \in \mathbb{R}\!\setminus\!\{0\}$ such that $\gamma + \alpha > 0$ and conditions \eqref{eq-Conn1}, \eqref{eq-Conn2} and \eqref{eq-Conn3} hold.
\end{enumerate}
\end{theorem}

\begin{proof}
Let $u, v \in (0,a]$ such that $u < v$. First notice that since $f$ is of bounded variation on $[u,v]$, see \cite[Theorems~2.21 and~2.24]{WZ}, the integration by parts yields
\begin{equation} \label{eq-int-by-par}
\int_{u}^{v} s^\gamma df(s)
= v^{\gamma} f(v) - u^{\gamma} f(u) - \gamma \int_{u}^{v} s^{\gamma -1} f(s) ds.
\end{equation}

Assume (\ref{th-acf-1}). Let $\gamma \in \mathbb{R}\!\setminus\!\{0\}$ be such that $\gamma + \alpha > 0$. Since by Definition~\ref{def:sv} the function $x\mapsto x^\gamma f(x)$ is regularly varying at $0_+$ with index $\gamma+\alpha$, Corollary~\ref{co-rvl} yields \eqref{eq-Conn2}.

Theorem~\ref{th-Kon} implies that \eqref{eq-KorCon1} and \eqref{eq-KorCon2} hold. Letting $u\downarrow 0$ in \eqref{eq-int-by-par} and using \eqref{eq-KorCon1} yields \eqref{eq-Conn1} and
\begin{equation} \label{eq-Conn-p1}
\frac{1}{v^{\gamma} f(v)} \int_{0}^{v} s^\gamma df(s) = 1 - \frac{\gamma}{v^{\gamma} f(v)} \int_{0}^{v} s^{\gamma -1} f(s) ds.
\end{equation}
Now letting $v\downarrow 0$ and using \eqref{eq-KorCon2} we deduce \eqref{eq-Conn3}, proving  (\ref{th-acf-2}).

The fact that (\ref{th-acf-2}) implies (\ref{th-acf-3}) is trivial. Now assume (\ref{th-acf-3}). Letting $u\downarrow 0$ in \eqref{eq-int-by-par} and using \eqref{eq-Conn1} yields \eqref{eq-KorCon1}, and we  again deduce \eqref{eq-Conn-p1}. Together \eqref{eq-Conn-p1} and \eqref{eq-Conn3} imply \eqref{eq-KorCon2} in Theorem~\ref{th-Kon}.  Thus,  (\ref{th-Kon-3}) in Theorem~\ref{th-Kon} holds and   (\ref{th-acf-1}) follows from Theorem~\ref{th-Kon}.
\end{proof}

Let $\gamma > 0$. With the substitution $t = v^\gamma$, conditions \eqref{eq-Conn1}, \eqref{eq-Conn2} and \eqref{eq-Conn3} can be rewritten as (see \cite[Theorem~12.11]{PM} for the change of variables formula in Riemann-Stieltjes integral)
\begin{gather*}
\int_{0}^{a^\gamma} t\mkern+1.5mu df(t^{1/\gamma}) \quad
\text{exists as an improper Riemann-Stieltjes integral at $0$}, \\
\lim_{t \downarrow 0} t f(t^{1/\gamma}) = 0, \\
 \lim_{t \downarrow 0} \frac{1}{t f(t^{1/\gamma})}
 \displaystyle\int_{0}^{t} s\mkern+1.5mu df\bigl(s^{1/\gamma}\bigr) = \frac{\alpha/\gamma}{1 + \alpha/\gamma}.
\end{gather*}
This observation and Theorem~\ref{th-acf} (with $\gamma$ being $1$ and $\alpha$ being $\alpha/\gamma$)  yield the following equivalence: The function $t\mapsto f(t^{1/\gamma})$ with $t\in(0,a^\gamma]$ is regularly varying at $0_+$ with index $\alpha/\gamma > -1$ if and only if conditions
\eqref{eq-Conn1}, \eqref{eq-Conn2}, \eqref{eq-Conn3} hold. Here it is convenient to read the last fraction in \eqref{eq-Conn3} as $(\alpha/\gamma)/\bigl(1+(\alpha/\gamma)\bigr)$.

The next corollary generalizes the preceding equivalence to any increasing bijection on $[0,a]$.

\begin{corollary} \label{cor:AtkCon}
Let $\alpha, a, b \in \mathbb{R}$ be such that $a, b > 0$ and $\alpha > -1$. Let $f:(0,b] \to \mathbb{R}_+$ be a function of bounded variation on every closed subinterval of $(0,b]$ and let $g:[0,b] \to [0,a]$ be an increasing bijection. The function $f\mkern-3mu\comp\mkern-3mu g^{-1}: (0,a] \to \mathbb{R}_+$ is regularly varying at $0_+$ with index $\alpha > -1$ if and only if the following three conditions are satisfied:
\begin{gather}\label{eq-WRn21}
\int_{0}^{b}\! g(s)\mkern+1.5mu df(s) \quad \text{exists as an improper Riemann-Stieltjes integral at $0$}, \\
\label{eq-WRn22}
\lim_{v \downarrow 0} f(v) g(v) = 0, \\
\label{eq-WRn23}
\lim_{v \downarrow 0} \frac{1}{f(v)g(v)}\! \int_{0}^{v}\! g(s)\mkern+1.5mu df(s) = \frac{\alpha}{1+\alpha}.
\end{gather}
\end{corollary}

\begin{proof}
Let $u, v \in (0,b]$ such that $u < v$. As in the preceding theorem we notice that since $f$ is of bounded variation on $[u,v]$ the integration by parts
(\cite[Theorem~2.21]{WZ})  yields
\begin{equation} \label{eq-int-by-par2}
\int_{u}^{v} g(s)\mkern+1.5mu df(s) = f(v)g(v) - f(u)g(u) - \int_{u}^{v} f(s)\mkern+1.5mu dg(s).
\end{equation}
In this proof we will also use that, since $g$ is a continuous increasing bijection, we have that $u \downarrow 0$ if and only if $g(u)\downarrow 0$.

Assume \eqref{eq-WRn21}, \eqref{eq-WRn22} and \eqref{eq-WRn23}. Letting $u\downarrow 0$ and using \eqref{eq-WRn21} and  \eqref{eq-WRn22} in \eqref{eq-int-by-par2} yields
\begin{equation} \label{eq-ribp1}
\int_{0}^{v} g(s)\mkern+1.5mu df(s) = f(v)g(v) - \int_{0}^{v} f(s)\mkern+1.5mu dg(s)
\end{equation}
for all $v \in (0,b]$. Therefore, for all $v \in (0,b]$ we have
\begin{equation} \label{eq-ribp}
\begin{split}
\frac{1}{f(v)g(v)} \int_{0}^{v} g(s)\mkern+1.5mu df(s)
& = 1 - \frac{1}{f(v)g(v)}\int_{0}^{v} f(s) dg(s) \\
   & = 1 - \frac{1}{f(v)g(v)}\int_{0}^{g(v)} f\bigl(g^{-1}(t)\bigr) dt,
\end{split}
\end{equation}
where, for the second equality, we used the change of variables formula in Riemann-Stieltjes integral, \cite[Theorem~12.11]{PM}.  Now \eqref{eq-WRn23} implies
\[
\frac{1}{1+\alpha} = \lim_{v \downarrow 0} \frac{1}{f(v)g(v)}\int_{0}^{g(v)} f\bigl(g^{-1}(t)\bigr) dt  = \lim_{u \downarrow 0} \frac{1}{u f\bigl(g^{-1}(u)\bigr)}\int_{0}^{u} f\bigl(g^{-1}(t)\bigr) dt.
\]
Since we assume $1+\alpha > 0$, Theorem~\ref{th-Kon} yields that $f\mkern-3mu\comp\mkern-3mu g^{-1}$ is regularly varying at $0_+$ with index $\alpha$.

To prove the converse assume that $f\mkern-3mu\comp\mkern-3mu g^{-1}$ is regularly varying at $0_+$ with index $\alpha > -1$. Then the function $x\mapsto xf\bigl(g^{-1}(x)\bigr)$ is regularly varying at $0_+$ with index $\alpha + 1 > 0$ and \eqref{eq-WRn22} follows from Corollary~\ref{co-rvl} after a change of variables in the limit.  By the change of variables formula for all $u \in (0,a]$ we have
\[
\int_{u}^{a} s\mkern+1.5mu d f\bigl(g^{-1}(s)\bigr) = \int_{g^{-1}(u)}^{b} g(t)\mkern+1.5mu df(t).
\]
Consequently, \eqref{eq-WRn21} follows from \eqref{eq-Conn1} in Theorem~\ref{th-acf} applied to $f\mkern-3mu\comp\mkern-3mu g^{-1}$ with $\gamma = 1$. Therefore, \eqref{eq-ribp1} and consequently \eqref{eq-ribp} both hold. Now \eqref{eq-WRn23} follows from \eqref{eq-KorCon2} in Theorem~\ref{th-Kon}  with $\gamma=1$.
\end{proof}

\subsection{Asymptotic equivalence of functions on a sequence} \label{subsec:AsEq}

In the next definition we extend the notation $\sim$ of asymptotic equivalence of functions to hold only on a  sequence.

\begin{definition}\label{Bichota}
Let $a \in \mathbb{R}_+$. For functions $f,g:[a,+\infty) \to \mathbb{R}_+$ we write
\[
f \ssim g \quad \text{at} \quad +\mkern -3mu\infty
\]
if and only if there exists an increasing sequence $(x_n)$ in $[a,+\infty)$ such that
\begin{equation*} 
\lim_{n\to+\infty} x_n  = +\infty \quad \text{and} \quad \lim_{n\to+\infty} \frac{f(x_n)}{g(x_n)} = 1.
\end{equation*}
For functions $f,g:(0,a] \to \mathbb{R}_+$ we write
\[
f \ssim g \  \text{at} \  0_+
\]
if and only if there exists a decreasing sequence $(x_n)$ in  $(0,a]$ such that
\[
\lim_{n\to+\infty} x_n  = 0 \quad \text{and} \quad \lim_{n\to+\infty} \frac{f(x_n)}{g(x_n)} = 1.
\]
\end{definition}

Recall, see \cite{Bel}, \cite[5.10.11]{TBB}, that for a function $\phi:[a,+\infty)$ a real number $L$ is a {\em cluster value} of $\phi$ at $+\infty$ if for every $\epsilon > 0$ and for every $X \in \mathbb{R}$ there exists $x > X$ such that $|\phi(x) - L| < \epsilon$. Similarly, for a function $\phi:(0, a]$ a real number $L$ is a {\em cluster value} of $\phi$ at $0_+$ if for every $\epsilon > 0$ and for every $\delta > 0$ there exists $x \in (0,\delta)$ such that $|\phi(x) - L| < \epsilon$.  Notice that $f \ssim g$ at $+\infty$ (at $0_+$) if and only if $1$ is a cluster value of the function $f/g$ at $+\infty$ (at $0_+$).

\begin{proposition}
Let $f$ and $g$ be regularly varying functions at $+\infty$ with indices $\alpha$ and $\beta$, respectively. If $f \ssim g$ at $+\infty$, then $\alpha = \beta$.
\end{proposition}
\begin{proof}
We will prove the contrapositive. Assume that $\alpha < \beta$. Since the function  $f(x)/g(x)$ is regularly varying with index $\alpha-\beta<0$ it follows from Corollary~\ref{co-rvl} that $\lim_{x\to+\infty}f(x)/g(x)=0$. Thus, $f \ssim g$ at $+\infty$ is not true. If $\alpha > \beta$ the preceding limit is $+\infty$, so $f \ssim g$ at $+\infty$ is not true in this case either.
\end{proof}

The converse of the preceding proposition is not true. For example, let $f$ be a slowly varying function at $+\infty$ and $g=2f$. Then $\alpha = \beta =0$, but $f \ssim g$ at $+\infty$ is clearly not true.

The following theorem extends \cite[Proposition~0.8(vi)]{Res} to the concept introduced in the previous definition. This theorem can be deduced from \cite[Corollary~7.66]{BIKS}. A direct proof is presented in~\cite[Theorem A.11, Corollary A.12]{CDT21}.

\begin{theorem} \label{cor:SVinv0}
Let $f$ and $g$ be strictly monotonic positive functions defined in a neighbourhood of $0_+$ and let $f$ be regularly varying at $0_+$ with a nonzero index.
\begin{enumerate}
\renewcommand*\theenumi{\alph{enumi}}
\renewcommand*\labelenumi{{\rm (\theenumi)}}
\item
If $f$ and $g$ are increasing with $0$ limit at $0_+$, then the inverses $f^{-1}$ and $g^{-1}$ are also increasing, defined in a neighbourhood of $0_+$ and the following equivalence holds
\[
f \ssim g \ \text{at} \ 0_+ \quad  \Leftrightarrow \quad f^{-1} \ssim g^{-1} \ \text{at} \ 0_+.
\]
\item
If $f$ and $g$ are decreasing and unbounded, then the inverses $f^{-1}$ and $g^{-1}$ are  decreasing, defined in a neighbourhood of $+\infty$ and the following equivalence holds
\[
f \ssim g \ \text{at} \ 0_+ \quad  \Leftrightarrow \quad f^{-1} \ssim g^{-1} \ \text{at} \ +\infty.
\]
\end{enumerate}
\end{theorem}

The following corollary is a consequence of the fact that the negation of $f \ssim g$  at  $0_+$ is the statement
\[
\left(\frac{f(x)}{g(x)} - 1 \right)^{-1} = O(1) \quad \text{as} \quad x \downarrow 0.
\]
Each of the two  statements in Theorem~\ref{cor:SVinv0} can be expressed using one of these negations. We state only the analogue of the last statement in Theorem~\ref{cor:SVinv0} since that is what is used in Theorem~\ref{thm:Reg_infty_PI}.

\begin{corollary} \label{cor:RVinv}
Let $f$ and $g$ be strictly monotonic positive functions defined in a neighbourhood of $0_+$ and let $f$ be regularly varying at $0_+$ with a nonzero index.  If $f$ and $g$ are decreasing and unbounded, then the inverses $f^{-1}$ and $g^{-1}$ are  decreasing, defined in a neighbourhood of $+\infty$ and the following equivalence holds
\[
\left(\frac{f(x)}{g(x)} - 1 \right)^{\!\!-1} = O(1) \ \text{as} \ x \downarrow 0 \quad  \Leftrightarrow \quad  \left(\frac{f^{-1}(y)}{g^{-1}(y)} - 1 \right)^{\!\!-1} = O(1) \ \text{as} \ y \to +\infty.
\]
\end{corollary}

Clearly $f\sim g$ at $+\infty$ implies $f\ssim g$ at $+\infty$. In the next example we will demonstrate that $f\ssim g$ at $+\infty$ does not imply $f\sim g$ at $+\infty$ even for smooth normalized slowly varying increasing functions $f$ and $g$ for which $f/g$ is normalized slowly varying function.

\subsection{Positively increasing functions}
The following class of functions was introduced as a generalization of regularly varying functions with positive index, see \cite[Section~3.1 and Definition~3.26]{BIKS}.

\begin{definition}\label{CieloDiablo}
Let $a \in \mathbb{R}_+$. A nondecreasing function $f: (0,a] \to \mathbb{R}_+$ is called {\it positively increasing at} $0$ {\it from the right} if there exists $\lambda \in (0,1)$ such that
\[
\limsup_{x\downarrow 0}\frac{f(\lambda x)}{f(x)} < 1.
\]
A function $g:[-a,0) \to\mathbb{R}_-$ is called {\it positively increasing at} $0$ {\it from the left} if the function $f(x) = -g(-x), x \in [-a,0)$, is positively increasing at $0$ from the right.

A function $g:[a,+\infty) \to\mathbb{R}_+$ is called {\it positively increasing at} $+\infty$ if the function $f(x) = 1/g(1/x), x \in (0,1/a]$, is positively increasing at $0$ from the right. A function $g:(-\infty,-a] \to\mathbb{R}_-$ is called {\it positively increasing at} $-\infty$ if the function $f(x) = -1/g(-1/x), x \in (0,1/a]$, is positively increasing at $0$ from the right.
\end{definition}

The relationship between regularly varying and positively increasing functions at $+\infty$, and analogously at $-\infty$, $0_+$ and $0_-$, is as follows. Each regularly varying function with positive index is positively increasing, while a regularly varying function with a nonpositive index is not positively increasing. In particular, a slowly varying function is not positively increasing. The exponential function $\exp$ is positively increasing at $+\infty$ but not regularly varying at $+\infty$. As was shown in~\cite[Example A.17]{CDT21} there exists a nondecreasing function  $f:[1,+\infty)\to\mathbb{R}_+$ which is neither positively increasing nor slowly varying at $+\infty$.

\section*{Acknowledgment}

The authors thank Aleksey Kostenko for fruitful discussions and literature hints. Volodymyr Derkach gratefully acknowledges financial  support by
the Ministry of Education and Science of Ukraine (project \# 0121U109525),  by the German Research Foundation (DFG), grant TR 903/22-1, by the Fulbright Program,
and the travel support by the David and Darla Kennerud Visiting Math Scholars Fund at Western Washington University.

The authors thank the anonymous reviewer for very careful reading of our paper, for a number of valuable
comments and several literature suggestions.


\begin{thebibliography}{99}
\bibitem{AP}
N.L.~Abasheeva, S.G.~Pyatkov,  Counterexamples in indefinite Sturm-Liouville problems. Siberian Advances in Mathematics. Siberian Adv. Math. 7 (1997), no.~4, 1-8.

\bibitem{AG}
N.I.~Akhiezer, I.M.~Glazman, \textit{Theory of Linear Operators in
Hilbert Space}, two volumes bound as one, Dover Publications, 1993.
\bibitem{Ako80}
R.V. Akopjan, On the regularity at infinity of the spectral function of a $J$-nonnegative operator, Izv. Akad. Nauk
Arm. SSR Ser. Mat. 15 (5) (1980) 357-364 (in Russian)

\bibitem{Atk85}
F.V.~Atkinson, Some further estimates for the Titchmarsh-Weyl $m$-coefficient, preprint, University of Toronto, 1985.


\bibitem{AI}
T.Ya.~Azizov,  I.S.~Iokhvidov, \textit{Linear Operators in Spaces
with an Indefinite Metric}, John Wiley \& Sons, 1990.


\bibitem{Bea85}
R.~Beals,  Indefinite Sturm-Liouville problems and half range completeness,
J.\ Differential Equations 56 (1985) 391–-407.

\bibitem{Bel}
C.L.~Belna, Cluster sets of arbitrary real functions: a partial survey. Real Anal. Exchange 1 (1976) 7--20.

\bibitem{Ben87}
C.~Bennewitz, The HELP inequality in the regular case, General
inequalities, 5 (Oberwolfach, 1986), 337--346, Internat.\ Schriftenreihe Numer.\ Math., 80, Birkhäuser, Basel, 1987.

\bibitem{Ben89}
C.~Bennewitz, Spectral asymptotics for Sturm-Liouville equations,
Proc.\ London Math.\ Soc.\ 59 (1989) 294--338.

\bibitem{BC06}
P.~Binding, B.~\'{C}urgus, Riesz bases of root vectors of indefinite Sturm-Liouville problems with eigenparameter dependent boundary conditions.~I. Operator theory and indefinite inner product spaces, 75--95, Oper. Theory Adv. Appl., 163, Birkhäuser, Basel, 2006.

\bibitem{BC09}
P.~Binding, B.~\'{C}urgus,
Riesz bases of root vectors of indefinite Sturm-Liouville problems with eigenparameter dependent boundary conditions.~II. Integral Equations Operator Theory 63 (2009), no. 4, 473--499.

\bibitem{BGT}
N.H.~Bingham; C.M.~Goldie; J.L.~Teugels,  \textit{Regular variation.} Encyclopedia of Mathematics and its Applications, 27. Cambridge University Press, Cambridge, 1989.

\bibitem{B74}
        J.~Bogn\'ar, \textit{Indefinite Inner Product Space},
        Springer-Verlag, Berlin, 1974.

\bibitem{BoPa}
V.V.~Buldygin; V.V.~Pavlenkov, Karamata theorem for regularly log-periodic functions, Ukrainian Math.\ J.\ 64 (2013) 1635--1657.

\bibitem{BIKS}
V.V.~Buldygin; K-H.~Indlekofer; O.I.~Klesov; J.G.~ Steinebach, Pseudo-regularly varying functions and generalized renewal processes. Probability Theory and Stochastic Modelling, 91. Springer, Cham, 2018.

\bibitem{ChisEve71}
  R.S. Chisholm and W.N. Everitt. On bounded integral operators in the space of integrable-square functions. Proc. Roy. Soc. Edinb. (A), 69 (1970/71), 199-204.


\bibitem{Cur}
B.~\'{C}urgus, On the regularity of the critical point infinity of
definitizable operators, Integral Equations Operator Theory 8 (1985)
462--488.

\bibitem{CD15}
B.~\'Curgus, V.~Derkach, Partially fundamentally reducible
operators in Krein spaces, Integral Equations
    Operator Theory 82 (2015) 469-518.

\bibitem{CDT21}
B.~\'Curgus, V.~Derkach, C~Trunk, Indefinite Sturm-Liouville operators in polar form. (2021)   arXiv:2101.00104v2 [math.SP]

\bibitem{CFK13}
B.~\'Curgus, A.~Fleige, A.~Kostenko, The Riesz basis property of an
indefinite Sturm-Liouville problem with non-separated boundary
conditions, Integral Equations Operator Theory 77 (2013) 533--557.

\bibitem{CL89}
B.~\'Curgus, H.~Langer,  A Krein space approach to symmetric
ordinary differential operators with an indefinite weight function,
J.\  Differential Equations 79 (1989) 31--61.

\bibitem{CN94}
B.~\'Curgus, B.~Najman,  A Krein space approach to elliptic eigenvalue problems with indefinite weights. Differential 88Integral Equations 7 (1994) 1241--1252.

\bibitem{CN95}
B.~\'Curgus, B.~Najman, The operator $(\sgn x){d^{2}}/{dx^{2}}$ is
similar to a self-adjoint operator in $L^{2}({\mathbb R})$,  Proc.\
Amer.\  Math.\  Soc.\ 123 (1995) 1125--1128.

\bibitem{CN96}
B.~\'Curgus, B.~Najman, Positive differential operators in Kreĭn space $L^2(\mathbb{R})$. Recent developments in operator theory and its applications (Winnipeg, MB, 1994), 95--104, Oper. Theory Adv. Appl., 87, Birkhäuser, Basel, 1996.

\bibitem{CN98}
B.~\'Curgus, B.~Najman, Positive differential operators in the Krein space $L^2(\mathbb{R}^n)$. Contributions to operator theory in spaces with an indefinite metric (Vienna, 1995), 113--129, Oper. Theory Adv. Appl., 106, Birkhäuser, Basel, 1998.

\bibitem{CRead}
B.~\'Curgus, T. Read,  Discreteness of the spectrum of second-order differential operators and associated embedding theorems. J.\ Differential Equations 184 (2002)  526--548.
\bibitem{DL77}
K.~Daho, H.~Langer,
Sturm--Liouville operators with an indefinite weight function. Proc.\ R.\ Soc.\ Edinb.\ Sect.\ A 87, 161--191 (1977)

\bibitem{D95}
  V.A.~Derkach,
    On Weyl function and generalized resolvents of a Hermitian
    operator in a Krein space,  Integral Equations Operator Theory
      23 (1995) 387--415.

\bibitem{D99}
V.A.~Derkach.
On generalized resolvents of Hermitian relations in Krein spaces,
 J.\ Math.Sci., 97 (1999), 
 4420--4460.

\bibitem{DHMS00}
V.A.~Derkach, S.~Hassi, M.M.~Malamud, H.S.V.~de Snoo, Generalized resolvents of symmetric operators and admissibility, Methods Funct. Anal. Topology 6 (2000), no.~3, 24--55.

\bibitem{DM91}
V.A.~Derkach, M.M.~Malamud, Generalized resolvents and the boundary
value problems for Hermitian operators with gaps, J.\ Funct.\ Anal.\ 95
(1991) 1--95.
\bibitem{DM17}
        V.\ Derkach, M.\ Malamud, \textit{Extension Theory of Symmetric Operators and Boundary Value Problems}, Transactions of Institute of Mathematics NAS of Ukraine, 104, Kyiv, 2017.

\bibitem{DSW20}
V.~Derkach, D.~Strelnikov, H.~Winkler, On a class of integral systems, arXiv:2010.01295,  2020

\bibitem{DT17}
V.~Derkach, C.~Trunk, Coupling of defnitizable operators in Krein spaces, Nanosystems:
Physics, Chemistry, Mathematics 8 (2017), 166--179.

\bibitem{EE87}
D.~Edmunds, W.~Evans, \textit{Spectral theory and differential operators.} Oxford:
Oxford University Press, 1987.

\bibitem{EvEv91}
W.D.~Evans, W.N.~Everitt, HELP inequalities for limit-circle and regular problems,
Proc.\ Roy.\  Soc.\  Lond.\  A 432 (1991), 367--390.

\bibitem{Ev72}
W. N. Everitt, On an extension to an integro-differential inequality of Hardy, Littlewood and Polya, Proc. Roy. Soc. Edinburgh A 69 (1972), 295--333.

\bibitem{FS00}
M.M. Faddeev, R.G. Shterenberg, On the similarity of some singular
differential operators to self-adjoint operators. (Russian) Zap.
Nauchn. Sem. S.-Peterburg. Otdel. Mat. Inst. Steklov. (POMI) 270
(2000), Issled. po Linein. Oper. i Teor. Funkts. 28, 336--349,
370--371; translation in J. Math. Sci. (N. Y.) 115 (2003) 2279--2286
\bibitem{Fl95}
A.~Fleige,  The turning point condition of Beals for indefinite Sturm-Liouville problems,
Math.\ Nachr.\ 172 (1995) 109--112.

\bibitem{Fl98}
A.~Fleige,  A counterexample to completeness properties for indefinite Sturm-Liouville
problems, Math.\ Nachr.\  190 (1998) 123--128

\bibitem{Fl08}
A.~Fleige, The Riesz basis property of an indefinite Sturm-Liouville problem with a non odd
weight function, Integral Equations Operator Theory 60 (2008) 237--246.

\bibitem{Fl15}
A.~Fleige, The Critical Point Infinity Associated with Indefinite Sturm-Liouville Problems. In: Alpay D.\ (eds) Operator Theory. Springer, Basel, 2015.

\bibitem{FN98}
A.~Fleige, B.~Najman, Nosingularity of critical points of some
differential and difference operators, Oper. Theory: Adv. Appl.,
vol. 102, Birkh\"auser, Basel, 1998.

\bibitem{Glaz50}
I.M.~Glazman, On the theory of singular differential operators.
 Uspekhi Matematicheskikh Nauk, 5(6):102--135, 1950.

\bibitem{GG}
V.I.~Gorbachuk,  M.L.~Gorbachuk, \textit{Boundary Value Problems for
Operator Differential Equations.} Kluwer Academic Publishers Group,  Dordrecht,
1991.

\bibitem{dH}
L.~de Haan, \textit{On regular variation and its application to the weak convergence of sample extremes.} Mathematical Centre Tracts, 32 Mathematisch Centrum, Amsterdam 1970.


\bibitem{HT}
A.~Horning, A.~Townsend, FEAST for differential eigenvalue problems. SIAM J. Numer. Anal. 58 (2020) 1239--1262.

\bibitem{J82}
P.~Jonas, Compact perturbations of definitizable operators II, J.~Operator Theory 8 (1982), 3--18.

\bibitem{J84}
P.~Jonas, Regularity criteria for critical points of definitizable operators. 179--195, Operator Theory: Advances and Applications  14,  Birkh\"{a}user, 1984.

\bibitem{JonLan79}
        P.\ Jonas, H.\ Langer,
        Compact perturbations of definitizable operators, J.\ Operator Theory 2 (1979), 63--77.

\bibitem{KacKr58}
I.S.~Kac, M.G.~Krein, Criteria for the discreteness of the spectrum of a singular string, (Russian) Izv.\ Vysš.\ Učebn.\ Zaved.\ Matematika 2  (1958) 136--153.


\bibitem{Kac73}
I.S.~Kac, A generalization of the asymptotic formula of V. A. Mar\v{c}enko for the spectral functions of a second order boundary value problem. (Russian) Izv. Akad. Nauk SSSR Ser. Mat. 37 (1973), 422–436; English transl. in: Math USSR Izv. 7 (1973), 424--436.

\bibitem{KaKr74}
I.S.~Kac, M.G.~Krein,~$R$-functions--analytic functions mapping
the upper halfplane into itself, Amer.\ Math.\ Soc.\ Transl.\ Ser.\ (2)
103 (1974) 1--18.

\bibitem{KaKr74B}
I.S.~Kac, M.G.~Krein, On the spectral functions of the string, Supplement II to the Russian edition of F.V.\ Atkinson, Discrete and continuous boundary problems, Mir, Moscow, 1968 (Russian) (English translation: Amer.\  Math.\ Soc.\  Transl.\ (2) 103 (1974), 19--102).

\bibitem{Kal74}
H.~Kalf, Remarks on some Dirichlet type results for semibounded Sturm-Liouville operators, Math.\ Ann.\  210 (1974) 197--205.


\bibitem{Kar00}
I.M.~Karabash, $J$-self-adjoint ordinary differential operators
similar to self-adjoint operators. Methods Funct.\ Anal.\ Topology (2000)  22--49.

\bibitem{Kar10}
I.M.~Karabash, A functional model, eigenvalues, and finite singular critical points for indefinite Sturm-Liouville operators, Oper. Theory: Adv. Appl 203, (2009), 247--287 (arXiv:0902.4900)

\bibitem{KaKost08}
I.M.~Karabash, A.~Kostenko,  Indefinite Sturm-Liouville operators
with the singular critical point zero, Proc.\ Roy.\ Soc.\ Edinburgh
Sect.\ A 138 (2008) 801--820.

\bibitem{KaKoMal09}
I.M.~Karabash, A.~Kostenko, M.M.~Malamud, The similarity problem for
$J$-nonnegative Sturm-Liouville operators, J.\ Differential Equations
246 (2009) 964--997.

\bibitem{KaMal07}
I.M.~Karabash, M.M.~Malamud, Indefinite Sturm-Liouville operators
with finite zone potentials, Operators and Matrices 1 (2007)
301--368.

\bibitem{KT09}  I.~Karabash, C.~Trunk, Spectral properties of singular
Sturm-Liouville operators, Proc.\ Roy.\ Soc.\ Edinburgh Sect.\ A, 139
(2009) 483--503.

\bibitem{Kar}
J.~Karamata, Sur un mode croissance r\'{e}guli\`{e}re des fonctions. Mathematica (Cluj) 4 (1930) 38--53.

\bibitem{Kas75}
Y.~Kasahara, \emph{Spectral theory of generalized second order differential operators and its applications to Markov processes}, Japan. J. Math. (N.S.) \textbf{1} (1975/76), no.1, 67--84.

\bibitem{Koch79}
A.N.~Kochubei,  Extensions of $J$-symmetric operators. (Russian)
Teor.\ Funkci\u{\i} Funkcional.\ Anal.\ i Prilozhen.\ 31 (1979) 74--80.

\bibitem{Kor}
J.~Korevaar, \textit{Tauberian Theory: A Century of Developments.} Grundlehren der Mathematischen Wissenschaften, Springer, 2004.

\bibitem{Kost11}
A.~Kostenko, The similarity problem for indefinite Sturm-Liouville operators with periodic coefficients. Oper. Matrices 5 (2011), 707--722.

\bibitem{Kost13}
A.~Kostenko, The similarity problem for indefinite Sturm-Liouville operators and the HELP inequality, Adv.\  Math.\  246 (2013) 368--413.

\bibitem{Kost14}
A.~Kostenko, On a necessary aspect for the Riesz basis property for indefinite Sturm-Liouville problems, Math. Nachr. 287, no. 14--15 (2014), 1710--1732 (arXiv:1202.2444)

\bibitem{La}
H.~Langer, Spectral functions of definitizable operators in
Krein spaces. Functional analysis (Dubrovnik, 1981) 1--46,
Lecture Notes in Math.\ 948, Springer, 1982.

\bibitem{Mar52}
V.A.~Mar\v{c}enko,  Some questions of the theory of one-dimensional linear differential operators of the second order. I. (Russian) Trudy Moskov. Mat. Ob\v{s}\v{c}. 1, (1952), 327--420.
\bibitem{Maz85}
V.G. Maz'ja, Sobolev spaces, Springer-Verlag, Berlin Heidelberg New York, 1985.

\bibitem{Muck72}
B. Muckenhoupt,  Weighted norm inequalities for the Hardy maximal function. Trans.
Am. Math. Soc. 165 (1972) 207--226.

\bibitem{OinOtel99}
R. Oinarov and M. Otelbaev, A criterion for a general Sturm-Liouville operator to have a
discrete spectrum, Differential Equations 24 (1988), 402--408.

\bibitem{Par03}
A.I.~Parfenov,  On an embedding criterion for interpolation spaces and application to indefinite spectral problems, Sib.\ Math.\ J.\ 44 (2003) 638--644

\bibitem{Par05}
A.I.~Parfenov, The \'Curgus condition in indefinite Sturm-Liouville problems, Sib.\ Adv.\ Math.\ 15 (2005) 68--103.

\bibitem{PM}
M.H.~Protter, C.B.~Morrey, {\it A first Course in Real Analysis}. Second edition. Undergraduate Texts in Mathematics. Springer-Verlag, New York, 1991.

\bibitem{Py1}
S.G.~Pyatkov, Interpolation of some function spaces and indefinite Sturm-Liouville problems. Differential and integral operators (Regensburg, 1995), 179--200, Oper. Theory Adv. Appl., 102, Birkhäuser, Basel, 1998.


\bibitem{Py3}
Pyatkov, S.G. Indefinite elliptic spectral problems. Sib. Math. J. 39 (1998) 358--372.

\bibitem{RS}
M.~Reed, B.~Simon: {\it Methods of Modern Mathematical Physics,
IV.~Analysis of Operators}, Academic Press, New York 1978.

\bibitem{RemSc}
C.~Remling, K.~Scarbrough, The essential spectrum of canonical systems, J.\ Approx.\ Theory  254 (2020) 105395, 11 pp.

\bibitem{Res}
S.I.~Resnick, \textit{Extreme values, regular variation, and point processes.} Applied Probability. A Series of the Applied Probability Trust, 4. Springer-Verlag, New York, 1987.

\bibitem{RomWor}
R.~Romanov, H.~Woracek,  Canonical systems with discrete spectrum, J.\ Funct.\ Anal.\ 278 (2020) 108318, 44 pp.

\bibitem{Shm74}
Ju.L.~Shmuljan,  Operator extension theory and spaces with indefinite metric. (Russian) Izv.\ Akad.\ Nauk SSSR Ser.\ Mat.\ 38 (1974) 896--908.

\bibitem{Sen}
E.~Seneta, \textit{Regularly Varying Functions}, Lecture Notes in Mathematics, Vol.\ 508. Springer-Verlag, Berlin-New York, 1976.

\bibitem{Simon2015}
B.~Simon, Operator theory. Compr. Course Anal., Part 4.
American Mathematical Society, Providence, RI, 2015.

\bibitem{St32}
M.H.~Stone,
\textit{Linear Transformations in Hilbert Space and their Applications
  to Analysis}, AMS Colloquium Publ., American Mathematical Soc., New York, 1979.

\bibitem{Stu72}
C.A. Stuart. ’The measure of non-compactness of some linear integral operators.’ Proc. Roy. Soc. Edinb. (A), 71 (1972/73), 167--179.
\bibitem{TBB}
B.S.~Thomson, J.B.~Bruckner, A.M.~Bruckner,
\textit{Elementary Real Analysis}, 2nd ed., www.classicalrealanalysis.com, 2008.

\bibitem{Tit62}
E.C.~Titchmarsh, \textit{ Eigenfunction Expansions Associated with Second-Order Differential Equations. Part I}, 2nd ed., Clarendon Press, Oxford, 1962.

\bibitem{Ves}
K.~Veseli\'c,  On spectral properties of a class of $J$-self-adjoint
operators. I, II. Glasnik Mat.\ Ser.\ III 7(27) (1972) 229--248; ibid.
7(27) (1972) 249--254.

\bibitem{Vol}
H.~Volkmer, Sturm–Liouville problems with indefinite weights and Everitt’s
inequality, Proc.\ Roy.\ Soc.\ Edinburgh Sect.\ A 126 (1996) 1097--1112.

\bibitem{Weyl12}
H.~Weyl, Das asymptotische Verteilungsgesetz der Eigenwerte linearer partieller Differentialgleichungen (mit einer Anwendung auf die Theorie der Hohlraumstrahlung), Math. Ann. 71 (1912), 441--479.

\bibitem{WZ}
R.L.~Wheeden, A.~Zygmund, \textit{Measure and integral. An Introduction to Real Analysis}, 2nd ed., Pure and Applied Mathematics. CRC Press, Boca Raton, 2015.

\end{thebibliography}
\end{document}